\documentclass{article}

\usepackage{graphicx}
\usepackage{amsmath,amssymb,amsfonts}
\usepackage{amsthm}
\usepackage{mathrsfs}
\usepackage{xcolor}
\usepackage{textcomp}
\usepackage{manyfoot}
\usepackage{booktabs}
\usepackage{psfig}
\usepackage{epsfig}

\newtheorem{theorem}{Theorem}

\newtheorem{definition}{Definition}
\usepackage{mathptmx}  
\usepackage{lipsum}
\usepackage{subcaption}
\usepackage{epstopdf}
\usepackage{algorithmic}
\usepackage{siunitx} 
\usepackage{color}
\newtheorem{lemma}[theorem]{Lemma}
\title{Dynamics of Performances in a Competitive Model}
\author{
  {Mahmoud Soufbaf}\footnote{School of Mathematics, Statistics and Computer Science, College of Science, University of Tehran, Tehran, Iran. soufbaf.mahmoud@ut.ac.ir},
  {Gholam Reza Rokni Lamouki}\footnote{School of Mathematics, Statistics and Computer Science, College of Science, University of Tehran, Tehran, Iran. rokni@ut.ac.ir (Corresponding Author)},
  {Khosro Tajbakhsh}\footnote{Faculty of Mathematical Sciences, Tarbiat Modares University, Tehran, Iran. khtajbakhsh@modares.ac.ir}
}
\date{}
\begin{document}

\maketitle

{\bf Abstract:} A competitive resource-consumer dynamical model is analyzed based on an integrated model of a competitive Lotka-Volterra model and a prey-predator Rosenzweig-MacArthur model that we call that LV-RM model throughout this paper. Resource growth in the absence of consumers is logistic, and competing consumers' type II Holling's functional response made the model structure more realistic. We used the normal form and the center manifold theorems for bifurcation analysis of the presented model, identified Hopf and zero-Hopf bifurcations and their directions, and discussed their biological interpretations. We hypothesized that differentiated time scales of the competing consumers' predatory that lead to asymmetry in competition are the mechanisms that promote coexistence through relaxation-oscillation dynamics. Though, other performance parameters of both competitors are the same. Graphical representation of variations of the first Lyapunov coefficient, after competition coefficients interplay, shows various dynamics with growing complexity from the periodic state towards chaotic motion like R\"ossler attractor. We presented simulations to visualize the theoretical results obtained through bifurcation analysis.\\


\section{Introduction}\label{sec1}

Competition among organisms for limited common resources is controversial among researchers, who apply many hypotheses and tests~{
  \cite{connor86} and \cite{hubbell2011}.} Generally, competition between species takes place through two main paths. The first is the exploitative/exploratory competition which involves indirect negative interactions. Here each consumer affects the other by reducing the abundance of the resource. The second is interference, which involves direct negative interactions between species due to mechanisms such as territoriality, predation, and chemical competition; see~\cite{amare2002}, {
  \cite{caseGilpin} }, \cite{shoner83}, and \cite{vance84}.
The classical competition theory predicts that competing species can only coexist through the mechanisms such as resource allocation in time and space, differences in the emergence and attack rates, as well as distribution and search efficiency~\cite{amare2007}\cite{vance84}. However, the microbial competition models studied in a chemostat show that competition eliminates the population that needs a higher concentration of nutrients to grow~\cite{smoth96}. It is common to believe that the diversity of food resources can explain the high biodiversity in animals~\cite{tilman}; however, the enriched biodiversity of plant species, which require the same physical and chemical resources, is unclear~\cite{grubb}. Researchers believe that the persistent disruption has created mosaics of localized resources. Then, due to spatial arrangement and partitioning of competitors on respective accessible resources, competitive exclusion avoids occur~\cite{vandenbosh}. The importance of spatial processes in population dynamics is a concern to ecologists. It is now well established that the qualitative behavior of interactions between species, such as host-parasite, prey-predator, and competitors in a heterogeneous environment, is quite different from the corresponding behavior in a homogeneous one~\cite{vandenbosh}. For instance,~\cite{pekas} suggested that parasitoids may be able to coexist on the same host species if they partition host resources according to size, age, and stage; or if their dynamics vary at spatial and temporal scales.
{
Some researchers proposed other mechanisms of coexistence of competing consumers from different taxa. For instance, Ayala (1971)  showed that the stable equilibrium between Drosophila willistoni and Drosophila pseudoobscura requires their relative fitnesses in competition to be frequency-dependent. They showed that the relative fitnesses of these two species are inversely related to their relative frequencies after about four generations. In competition, Drosophila pseudoobscura oscillates with a frequency 35\%  of the frequency of Drosophila willistoni. Livingston et al. (2012) showed that two Pseudomonas bacterial strains can stably coexist in a landscape if they trade off between competitive and colonization abilities. Also, local competition between relatives for resources reduces kin selection for altruism in various species of many taxa such as pathogenic bacteria and social insects,~\cite{grifin}, \cite{west}. 
}

However, reducing niche overlap mainly by weaker competitors is suggested as a priory mechanism promoting the spatio-temporal coexistence of both species in a community \cite{mills}. There are many ways to incorporate non-uniformity into models of the consumer population. The time delay in the consumer response to the change in resource density and the relatively slow resource population dynamics with certain types of functional responses cause the limit cycle{
  ; however, \cite{huisman2001} regards competition as a process that may generate oscillations and chaotic fluctuations.} 
Generally, the biological model uses the linear functional response for consumers~\cite{fernandez}. But, in many communities, such as insect parasitoids, the governing feeding mechanism is of type II Holling functional response. This kind of growth can better realize the population dynamics of interacting species. The shape of type II functional response is an asymptotic curve that decelerates prey/host density increases due to the handling time. In this form, the asymptotic level reflects the maximum attack rate of the predator/parasitoid in consuming prey/host. This functional response, which means parasitoids have limited fixed time for searching hosts, leads to inverse density-dependent parasitism~\cite{fernandez}. Therefore, applying this type of functional response is compatible with competition interactions among parasitoids for a limited density of the hosts.\\
{
Even though population size and dynamics are crucial fitness criteria of an organism, visualizing temporal dynamics of these correlates versus depicting dynamics of other fitness characteristics such as feeding rate, reproduction rate, and predation ability need to be synchronized due to considerable differences between times scales of each one. 
In fact, if we accept correlations between different fitness parameters of an organism, one would expect to show roughly the same dynamic patterns of all fitness correlations. Then, we can visualize the population dynamics as any other fitness dynamics but with different time scales. 
For instance, \cite{collier2001} used two parasitoids~\footnote{Eretmocerus eremicus and Encarsia sophia} 
of the whitefly pest, Bemisia tabaci, as competitors. They suggested that host feeding on parasitized hosts may be a generally important mechanism of interaction in parasitoid communities where host-feeding parasitoids have disproportionately established and led to successful control relative to non-host-feeding species. 
They claimed that lethal interference competition in insect parasitoids, including facultative hyperparasitism, multiparasitism, and intraguild predation, lead to the same consequences for individuals, such that one kills another. 
The population consequences of these mechanisms of lethal interference competition can also be quite similar. However, the outcomes of this interaction are not entirely clear~\cite{collier2001}.
  }\\

\begin{figure}[h]
	\centering
	\includegraphics[width=0.5\linewidth]{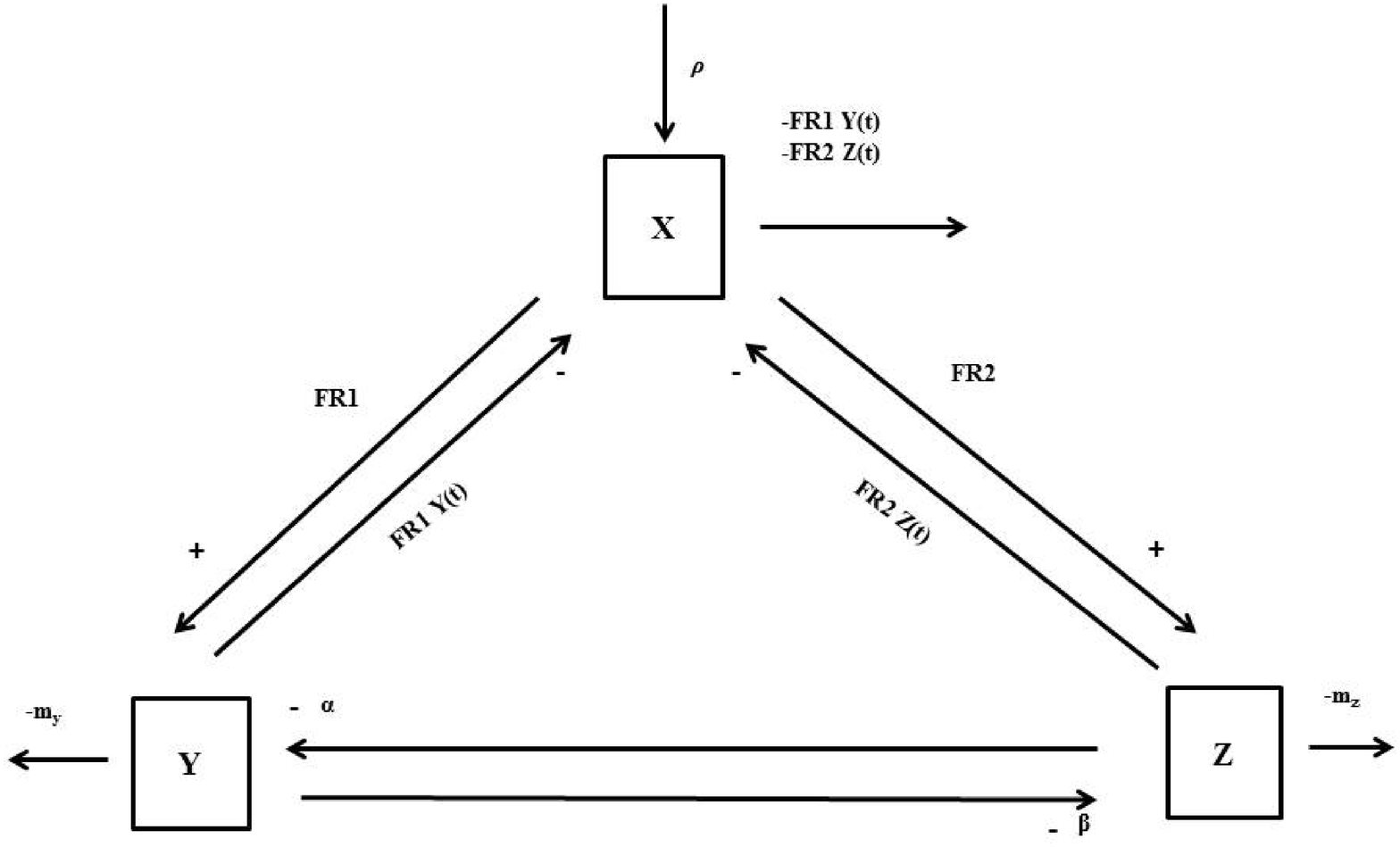}
	\caption{Flowchart diagram of the competitive system based on LV-RM model.}
	\label{fig1}
\end{figure}

In this study, we show that difference in predatory time scales is a mechanism that governs the coexistence of two similar competing consumers. One could consider this mechanism as a form of reducing nich overlap, as suggested by \cite{mills}. We have two time scales of competition between two predators. The first one, denoted by $\tau_1$, is the time scale of direct interaction and resource consumption. This time scale can range from hours to weeks. The related time scale depends on the predators' feeding habits. The second one, denoted by $\tau_2$, is the time scale of reproduction. We assume that $\tau_1<<\tau_2$. This inequality means that the population size of each predator may not change significantly during the resource consumption phase, but it may be affected by the cumulative outcome of competition when the reproduction phase occurs. Therefore, we have two coupled dynamics: one for resource consumption and one for reproduction. The resource in our model is the prey. In another view, predator and prey compete with each other to succeed in predation and escape from predation, respectively. Our biological system consists of three species with one prey and two predators. The two predators are competing with each other following the Lotka-Volterra (LV) traditional two-dimensional competitive model; each one acts over the prey based on the general prey-predator Rosenzweig-MacArthur (RM) model. The presented model here, model \eqref{eq1}, is integrated of two prey-predator models of RM class with the same prey, and we call it LV-RM model in this paper. The flow chart in Figure~\ref{fig1} illustrates the interactions of these three species. Based on this flow chart and employing the functional response of Holling type II, the equations of the governing system introduced in~\cite{soufbaf} restates as follows in equation~\eqref{eq1}.
\begin{equation}\label{eq1}
	\begin{split}
		&\frac{dX(t)}{dt}=X(t)(\rho(1-\frac{X(t)}{K})-\frac{\gamma_{y} Y(t)}{1+\gamma_{y} \mu_{y} X(t)}-\frac{\gamma_{z} Z(t)}{1+\gamma_{z} \mu_{z} X(t)})\\
		&\frac{dY(t)}{dt}=Y(t)(\frac{\gamma_{y} X(t)}{1+\gamma_{y} \mu_{y} X(t)}-m_y-\alpha_{yz} Z(t))\\
		&\frac{dZ(t)}{dt}=Z(t)(\frac{\gamma_{z} X(t)}{1+\gamma_{z} \mu_{z} X(t)}-m_z-\alpha_{zy} Y(t))
	\end{split}
\end{equation} 

Here, we assume that the highly fluctuating time series of interacting species in either prey-predator or competition models are not population sizes of the species but are the per capita successes of interacting species in predation and escaping from predation for predators and prey, respectively. Therefore, in model \eqref{eq1}, the prey resistance to predation, $X$, is hypothesized to be averaged per capita resistances $\zeta^1_t(N_1)$ over the population manifold of the prey $\Omega_{N_1(t)}$ both are time dependent. Similarly, the predation success of two predators, $Y$ and $Z$, are hypothesized to be averaged per capita predation success $\zeta^2_t(N_2)$ and  $\zeta^3_t(N_3)$ over the population manifolds of the consumers $\Omega_{N_2(t)}$ and $\Omega_{N_3(t)}$, respectively. The governing functionals are as follows in equation~\eqref{eq5}.
Here, $\omega_{N_k}$ is the distribution differential form over the population manifold $\Omega_{N_k(t)}$ for $k=1,2,3$.
\begin{equation}\label{eq5}
	\begin{split}
		&X(t)=\frac{1}{N_1(t)} \times \int_{\Omega_{N_1(t)}} \zeta^1_t(N_1)\omega_{N_1}, \quad N_1(t)=\int_{\Omega_{N_1(t)}}\omega_{N_1},\\
		&Y(t)=\frac{1}{N_2(t)} \times \int_{\Omega_{N_2(t)}} \zeta^2_t(N_2)\omega_{N_2}, \quad N_2(t)=\int_{\Omega_{N_2(t)}}\omega_{N_2},\\
		&Z(t)=\frac{1}{N_3(t)} \times \int_{\Omega_{N_3(t)}} \zeta^3_t(N_3)\omega_{N_3}, \quad N_3(t)=\int_{\Omega_{N_3(t)}}\omega_{N_3}.
	\end{split}
\end{equation}

Meanwhile, we assume that dynamics of population sizes of $N_1(t)$, $N_2(t)$, and $N_3(t)$ in equation~\eqref{eq5} are governed by standard logistic growth equations as shown in equation~\eqref{eq12}.
\begin{equation}\label{eq12}
	\begin{split}
		&\dot{N_1}(t)=r_{N_1} \times N_1(t) \times(1-\frac{1}{K_{N_1} \times X(t)}),\\
		&\dot{N_2}(t)=r_{N_2} \times N_2(t) \times(1-\frac{1}{K_{N_2} \times Y(t)}),\\
		&\dot{N_3}(t)=r_{N_3} \times N_3(t) \times(1-\frac{1}{K_{N_3} \times Z(t)}),
	\end{split}
\end{equation}
where $K_{\iota}$ and $r_{\iota}$, for ${\iota}=N_1, N_2, N_3$, are carying capacity and population growth rate of the resource and consumers, respectively.    
To show differences in predatory time scales, we hypothesized that the second competing consumer waits throughout the handling time of the first and, with an enhanced searching efficiency, consumes the resource by the lowest possible interference. Accordingly, per capita, the functional response ratio $\mathrm{FRR}$ of two competing consumers is as follows.
\begin{equation}\label{eq0}
	\begin{split}
		& \mathrm{FRR}=\cfrac{\; \; \cfrac{c \times \Xi \times \Upsilon\; \; }{1+a \times \Xi}}{\cfrac{d \times \Xi \times \Pi}{1+b \times \Xi}}
	\end{split}
\end{equation}
\begin{align}
	& \Xi      = (\text{prey resistance to predation}),       && (a,b,c,d)=(\text{regulatory parameters})    \nonumber\\
	& \Upsilon = (\text{predation success of predator 1}),  && \Pi= (\text{predation success of predator 2}).\nonumber
\end{align}
 
When $\mathrm{FRR}<1$, the first consumer's attack rate accelerates, while for $\mathrm{FRR}>1$, the second consumer benefits from a higher attack rate than the first consumer. When $\mathrm{FRR}=1$, consumers compete with each other with symmetric attack rates. Completing the handling time, the first competitor ends the successful predatory activity. Due to absorbing food and providing enough energy, the first competitor is in a higher position than the second competitor. The second competitor is hungry and suffers from energy deficiency. So the symmetry of the competition is lost. The highness of the first competitor is in terms of the innate ability to attack and other efficiency components, not predation motivation. To return the competitive system to symmetry, we strengthen the second competitor. This strengthening continues to some extent until coexistence is assured. In other words, during the simulation process, the activity of equalizing the competitors' power is continuously carried out by manipulating the attack rate. Competitive exclusion and coexistence are important from applied and theoretical viewpoints. We examine these two issues in this paper by incorporating type II functional response and different predatory time scales for ecologically similar competing consumers. Both of these issues predict oscillations in the system such that the whole system works around limit cycles. We study the limit cycles and their features by identifying the Hopf and zero-Hopf bifurcation. We use normal form and the center manifold theory to study the underlying dynamics of the system in a lower dimension by performing
codimension 1 and 2  bifurcation~\cite{Guckenheimer} and~\cite{khalil}. Table~\ref{tab:table1} presents the details of the parameters involved in the system.\\

Taking $0$ of accumulative predation success of either $Y$ or $Z$ in \eqref{eq1} leads to a two-dimensional prey-predator model of RM type while declining $X$ to $0$ leads to the extinction of all three species, $X$, $Y$, and $Z$.
However, taking $0$ of accumulative predation success of $Y$ and $Z$ in \eqref{eq1} leads to a one-dimensional prey-predator model of the logistic type.
Note that $\{\zeta^i_t(N_i)\}_{i=1}^3$ are the spatio-temporal instantaneous successes in prey-predator dynamics for species $i=1,2,3$ and their averaged integration over their population manifolds are their temporal averaged successes. Zero value for $\zeta^i_t(N_i)$ at a time-individual point would not impose extinction. But, zero value $\iota$, for a given $\iota=X, Y, Z$, means zero value for full-scale time-individual points over the related spatio-temporal population manifold, which imposes extinction for that species.  

\begin{table}[h]
		\caption{Variables and parameters used in the model \eqref{eq1}: $\theta$ and $\chi$ represent time and number of individual units, respectively. Time scales are different among different parameters.}
		\label{tab:table1}
		\begin{tabular}{l|S|l}
			\hline
			\textbf{Variable and Parameters} & \textbf{Notion} & \textbf{Unit}\\
			\hline
			Number ofsuccesses in escaping the predation & $X$ & $\chi$ \\
			Number of successes in the predation & $Y$ & $\chi$ \\
			Number of successes in the predation & $Z$ & $\chi$ \\
			The intrinsic growth rate of resource & $\rho$ & $\theta^{-1}$ \\
			Carrying capacity of the environment & K & $\chi$ \\
			The death rate of consumer 1 or equivalently lost predation & $m_y$ &  $\theta^{-1}$ \\
			The death rate of consumer 2 or equivalently lost predation & $m_z$ & $\theta^{-1}$ \\
			The attack rate of consumer 1 & $\gamma_{y}$ & $\chi^{-1} \theta^{-1}$\\
			The attack rate of consumer 2 & $\gamma_{z}$ & $\chi^{-1} \theta^{-1}$\\
			The handling time of consumer 1 & $\mu_{y}$ &$\theta$ \\
			The handling time of consumer 2 & $\mu_{z}$ &$\theta$ \\
			Competition effect from consumer 2 on consumer 1 & $\alpha_{yz}$ &$\chi^{-1} \theta^{-1}$ \\
			Competition effect from consumer 1 on consumer 2 & $\alpha_{zy}$ &$\chi^{-1} \theta^{-1}$\\
			\hline
		\end{tabular}
\end{table}
The current paper presents the dynamics of the system~\eqref{eq1}. We follow the procedure of non-dimensionalization by considering the time reparameterization and change of variables in equations~\eqref{eq2}.
\begin{equation} \label{eq2}
		\hat{X}=\frac{1}{K}X,\quad\hat{Y}=\frac{\gamma_{y}}{\rho}Y,\quad\hat{Z}=\frac{\gamma_{z}}{\rho}Z,\quad \hat{t}=t\rho.
\end{equation}
We consider the following reparameterization.
\begin{equation}\label{eq4}
		a=\gamma_{y} \mu_{y} K,\quad b=\gamma_{z} \mu_{z} K,\quad c=\frac{K \gamma_{y}}{\rho},\quad d=\frac{K \gamma_{z}}{\rho},
		\mu=\frac{m_y}{\rho},\quad \nu=\frac{m_z}{\rho},\quad  \alpha=\frac{\alpha_{yz}}{\gamma_{z}},\quad \beta=\frac{\alpha_{zy}}{\gamma_{y}}
\end{equation}
We let $\varphi=\{a,b,c,d,\mu,\nu,\alpha,\beta\}$. We rewrite the equations~\eqref{eq1} in terms of $(\hat{X}, \hat{Y}, \hat{Z}; \hat{t})$ and $\varphi$; then, revert $(\hat{X},\hat{Y},\hat{Z};\hat{t})$ to $(X,Y,Z,;t)$ for simplicity in notations. Thus, the equations \eqref{eq1} can be rewritten as \eqref{eq6}. 
\begin{equation}\label{eq6}
	\begin{split}
		&\frac{dX}{dt}=X(t)(1-X(t))-\frac{X(t)Y(t)}{1+aX(t)}-\frac{X(t)Z(t)}{1+bX(t)}\\
		&\frac{dY}{dt}=\frac{cX(t)Y(t)}{1+aX(t)}-\mu Y(t)-\alpha Y(t)Z(t)\\
		&\frac{dZ}{dt}=\frac{dX(t)Z(t)}{1+bX(t)}-\nu Z(t)-\beta Z(t)Y(t)
	\end{split}
\end{equation}
The equilibria of the system~\eqref{eq6} are a zero equilibrium (0,0,0) and a non zero positive equilibrium $(X^*,Y^*,Z^*)$ where
\begin{equation}\label{eq7}
	\begin{split}
		X^*(\varphi)&=-\frac{-b \upsilon +\xi +\upsilon }{3 b \upsilon }+\frac{\sqrt[3]{\sqrt{4 P^3+Q^2}+Q}}{3 \sqrt[3]{2} b \upsilon }-\frac{\sqrt[3]{2} P}{3 b \upsilon  \sqrt[3]{\sqrt{4 P^3+Q^2}+Q}}\\
		Y^*(\varphi)&=\frac{-a \alpha  d \left(X^*\right)^2+a \alpha  d X^*+a \mu  \nu -b \mu  \nu -c \nu +\alpha  d+d \mu -\alpha  d X^*}{-a \beta  \mu +b \beta  \mu +\beta  c+\alpha  d}\\
		Z^*(\varphi)&=\frac{-a \mu  \nu -b \beta  c \left(X^*\right)^2+b \beta  c X^*+b \mu  \nu +\beta  c+c \nu -\beta  c X^*-d \mu }{a \alpha  \nu -\alpha  b \nu +\beta  c+\alpha  d}
	\end{split}
\end{equation}
Here, $\upsilon=\upsilon(\varphi)=a \alpha  \beta$, $\xi=\xi(\varphi)=\alpha  b \beta$,
$P=P(\varphi)$
\footnote{$P(\varphi)=3 b \upsilon  (-a \beta  \mu +\alpha  \beta -\alpha  b \nu +\beta  c+\alpha  d-\xi -\upsilon )-(-b \upsilon +\xi +\upsilon )^2$.},
and
$Q=Q(\varphi)$
\footnote{$Q(\varphi) =-2 a \alpha  \beta  \upsilon ^2-\frac{2 \alpha  \beta  b^2 \xi  \upsilon }{a}+9 a \beta  b^2 \mu  \upsilon ^2+2 a b^2 \xi  \upsilon ^2-9 a \beta  b \mu  \upsilon ^2+3 a b \xi  \upsilon ^2-3 a \xi  \upsilon ^2-3 \alpha  \beta  b^2 \xi  \upsilon -9 \alpha  b^2 \nu  \xi  \upsilon + 9 \alpha  b^3 \nu  \upsilon ^2+18 \alpha  b^2 \nu  \upsilon ^2+18 \beta  b^2 \mu  \upsilon ^2-9 \beta  b^2 c \upsilon ^2-9 \alpha  b^2 d \upsilon ^2+3 b^2 \xi  \upsilon ^2+3 \alpha  \beta  b \xi  \upsilon +9 \beta  b c \xi  \upsilon +9 \beta  b c \upsilon ^2+9 \alpha  b d \xi  \upsilon +9 \alpha  b d \upsilon ^2+12 b \xi  \upsilon ^2+3 \xi  \upsilon ^2$.}
 are new parameters produced by combinations of elemnets in $\varphi$. 
The nonzero equilibrium exists when $X^*\geq 0$, $Y^*\geq 0$ and $Z^*\geq 0$. The positivity conditions of this equilibrium as are assured in some regions in $\varphi-$space of parameters. For example,
\begin{align}
  &{\text{Region}}\; \mathrm{I}:
  \begin{cases}
        a>1, \quad b\geq a,\quad 0<\nu \leq \frac{\alpha  (-d)-d \mu }{a \mu -b \mu -c},& \\[10pt]
	0<X^*(\varphi)\leq \frac{1}{2} \sqrt{\frac{\alpha  a^2 d+4 a^2 \mu  \nu -4 a b \mu  \nu -4 a c \nu +2 \alpha  a d+4 a d \mu +\alpha  d}{a^2 \alpha  d}}+\frac{a-1}{2 a}.
  \end{cases}\\
  &{\text{Region}}\; \mathrm{II}:
  \begin{cases}
        b>1,\quad a\geq b,\quad 0<\mu \leq \frac{\beta  c+c \nu }{a \nu -b \nu +d},& \\[10pt]
	0<X^*(\varphi)\leq \frac{1}{2} \sqrt{\frac{-4 a b \mu  \nu +\beta  b^2 c+4 b^2 \mu  \nu +2 \beta  b c+4 b c \nu -4 b d \mu +\beta  c}{b^2 \beta  c}}+\frac{b-1}{2 b}.
  \end{cases}
\end{align}
There are nine more regions in parameter space within which the equilibrium $(X^*, Y^*, Z^*)$ is positive.  Following our hypothesis, we introduce two exactly similar competing consumers to depict dynamical properties after competitive effects, i.e., the interplay of competition coefficients $\alpha$ and $\beta$. Other parameters are valued symmetrically as $a=b,c=d$, and $\mu = \nu$.

\section{Dynamical Analysis}
The model \eqref{eq6} satisfies the invariant requirement of the theoretical approach in ecology; namely, trajectories of initial value problems fulfill the existence and uniqueness condition and will never exit from the positive region.
The positive region in $(X,Y,Z)-$space and its boundary are defined by
\begin{align}
  &\Omega=\{(X,Y,Z)|\; X\geq 0, Y\geq 0, Z\geq 0 \},\quad   \partial\Omega=S_X\cup S_y\cup S_Z \nonumber\\
  & S_X=\{(0,Y,Z):\; Y,Z\geq 0\},\quad S_Y=\{(X,0,Z):\; X,Z\geq 0\},\quad S_Z=\{(X,Y,0):\; X,Y\geq 0\}.
\end{align}
The following theorem shows the consistency of the model.\\

\begin{theorem}\label{th:cns}
      The positive region $\Omega$ is invariant under the dynamics of the model \eqref{eq6}.
\end{theorem}
\proof
Suppose $n$ is the outward unit normal vector on $\partial\Omega$. Based on the partitioning $\partial\Omega=S_X\cup S_y\cup S_Z$, the outward unit normal vectors of $\partial \Omega$ are $n_X=-e_1=(-1,0,0)$ for $S_X$, $n_Y=-e_2=(0,-1,0)$ for $S_Y$, and $n_Z=-e_3=(0,0,-1)$ for $S_Z$.
The positive region $\Omega$ is invariant under the dynamics of equations~\eqref{eq6} if and only if the inner product $n\cdot (\dot{X},\dot{Y},\dot{Z})$ is always non-positive on $\partial\Omega$. 
This requirement is satisfied by the following equations.
\begin{align}
  n_X\cdot (\dot{X},\dot{Y},\dot{Z})|_{S_X}\leq 0,\quad
  n_Y\cdot (\dot{X},\dot{Y},\dot{Z})|_{S_Y}\leq 0,\quad
  n_Z\cdot (\dot{X},\dot{Y},\dot{Z})|_{S_Z}\leq 0.
\end{align}
In fact, subsets $S_X$, $S_Y$, and $S_Z$ are invariant. The biological interpretation of these invariant sets is that the existence of each species in the model~\eqref{eq6} is independent of the existence of the other two species.
\hfill $\square$\\

The stability of zero equilibrium of the system~\eqref{eq6} is analyzed using the Jacobian matrix, $J_0$.

\begin{equation}
	J_0=	\begin{matrix}
		\left(
		\begin{array}{ccc}
			1 & 0 & 0 \\
			0 & -\mu & 0 \\
			0 & 0 & -\nu \\
		\end{array}
		\right)
	\end{matrix}
\end{equation}

Therefore, the characteristic equation and its roots are as follows \eqref{eq52}
\begin{equation}\label{eq52}
	-\lambda ^3+\lambda ^2 (-\mu -\nu +1)+\lambda  (\mu  (-\nu )+\mu +\nu )+\mu  \nu,\quad \lambda_1 = 1,\lambda_2=-\mu,\lambda_3=-\nu		
\end{equation}

The zero equilibrium has two negative real roots and one positive root. Thus, the origin is an unstable saddle equilibrium with index 1. This result is plausible since the exclusion of both predators, $Y$ and $Z$, in ~\eqref{eq6} is related to negative roots of \eqref{eq52}, and the positive root is related to the prey, $X$. The stability of the nonzero positive equilibrium and its local bifurcations along the model~\eqref{eq6} requires some algebra. The expansion of equations~\eqref{eq6} up to the fourth order is presented in equation~\eqref{eq9} using the change of variables $x=X-X^*$, $y=Y-Y^*$, and $z=Z-Z^*$.
The new equations in terms of new variables $x, y, z$ are as follows. 

\begin{equation}\label{eq9}
	\begin{split}
		\dot{x}(t)=&b_{10} x(t)^4 y(t)+b_9 x(t)^3 y(t)+b_8 x(t)^2 y(t)+b_7 x(t) y(t)+b_{14} x(t)^4 z(t)+b_{13} x(t)^3 z(t)\\
		&+b_{12} x(t)^2 z(t)+ b_{11} x(t) z(t)+b_4 x(t)^4+b_3 x(t)^3+b_2 x(t)^2+b_1 x(t)+b_5 y(t)+b_6 z(t):=F^1_{\phi}(x,y,z)\\
		\dot{y}(t)=&c_9 x(t)^4 y(t)+c_8 x(t)^3 y(t)+c_7 x(t)^2 y(t)+c_6 x(t) y(t)+c_5 x(t)^4+c_4 x(t)^3+c_3 x(t)^2+c_2 x(t)\\
		&+c_{11} y(t) z(t)+c_1 y(t)+c_{10} z(t):=F^2_{\phi}(x,y,z)\\
		\dot{z}(t)=&d_{10} x(t)^4 z(t)+d_9 x(t)^3 z(t)+d_8 x(t)^2 z(t)+d_7 x(t) z(t)+d_5 x(t)^4+d_4 x(t)^3+ d_3 x(t)^2\\
		&+d_2 x(t)+d_{11} y(t) z(t)+d_6 y(t)+d_1 z(t):=F^3_{\phi}(x,y,z)
	\end{split}
\end{equation}
Here, $\phi=(\mathbf{b,c,d})$, and $\mathbf{b}=\mathbf{b}(\varphi)=\{b_i(\varphi)\}_{i=1}^{14}$, $\mathbf{c}=\mathbf{c}(\varphi)=\{c_i(\varphi)\}_{i=1}^{11}$, and $\mathbf{d}=\mathbf{d}(\varphi)=\{d_i(\varphi)\}_{i=1}^{11}$ are new parameters produced by complex combinations of elemnets in $\varphi$ through long formula. We do not state these formula for simplicity, but we used the exact relations in further calculations and numerical computations. See equation~\eqref{eq4} for details of $\varphi$.\\

Let $\delta=(b_1,b_5,b_6,c_1,c_2,c_{10},d_1,d_2,d_6)=\mathrm{P}\phi$ for the projection matrix $\mathrm{P}$.
The Jacobian matrix, $J_1$, and the characteristic polynomial, $P(\lambda)$, of the system~\eqref{eq9} around $(x,y,z)=(0,0,0)$, equally for the positive equilibrium $E^*=\left(X^*, Y^*, Z^* \right)$ of the original system \eqref{eq6}, are as follow.
\begin{align}
	&J_1(\delta)=	\begin{matrix}
		\left(
		\begin{array}{ccc}
			b_1 & b_5 & b_6 \\
			c_2 & c_1 & c_{10} \\
			d_2 & d_6 & d_1 \\
		\end{array}
		\right)
	\end{matrix}.\\
		&P(\lambda)=+\lambda ^3-\lambda ^2 \left(b_1+c_1+d_1\right)+\lambda   \left(b_1 c_1 -b_5 c_2+b_1 d_1 -b_6 d_2 +c_1 d_1-c_{10} d_6\right)\\
	  &-\left(b_1 c_1 d_1-b_5 c_2 d_1-b_6 c_1 d_2+b_5 c_{10} d_2+b_6 c_2 d_6-b_1 c_{10} d_6\right):=\lambda^3+\mathrm{a}(\delta)\lambda^2+\mathrm{b}(\delta)\lambda+\mathrm{c}(\delta).
\end{align}
Let $F_{\phi}=(F^1_{\phi}, F^2_{\phi}, F^3_{\phi})^T$, $\chi=(x,y,z)^T$, and $\hat{F}_{\phi}=F_{\phi}-J_1(\delta)\chi$. Then, we have
\begin{equation}
  \dot{\chi}=J_1(\delta)\chi+\hat{F}_{\phi}(\chi).
\end{equation}

\begin{lemma}\label{lem1}
  The jacobian matrix $J_1$ of the linearization of the system~\eqref{eq6} at the positive equilibrium $E^*$ has a pair of imaginary eigenvalues $\lambda_{1,2}=\pm i\omega_0$ and a real non-zero eigenvalue  $\lambda_3 =\mathrm{a}_0$ if and only if $\delta=\delta_0\in\Omega_0$ where

  \begin{align}\label{eq10}
    \Omega_0=\{\delta:\; \mathrm{a}(\delta)\neq 0,\; \mathrm{b}(\delta)>0,\; \mathrm{c}(\delta)-\mathrm{a}(\delta)\mathrm{b}(\delta)=0\}
	\end{align}
 
\end{lemma}
\proof
 We look for the region of parameters for which $P(\lambda)$ has a pair of complex conjugate eigenvalues $\lambda_{1,2}=\pm i\omega_0$ and a real eigenvalue $\lambda_3=\mathrm{a}_0$. This setting requires that
\begin{align}
  P(\lambda)=(\lambda+\mathrm{a})(\lambda^2+\omega^2)+(\mathrm{b}-\omega^2) \lambda +(\mathrm{c}-\mathrm{a} \omega^2).
\end{align}

By some simple algebra for $P(\lambda)=0$ we have $\lambda_3=\mathrm{a}_0$, and $\lambda_{1,2}=\pm i\omega_0=\pm i\sqrt{\mathrm{b}_0}$ when $\mathrm{a}_0=\mathrm{a}(\delta_0),\; \mathrm{b}_0=\mathrm{b}(\delta_0),\; \mathrm{c}_0=\mathrm{c}(\delta_0)$, and $\delta_0\in\Omega_0$.
\hfill $\square$\\

We can choose $\phi_0$ such that $\delta_0=\mathrm{P}\phi_0\in\Omega_0$.
Suppose $\delta_0\in\Omega_0$ and $J_1(0)=J_1(\delta_0)$. For the given $\epsilon-$ perturbation $\delta(\epsilon)=\delta_0+O(\epsilon)$, an $\epsilon-$perturbation of matrix $J_1(0)$ is the matrix $J_1^{\epsilon}=J_1\left(\delta_0+O(\epsilon)\right)$. Then $\mathrm{a}(\epsilon)=\mathrm{a}_0+O(\epsilon)$, $\mathrm{b}(\epsilon)=\mathrm{b}_0+O(\epsilon)$, and $\mathrm{c}(\epsilon)=\mathrm{c}_0+O(\epsilon)$. 

\begin{theorem}\label{th1}
  For small enough real number $\epsilon$, with $|\epsilon|<<1$ the jacobian matrix $J_1^{\epsilon}$ of the linearization of the system~\eqref{eq6} at the positive equilibrium $E^*$ has a pair of complex conjugate eigenvalues $\lambda_{1,2}=\mu(\epsilon)\pm i\omega(\epsilon)$ and a real non-zero eigenvalue  $\lambda_3(\epsilon) \in \mathbb{R}$. For $\epsilon=0$, we have $\lambda_3(0)=\mathrm{a}_0\neq 0$, $\mu(0)=0$, $\mu'(0)\neq 0$, and $\omega(0)=\omega_0>0$. 
\end{theorem}
\proof

We look for a perturbation $\delta(\epsilon)=\delta_0+O(\epsilon)$ such that $P(\lambda)$ has a pair of complex conjugate eigenvalues $\lambda_{1,2}=\mu(\epsilon)\pm i\omega(\epsilon)$ and a real eigenvalue $\lambda_3(\epsilon)$.
Based on the continuity of the roots of a polynomial with respect to its coefficients, this setting requires that
\begin{align}
  P(\lambda)=\left(\lambda+(\mathrm{a}+2\mu)\right)\left(\lambda^2-2\mu\lambda+(\mu^2+\omega^2)\right),
\end{align}
subject to the following conditions for $|\epsilon|<<1$.
	\begin{align}\label{eq10-2}
          & \mathrm{a}=\mathrm{a}_0+O(\epsilon),\; \mathrm{b}=\mathrm{b}_0+O(\epsilon),\; \mathrm{c}=\mathrm{c}_0+O(\epsilon),\nonumber\\
          &   \omega^2(\epsilon)=\mathrm{b}_0+O(\epsilon),\; \mu=O(\epsilon). \nonumber
	\end{align}
These conditions fulfill the requirements. Moreover, we have $\mu'(0)\neq 0$. \hfill $\square$\\

Based on Theorem~\ref{th1}, the positive equilibrium $E^*$ of the system~\eqref{eq6} undergoes a generic Hopf bifurcation at $\delta_0\in\Omega_0$ in equation (\ref{eq10}).
According to center manifold theory, there exists a two-dimensional center manifold of system~\eqref{eq9} at $E^*$ such that system~\eqref{eq9} on this manifold has a center-type equilibrium which is a weak focus at the bifurcation point $\mu=0$. That is a hyperbolic stable focus for $\mu<0$ and is a hyperbolic unstable focus for $\mu>0$. The change in the stability appears under small perturbation of the parameters of the system~\eqref{eq6} around any $\delta_0\in\Omega_0$. On one side, a limit cycle appears depending on the sign of the first Lyapunov coefficient. Then, the system~\eqref{eq9} experiences a Hopf bifurcation.. 
\hfill $\square$  \\

We have the possibility of Hopf bifurcation when $\lambda_3\neq 0$ and Zero-Hopf bifurcation when $\lambda_3=0$. First we consider the Hopf bifurcation which divides into two cases of $\lambda_3>0$ and $\lambda_3<0$.

\subsection{The case of Hopf bifurcation}

Suppose we choose $\phi_0$ such that $\delta_0=\mathrm{P}\phi_0\in\Omega_0$, and $\phi(\epsilon)=\phi_0+O(\epsilon), \delta(\epsilon)=\delta_0+O(\epsilon)$.
The center manifold theory reduces the stability analysis of the system around the origin to the stability analysis of a system whose order is equal to the number of eigenvalues of the Jacobian matrix with the zero real part~\cite{Carr}. After applying the Hopf conditions in the system~\eqref{eq9}, we construct the matrix $M_1$ by the eigenvectors related to the hyperbolic eigenvalue, the imaginary part, and the real part of complex eigenvalues.
We let $V=(u,v,w)^T$; then, the transformation $V=M_1^{-1}\chi$ changes the system variable to the standard form as follows.
\begin{align}
        \dot{V}=M_1^{-1} J_1(\delta) M_1 V + M_1^{-1}\hat{F}_{\phi}\left( M_1 V \right).
\end{align}
Here, $M_1=M_1(\phi)$ is a square matrix of dimension $3$. Its arrays depends on parameters $\phi$. We do not state these formula of $M_1(\phi)$ for simplicity, but we used the exact relations in further calculations and numerical computations. For admissible $\phi_0$ we have $\delta_0\in\Omega_0$, and
\begin{align}\label{stdf}
  \begin{pmatrix} \dot{u}\\\dot{v}\\\dot{w}\end{pmatrix}=
    \begin{pmatrix}
      0 & -\omega_0 & 0\\
      \omega_0 & 0 & 0\\
      0       &0   &\mathrm{a}_0
    \end{pmatrix}
    \begin{pmatrix} u\\v\\w\end{pmatrix} +
          \begin{pmatrix} f(u,v,w)\\g(u,v,w)\\h(u,v,w) \end{pmatrix}
\end{align}

After converting the system~\eqref{eq9} to the standard form~\eqref{stdf}, due to the existence of a pair of purely imaginary roots, the center manifold is the solution of a partial differential equation.
See \cite{khalil} for details of center manifold reduction in finite dimensional case. An approximation of $H$ up to the third order is 
\begin{equation}\label{eq13}
	w=H(u,v)=\theta  u^2+\kappa  u v+\iota  v^2+O_3(u,v),
\end{equation}
which represent the two-dimensional reduced system up to the fourth order; see~\cite{Carr} for center manifold theory and its calculations.
Then, the first Lyapunov coefficient is as follows; see~\cite{Guckenheimer}.

\begin{equation}\label{eq16}
	\begin{split}
		\sigma&=\frac{1}{16 \omega_0 }\bigg(\mathrm{\mathbf f}_{uv} (\mathrm{\mathbf f}_{uu}+\mathrm{\mathbf f}_{vv})-\mathrm{\mathbf f}_{uu} \mathrm{\mathbf g}_{uu}+\mathrm{\mathbf f}_{vv} \mathrm{\mathbf g}_{vv}-\mathrm{\mathbf g}_{uv} (\mathrm{\mathbf g}_{uu}+\mathrm{\mathbf g}_{vv})\bigg)+\frac{1}{16} \bigg(\mathrm{\mathbf f}_{uuu}+\mathrm{\mathbf f}_{uvv}+\mathrm{\mathbf g}_{uuv}+\mathrm{\mathbf g}_{vvv}\bigg)
	\end{split}
\end{equation}

Here $\mathrm{\mathbf f}(u,v)=f\left(u,v,H(u,v)\right)$ and $\mathrm{\mathbf g}(u,v)=g\left(u,v,H(u,v)\right)$.  
The equation~\eqref{eq16} for the first Lyapunov coefficient is a polynomial with a long expression. For simplicity, we show the polynomial as a function of competitive parameters $\alpha$ and $\beta$ as follows.
\begin{equation}\label{eq19}
	\begin{split}
		\sigma=&\sigma(\alpha,\beta),\qquad \alpha>0, \beta>0
	\end{split}.
\end{equation}

Then, to study the Lyapunov coefficient sign and decide on the type of Hopf bifurcation, except for two competitive parameters $\alpha$ and $\beta$, others are valued by the symmetric values as $a=b,c=d$, and $\mu=\nu$. Four cases must be considered as $\lambda_3<0\land \sigma<0$, $\lambda_3<0\land\sigma>0$, $\lambda_3>0\land \sigma<0$, and $\lambda_3>0\land \sigma>0$. \\

{\%bf
MacKevin et al. (1998) suggested that knowledge about interaction strength is vital in food webs study. They claimed that weak interactions may bind natural communities together and distributions of interaction strength are skewed towards weak interactions at most. Following this idea, we introduced two new parameters as relative competition strength.
Together with FRR, we use these three parameters to study asymmetricity in two–species competition.
  }
We varied parameter values in due course to compare the induced asymmetry in competition by introducing the following coefficients. Superscripts $-$ and $+$ denote the states before and after bifurcation.
\begin{equation}\label{eq17}
	\begin{split}
		&\eta_{yz}=\alpha / \beta,\qquad \eta_{zy}=\beta / \alpha\\
	\end{split}
\end{equation}

\subsubsection{Hopf case 1,  $\lambda_3<0, \sigma<0$}
\begin{figure}[h!]
	\centering
	\begin{subfigure}[b]{0.3\linewidth}
		\includegraphics[width=\linewidth]{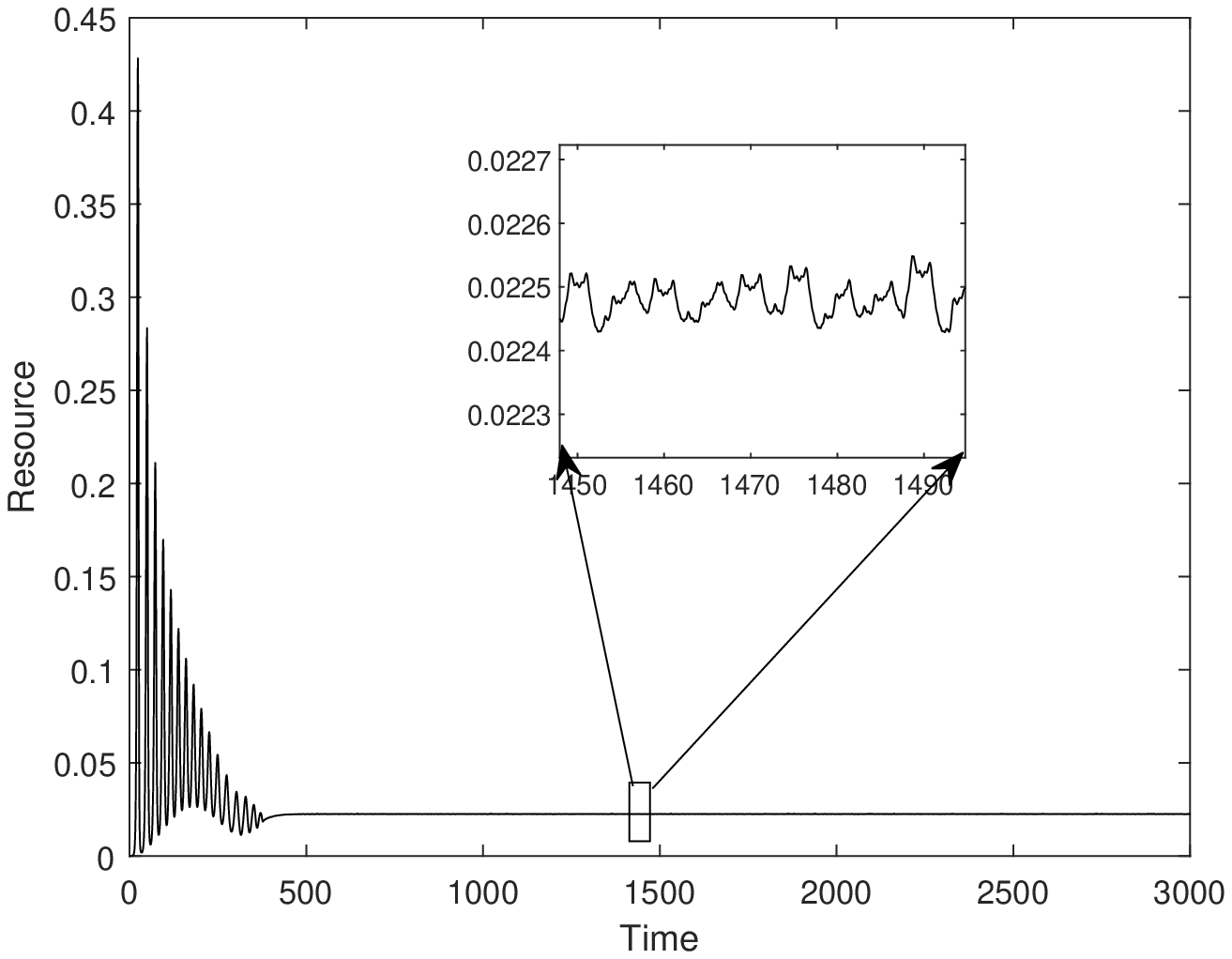}
		\caption{Resource success}
	\end{subfigure}
	\begin{subfigure}[b]{0.3\linewidth}
		\includegraphics[width=\linewidth]{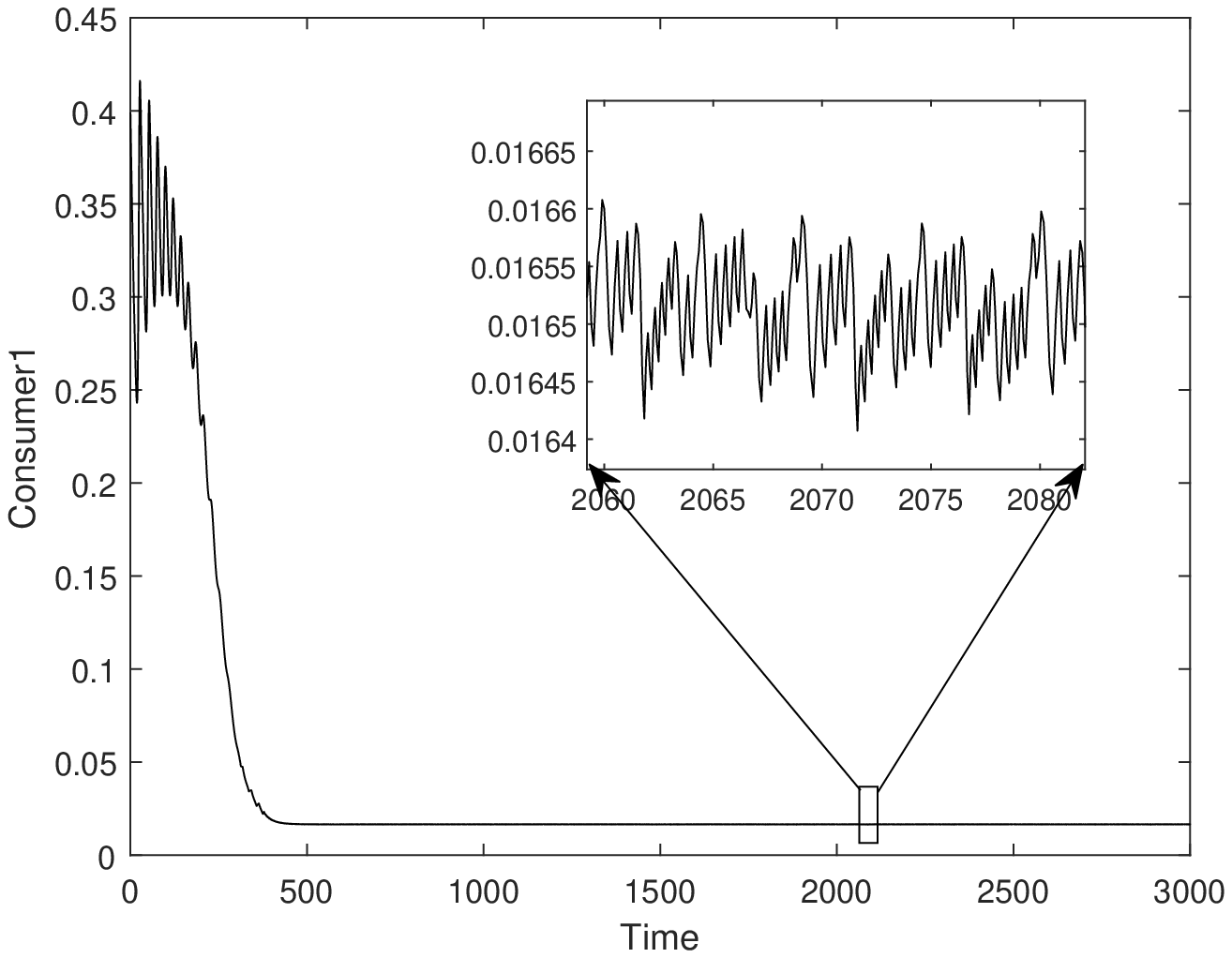}
		\caption{Consumer 1 success}
	\end{subfigure}
	\begin{subfigure}[b]{0.3\linewidth}
		\includegraphics[width=\linewidth]{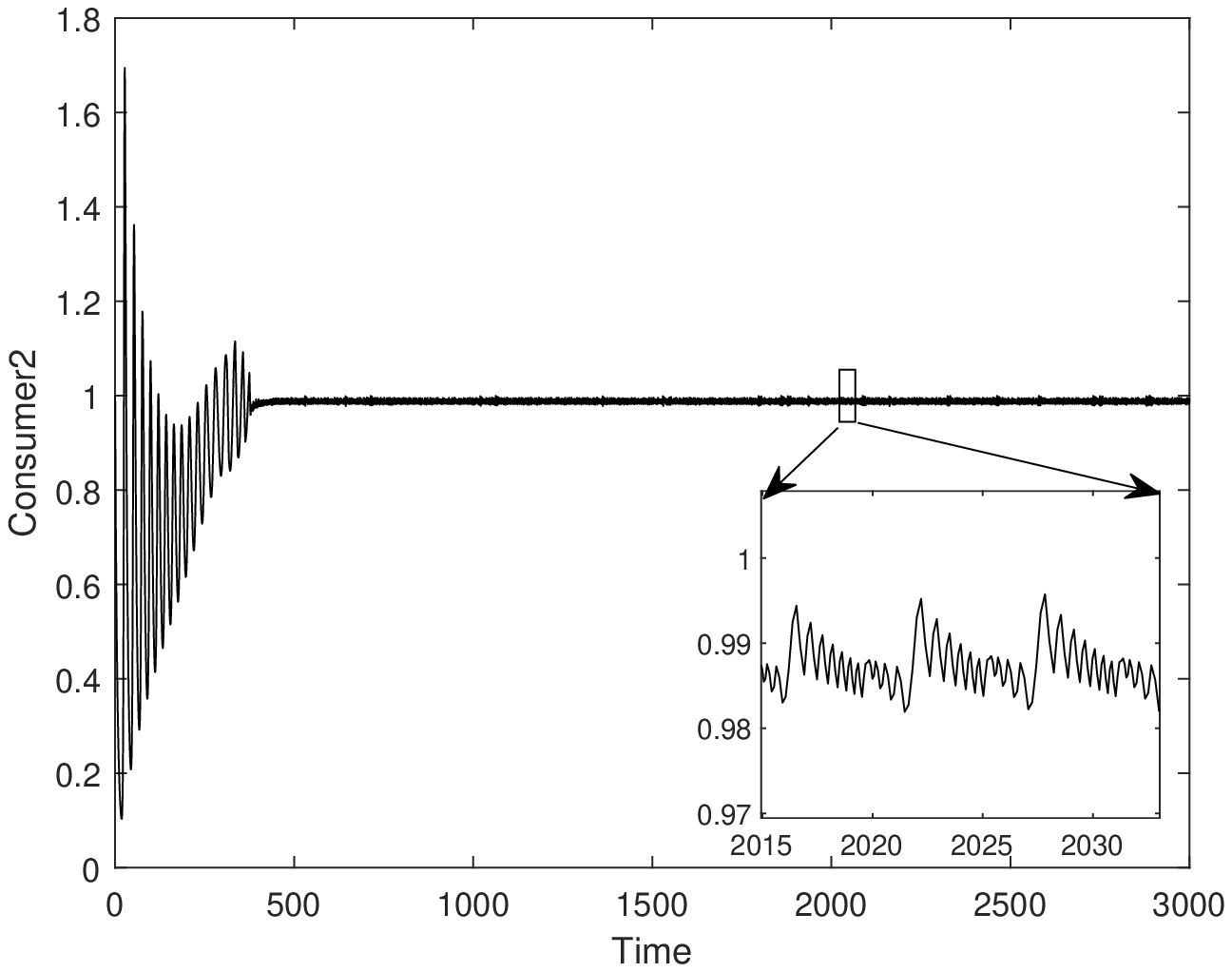}
		\caption{Consumer 2 success}
	\end{subfigure}
	\begin{subfigure}[b]{0.3\linewidth}
		\includegraphics[width=\linewidth]{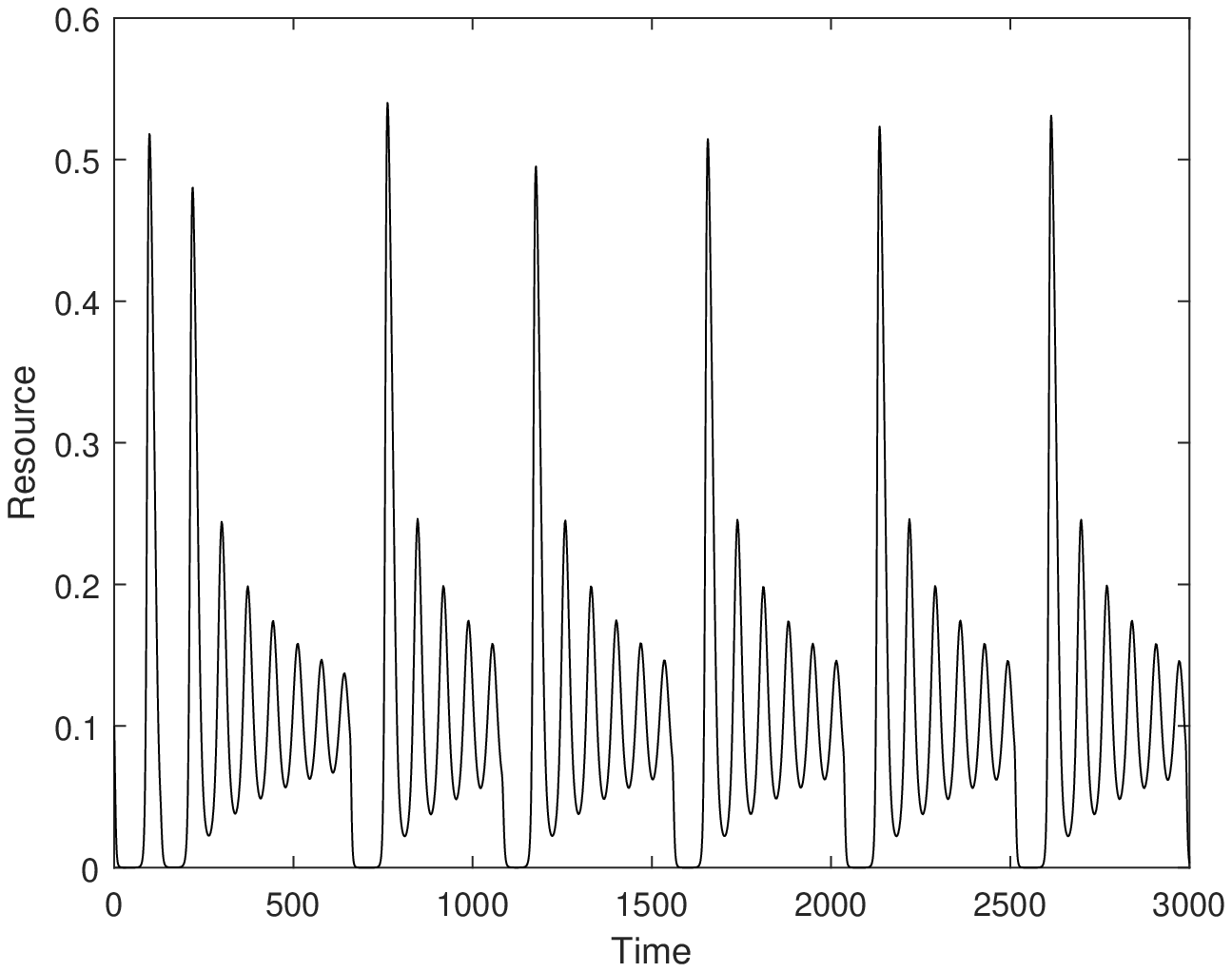}
		\caption{Resource success}
	\end{subfigure}
	\begin{subfigure}[b]{0.3\linewidth}
		\includegraphics[width=\linewidth]{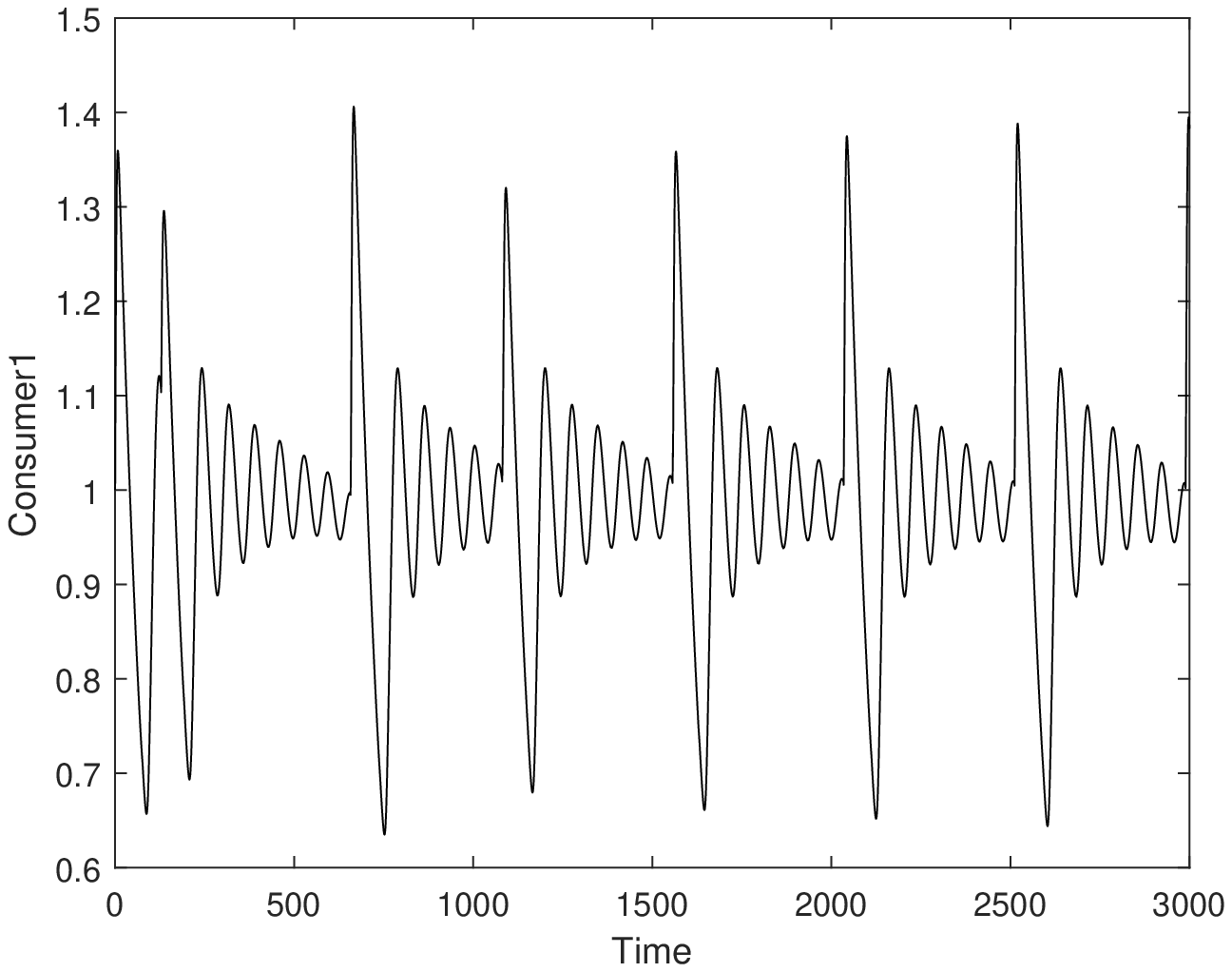}
		\caption{Consumer 1 success}
	\end{subfigure}
	\begin{subfigure}[b]{0.3\linewidth}
		\includegraphics[width=\linewidth]{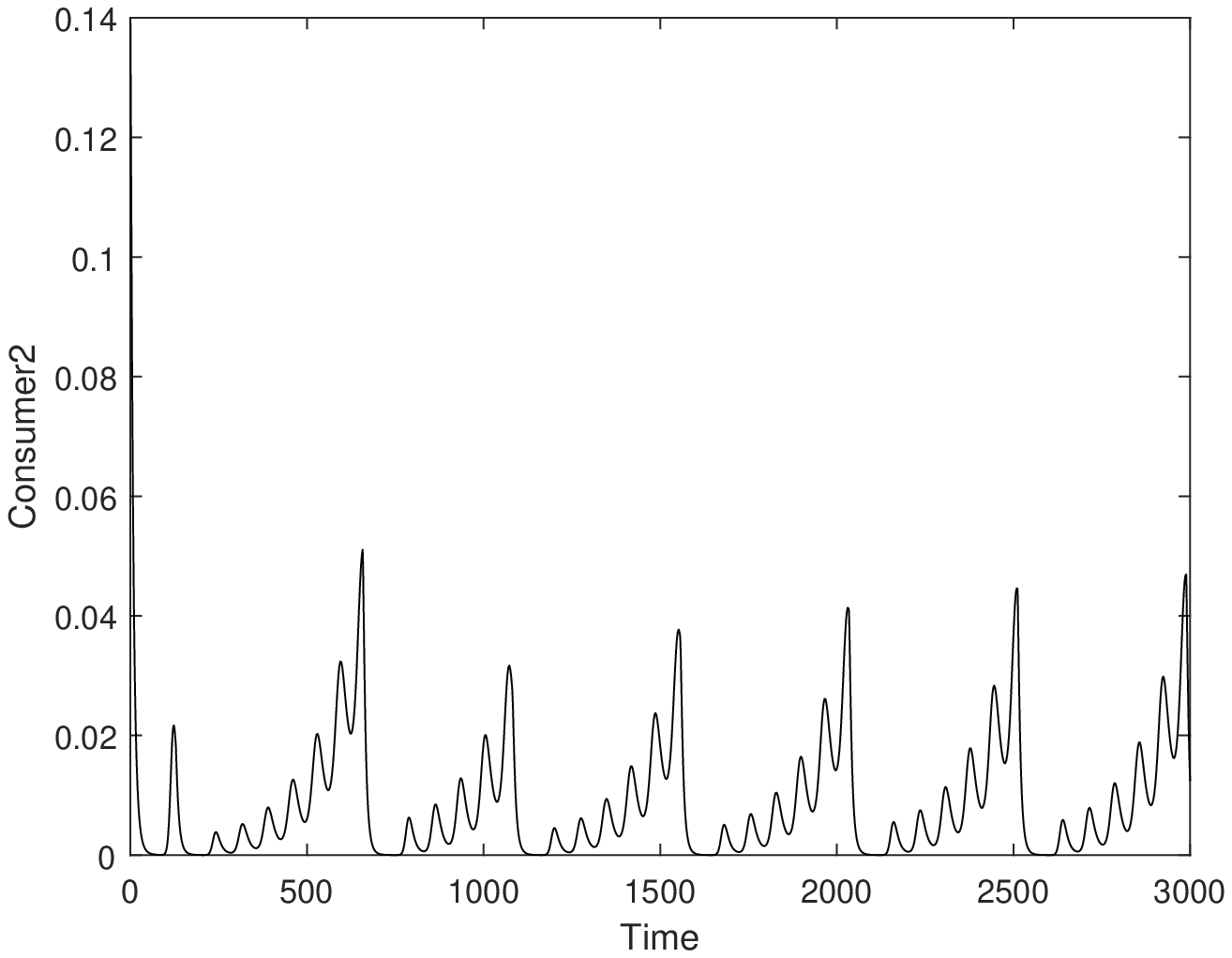}
		\caption{Consumer 2 success}
	\end{subfigure}
	\caption{Temporal dynamics of system \eqref{eq6} for the Hopf bifurcation case 1,  $\lambda_3<0, \sigma<0$, befor bifurcation point (a,b,c) and after bifurcation point (d,e,f). A relaxation-oscillation and a repeated damping oscillation pattern emerge before and after bifurcation point, respectively}
	\label{fig2}
\end{figure}

Figure \ref{fig2} (a, b, c) shows the temporal dynamics of the system \eqref{eq6} success before supercritical Hopf bifurcation ($\mu=-0.00134$) that represent relaxation-oscillation of the system for competitive coefficients $\alpha =1\times 10^{-11}$ and $\beta =0.34$, undergoing negative first Lyapunov coefficient, $\sigma = -0.049$ and $\lambda_3=-5.16 \times 10^{-11}$. Initial values for this simulation is $(0.0001,0.399,1)$ and the $\mathrm{FRR}$' threshold value is $0.00001$ with +2 improvement of attack rate after relaxation. In this case, consumer 2 is more effective in competition than consumer 1 with $\eta_{zy}^-= 3.4 \times 10^{10}$. Parameters symmetric values are $a=b=1.2$, $c=d=0.5$, and $\mu = \nu = 0.03$ in this simulation. Figure \ref{fig2} (d, e, f) shows the temporal dynamics of system \eqref{eq6} components for competitive coefficients $\alpha =0.9\times 10^{-12}$ and $\beta =0.091$, undergoing negative first Lyapunov coefficient, $\sigma = -0.231$ and $\lambda_3=-2.655 \times 10^{-11}$ and after supercritical Hopf bifurcation ($\mu=+0.00204$). Initial values for this simulation is $(0.1,1,0.14)$ and the $\mathrm{FRR}$' threshold value is 0.01 with +6 improvement of attack rate after relaxation. Again in this case, consumer 2 is more effective in competition than consumer 1 with $\eta_{zy}^+= 1.01 \times 10^{11}$. Parameters symmetric values are $a=b=1.2$, $c=d=0.11$, and $\mu = \nu = 0.01$ in this simulation. The figure \ref{fig2} (a, b, c) shows relaxation-oscillations and after the bifurcation point, the repeated stable limit cycles appear, figure \ref{fig2} (d, e, f). Negative $\lambda_3$ guarantees the absorption towards  equilibrium inside limit cycles. Resource-consumer 1 phase portrait after bifurcation point exhibits a kind of attractors similar to the R\"ossler one \cite{rossler76}. Through this bifurcation the value of $\eta_{i,j}$ changes from $3.4 \times 10^{10}$ to $1.01 \times 10^{11}$ so that the ratio of changes in asymmetry of the competition is $\eta_{zy}^+/\eta_{zy}^-=2.97$.

\begin{figure}[h!]
	\centering
	\begin{subfigure}[b]{0.3\linewidth}
		\includegraphics[width=\linewidth]{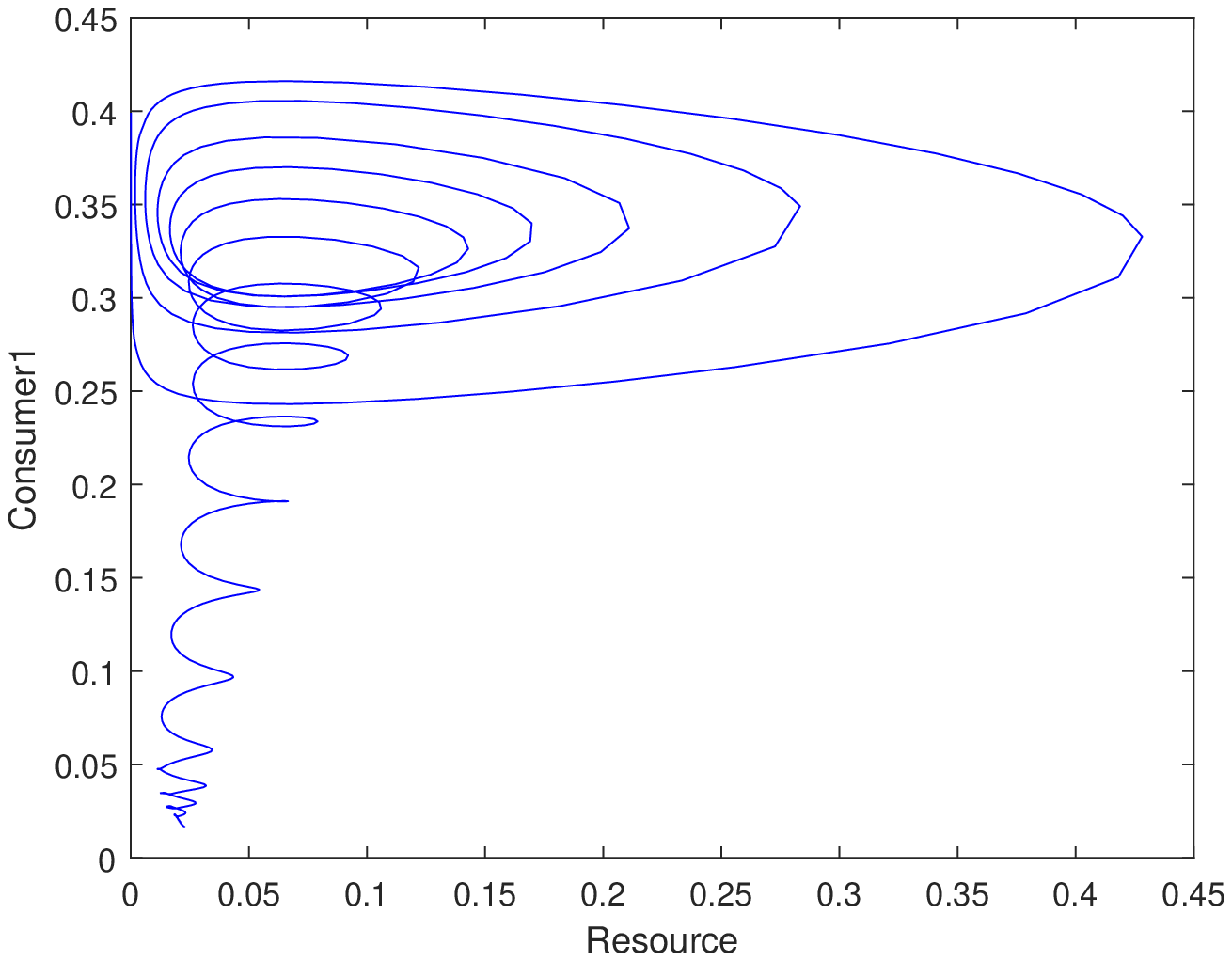}
		\caption{Resource- Consumer 1}
	\end{subfigure}
	\begin{subfigure}[b]{0.3\linewidth}
		\includegraphics[width=\linewidth]{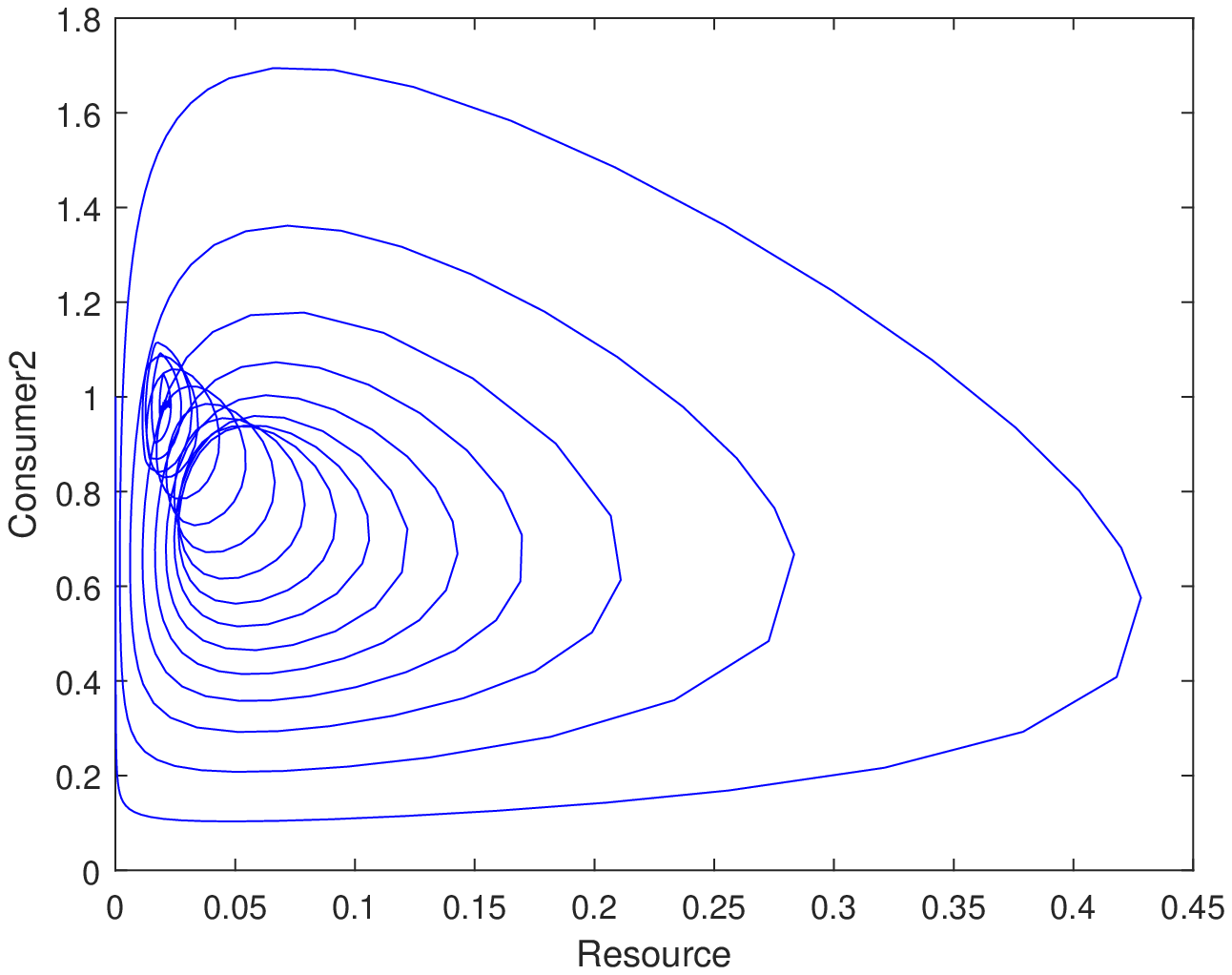}
		\caption{Resource- Consumer 2}
	\end{subfigure}
	\begin{subfigure}[b]{0.3\linewidth}
		\includegraphics[width=\linewidth]{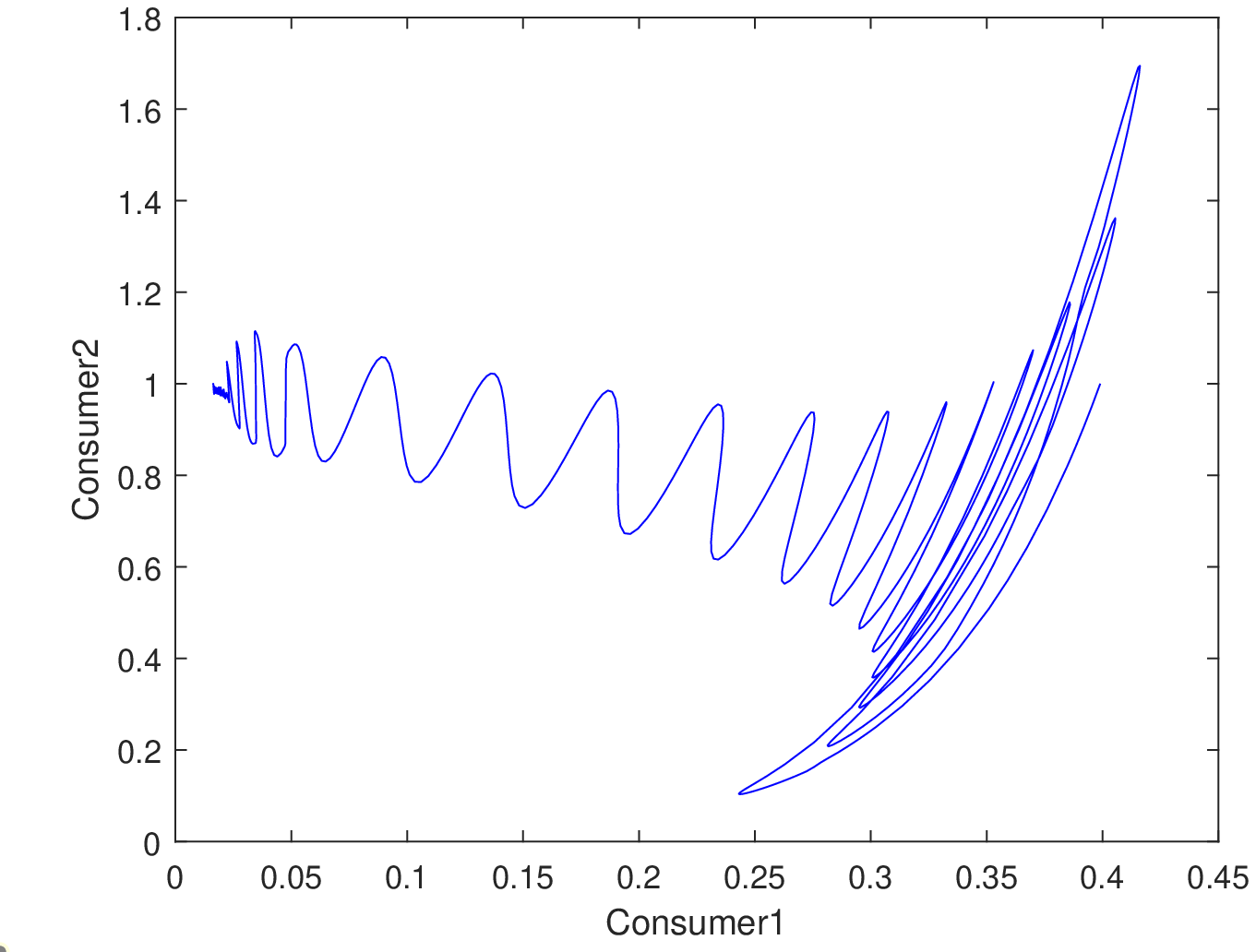}
		\caption{Consumer 1- Consumer 2}
	\end{subfigure}
	\begin{subfigure}[b]{0.3\linewidth}
		\includegraphics[width=\linewidth]{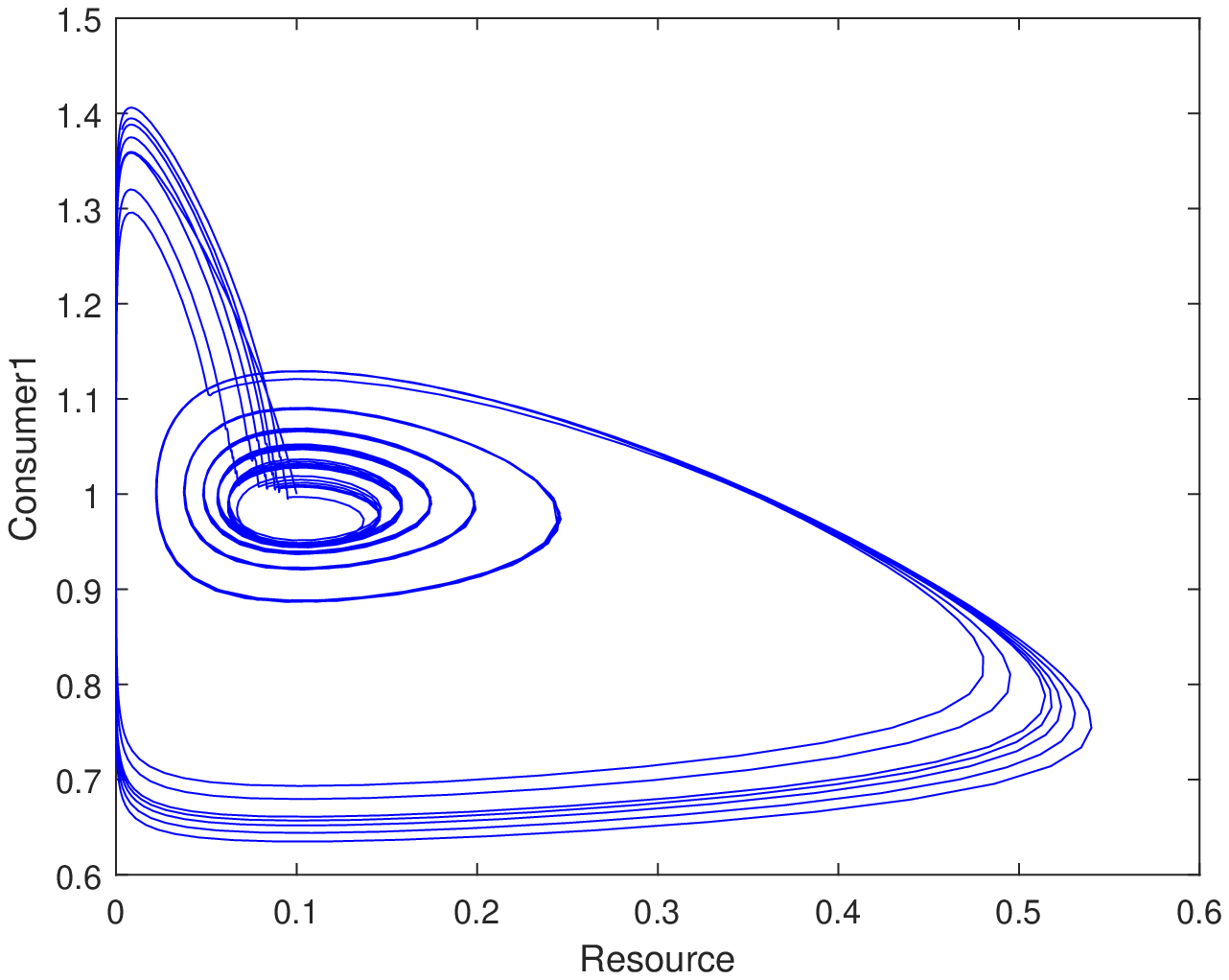}
		\caption{Resource- Consumer 1}
	\end{subfigure}
	\begin{subfigure}[b]{0.3\linewidth}
		\includegraphics[width=\linewidth]{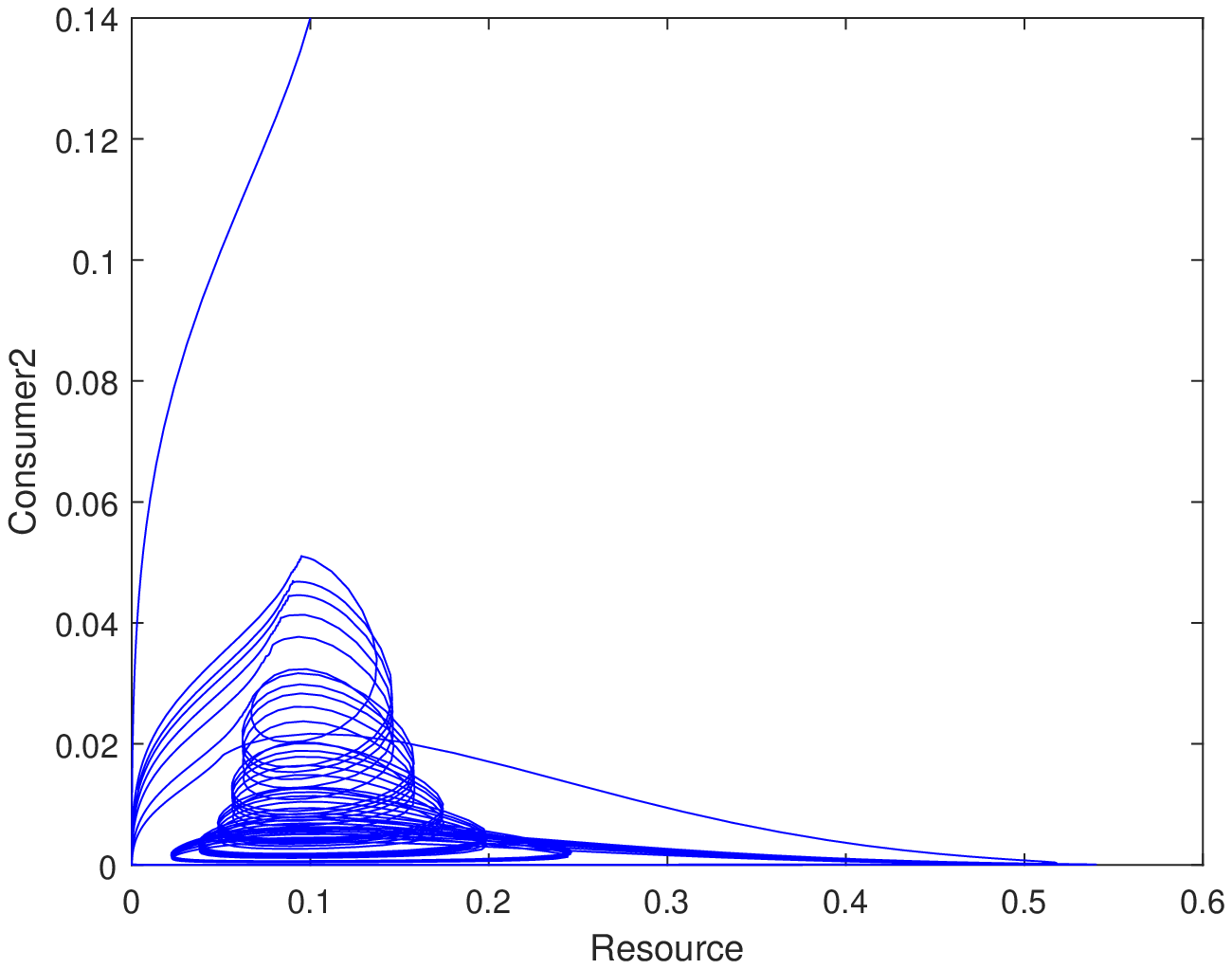}
		\caption{Resource- Consumer 2}
	\end{subfigure}
	\begin{subfigure}[b]{0.3\linewidth}
		\includegraphics[width=\linewidth]{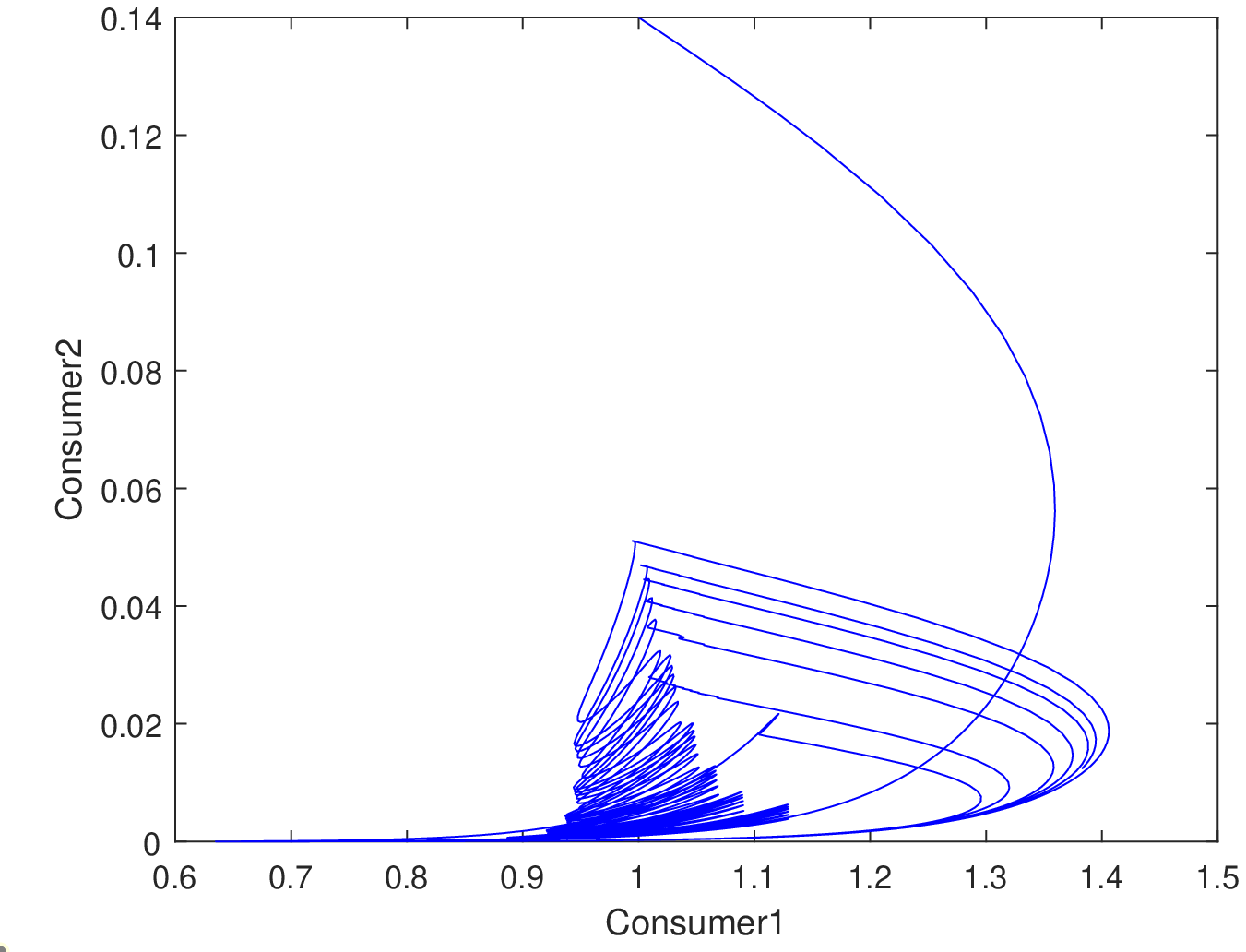}
		\caption{Consumer 1- Consumer 2}
	\end{subfigure}
	\caption{Two dimensional phase diagram of system \eqref{eq6} for the Hopf bifurcation case 1,  $\lambda_3<0, \sigma<0$, befor bifurcation point (a,b,c) and after bifurcation point (d,e,f)}
	\label{fig3}
\end{figure}

Figure \ref{fig3} shows two dimensional phase diagrams of the system \eqref{eq6} for the Hopf bifurcation of the supercritical type. Before bifurcation, the resource may experience extinction due to Alee effect. However, consumer 1 is close to the extinction in consumer 1-consumer 2 interaction. After bifurcation, with increasing resource resistance to predation, predation success of consumer 1 increases slowly, and in a turning point, predation success of consumer 1  increases with the resource decreased resistance to predatory.  

\begin{figure}[h!]
	\centering
	\begin{subfigure}[b]{0.4\linewidth}
		\includegraphics[width=\linewidth]{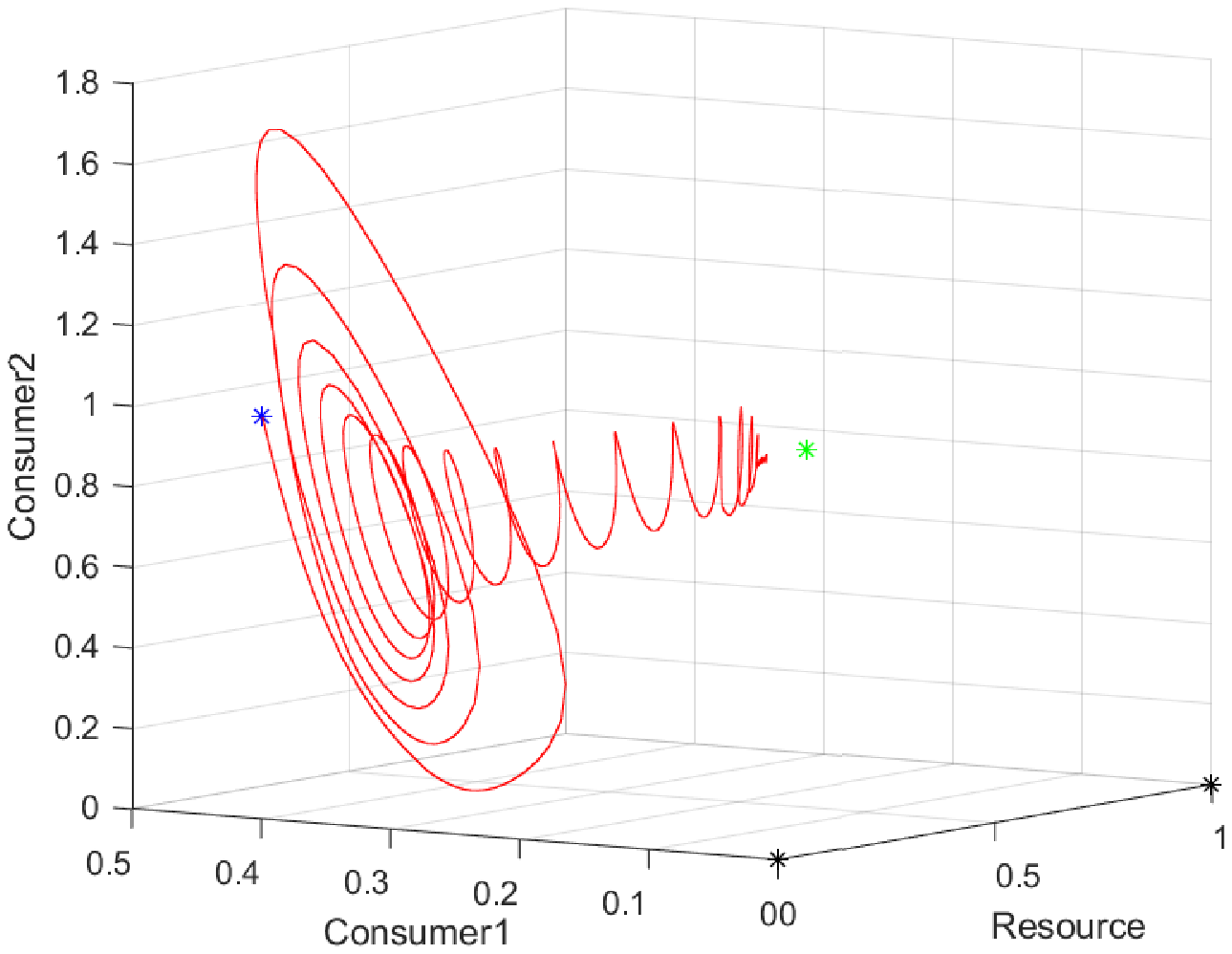}
		\caption{Before Hopf bifurcation}
	\end{subfigure}
	\begin{subfigure}[b]{0.4\linewidth}
		\includegraphics[width=\linewidth]{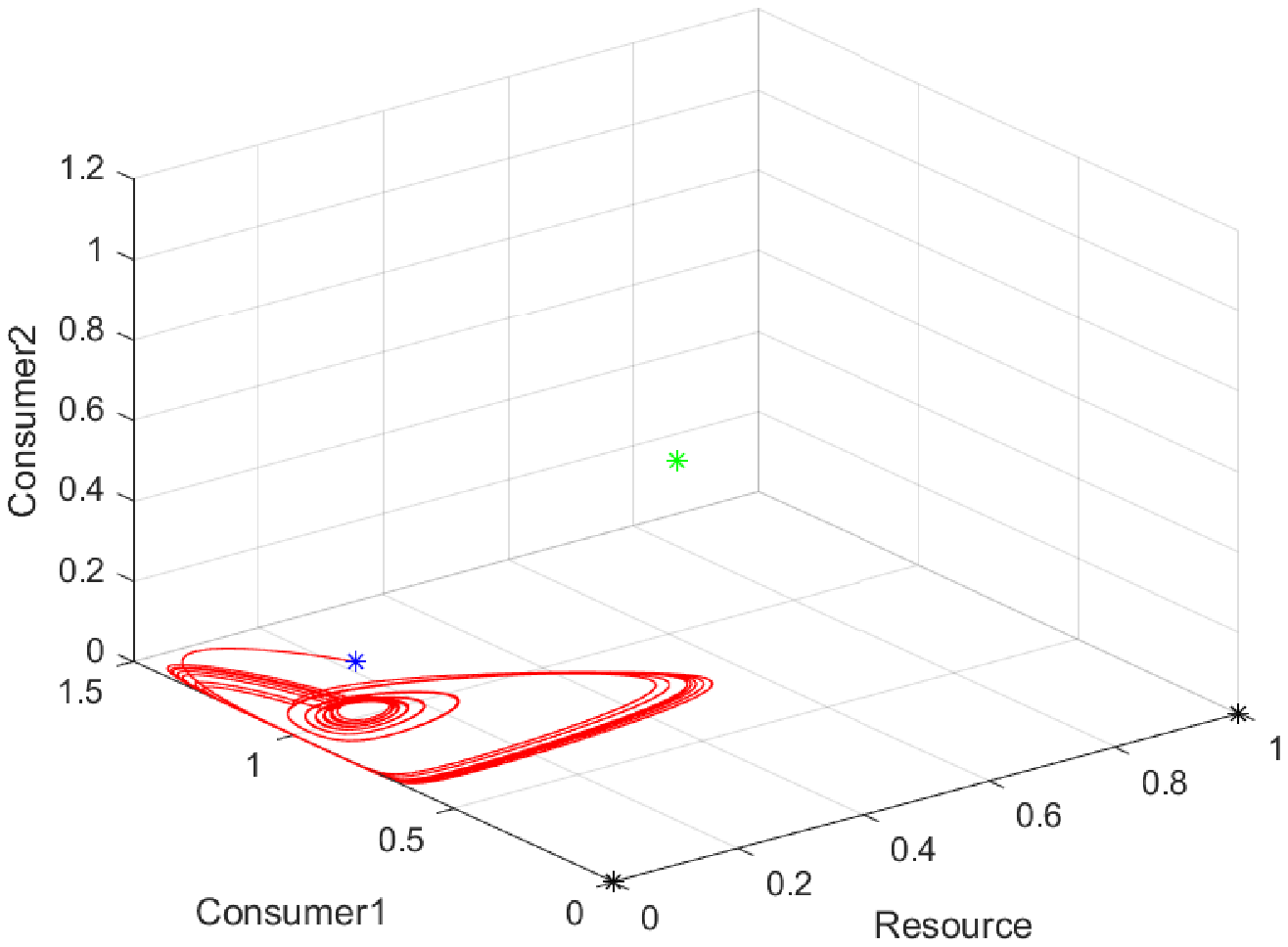}
		\caption{After Hopf bifurcation}
	\end{subfigure}
	\caption{Three dimensional phase diagram of system~\eqref{eq6} for the Hopf bifurcation case 1,  $\lambda_3<0, \sigma<0$. (a) Before bifurcation point. (b)  after bifurcation point. Blue, green and black stars are initial values, positive, and boundary equilibrium, respectively.}
	\label{fig4}
\end{figure}

Accordingly, figure \ref{fig4} shows the three dimensional phase space of the system \eqref{eq6} for the Hopf bifurcation of the supercritical type.  

\subsubsection{Hopf case 2, $\lambda_3>0, \sigma<0$}
\begin{figure}[h!]
	\centering
	\begin{subfigure}[b]{0.3\linewidth}
		\includegraphics[width=\linewidth]{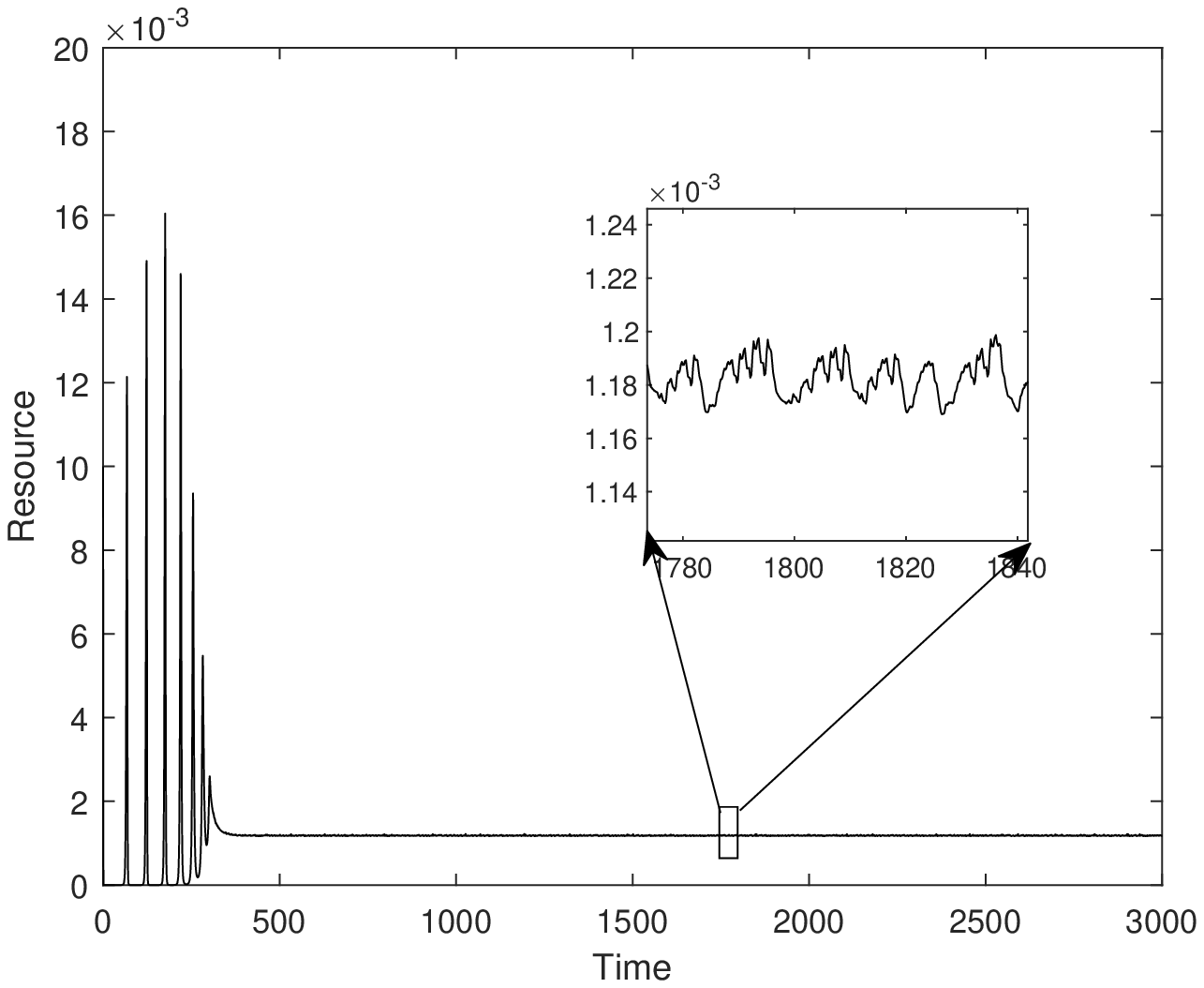}
		\caption{Resource success}
	\end{subfigure}
	\begin{subfigure}[b]{0.3\linewidth}
		\includegraphics[width=\linewidth]{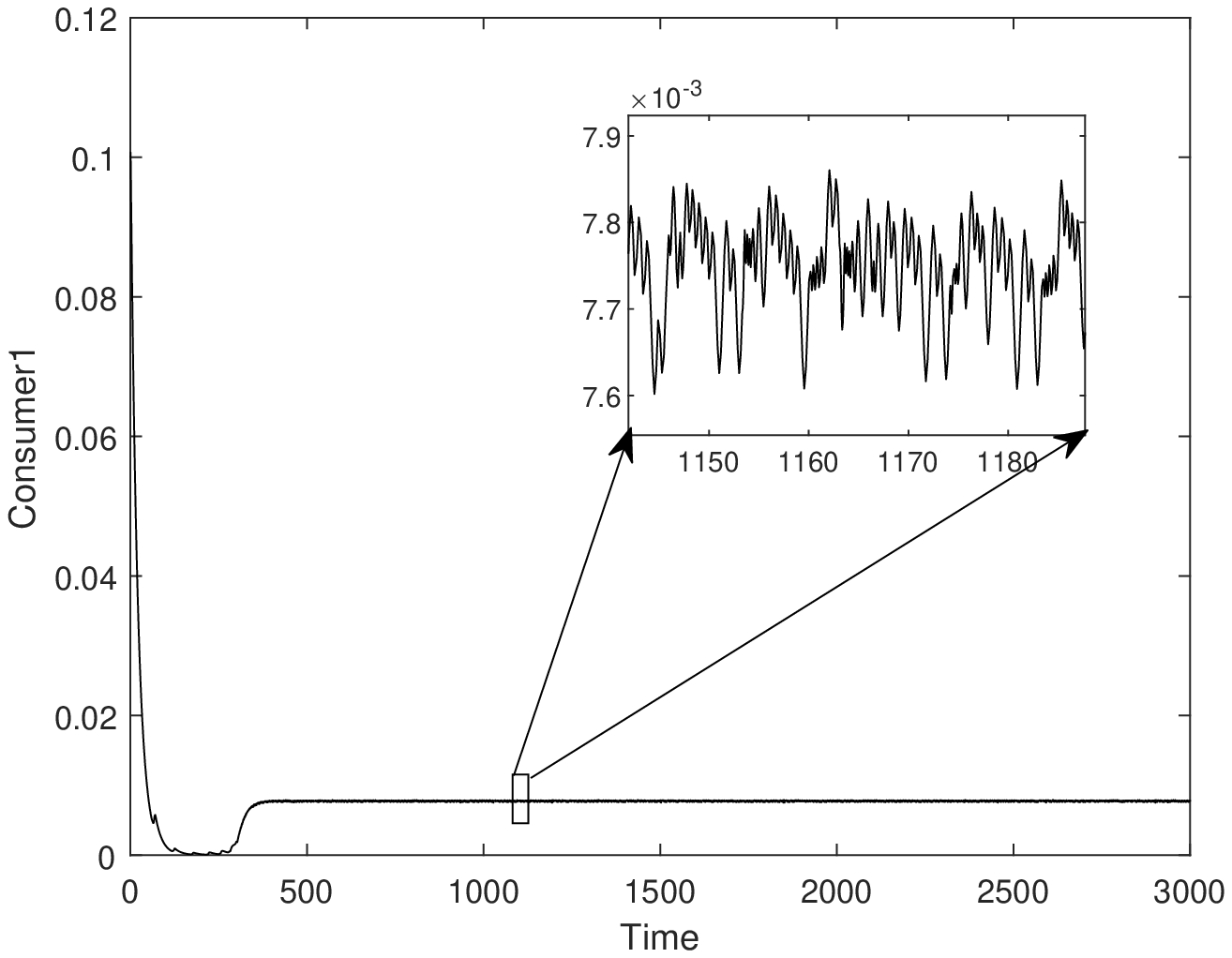}
		\caption{Consumer 1 success}
	\end{subfigure}
	\begin{subfigure}[b]{0.3\linewidth}
		\includegraphics[width=\linewidth]{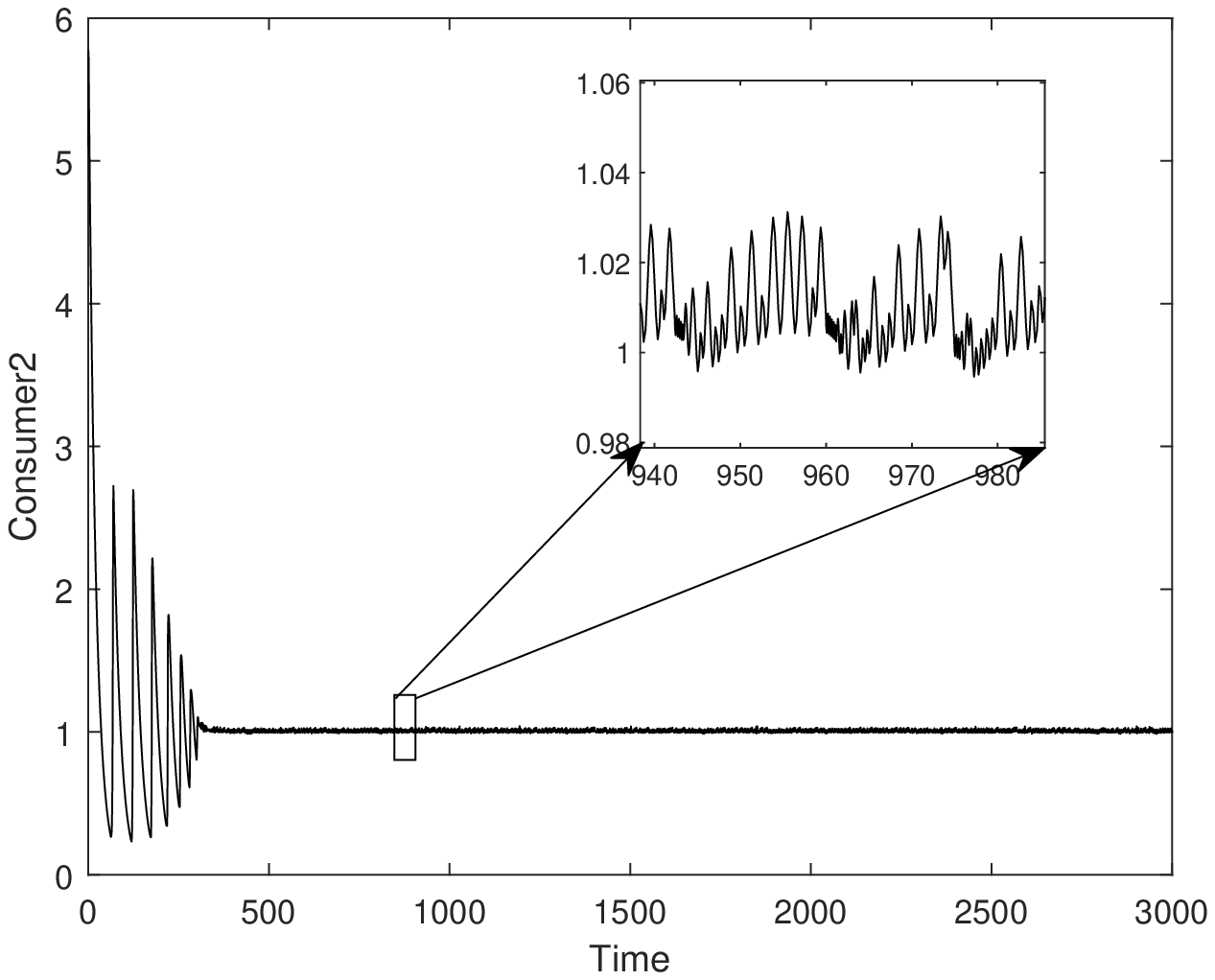}
		\caption{Consumer 2 success}
	\end{subfigure}
	\begin{subfigure}[b]{0.3\linewidth}
		\includegraphics[width=\linewidth]{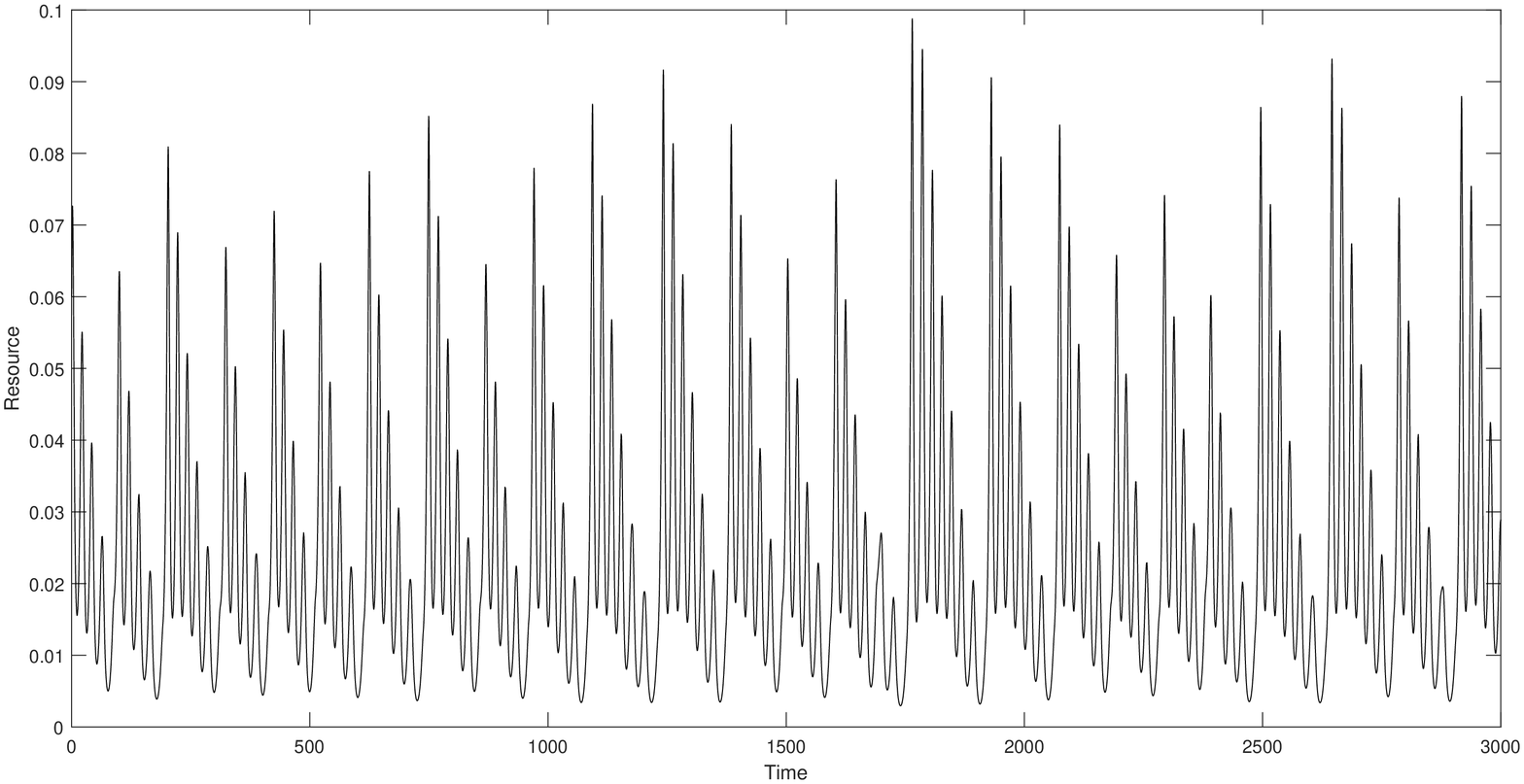}
		\caption{Resource success}
	\end{subfigure}
	\begin{subfigure}[b]{0.3\linewidth}
		\includegraphics[width=\linewidth]{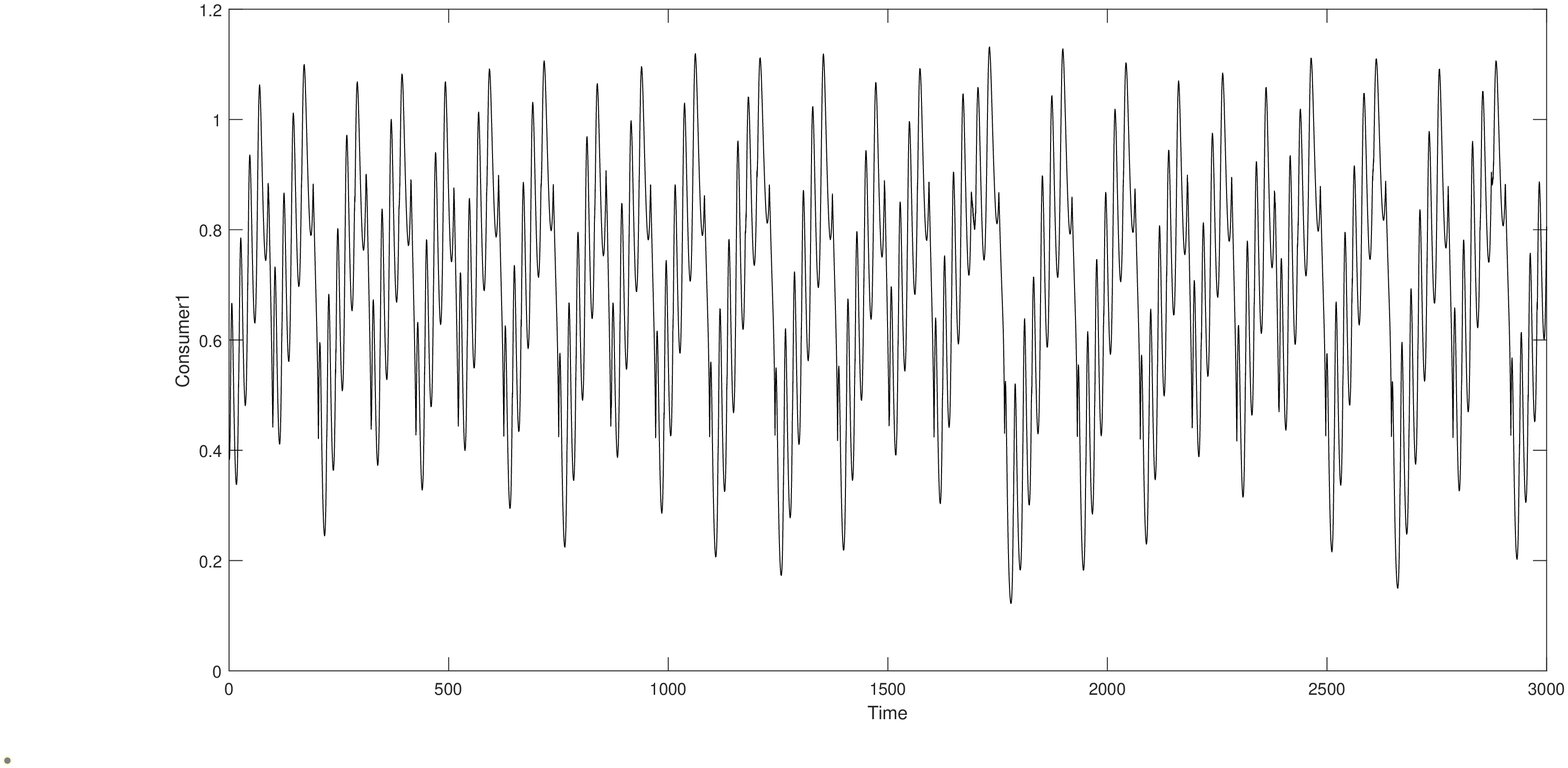}
		\caption{Consumer 1 success}
	\end{subfigure}
	\begin{subfigure}[b]{0.3\linewidth}
		\includegraphics[width=\linewidth]{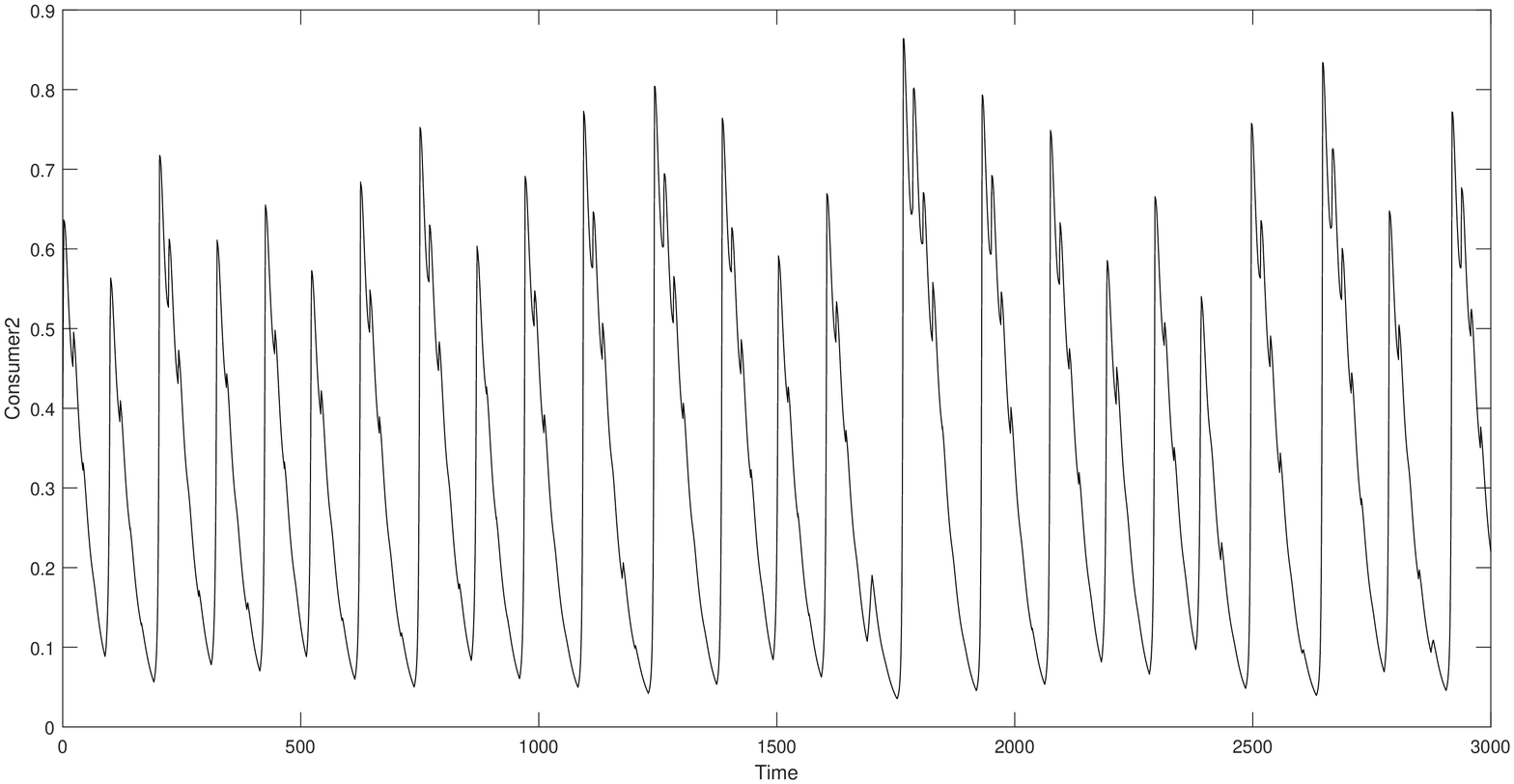}
		\caption{Consumer 2 success}
	\end{subfigure}
	\caption{Temporal dynamics of the system \eqref{eq6} for the Hopf bifurcation case 2, $\lambda_3>0, \sigma<0$, before bifurcation point (a,b,c) and after bifurcation point (d,e,f)}
	\label{fig5}
\end{figure}

\begin{figure}[h!]
	\centering
	\begin{subfigure}[b]{0.3\linewidth}
		\includegraphics[width=\linewidth]{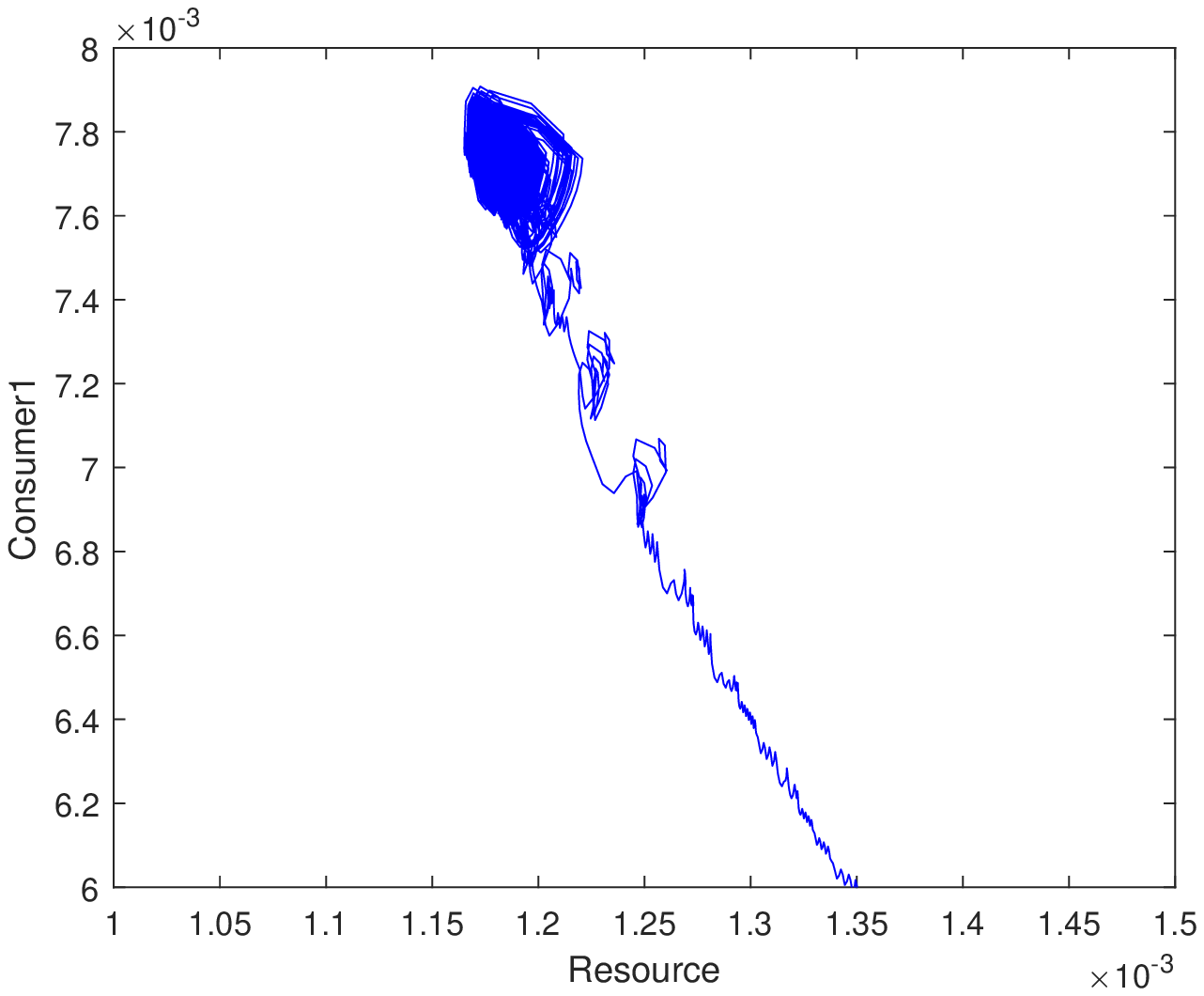}
		\caption{Resource-Consumer 1}
	\end{subfigure}
	\begin{subfigure}[b]{0.3\linewidth}
		\includegraphics[width=\linewidth]{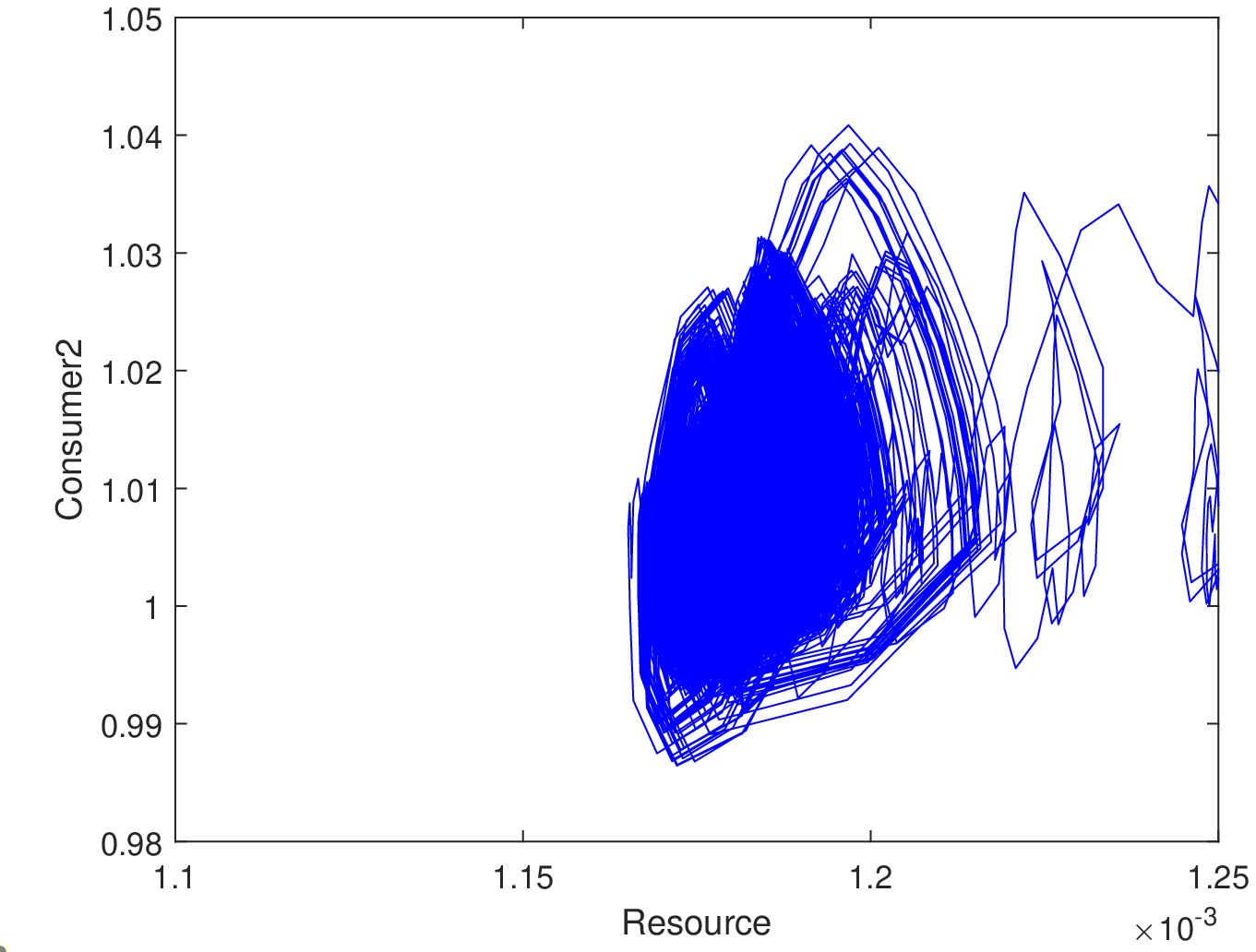}
		\caption{Resource-Consumer 2}
	\end{subfigure}
	\begin{subfigure}[b]{0.3\linewidth}
		\includegraphics[width=\linewidth]{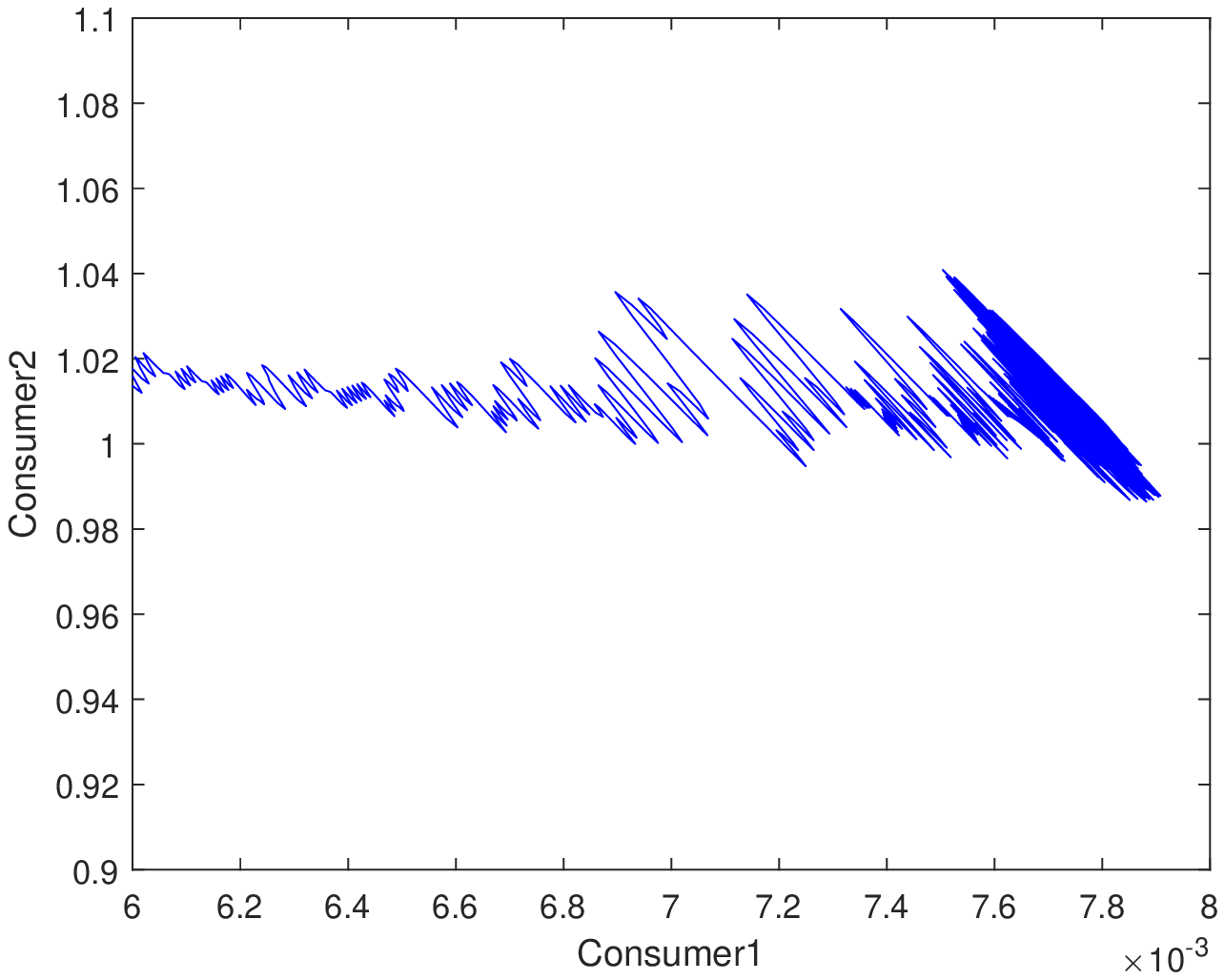}
		\caption{Consumer 1-Consumer 2}
	\end{subfigure}
	\begin{subfigure}[b]{0.3\linewidth}
		\includegraphics[width=\linewidth]{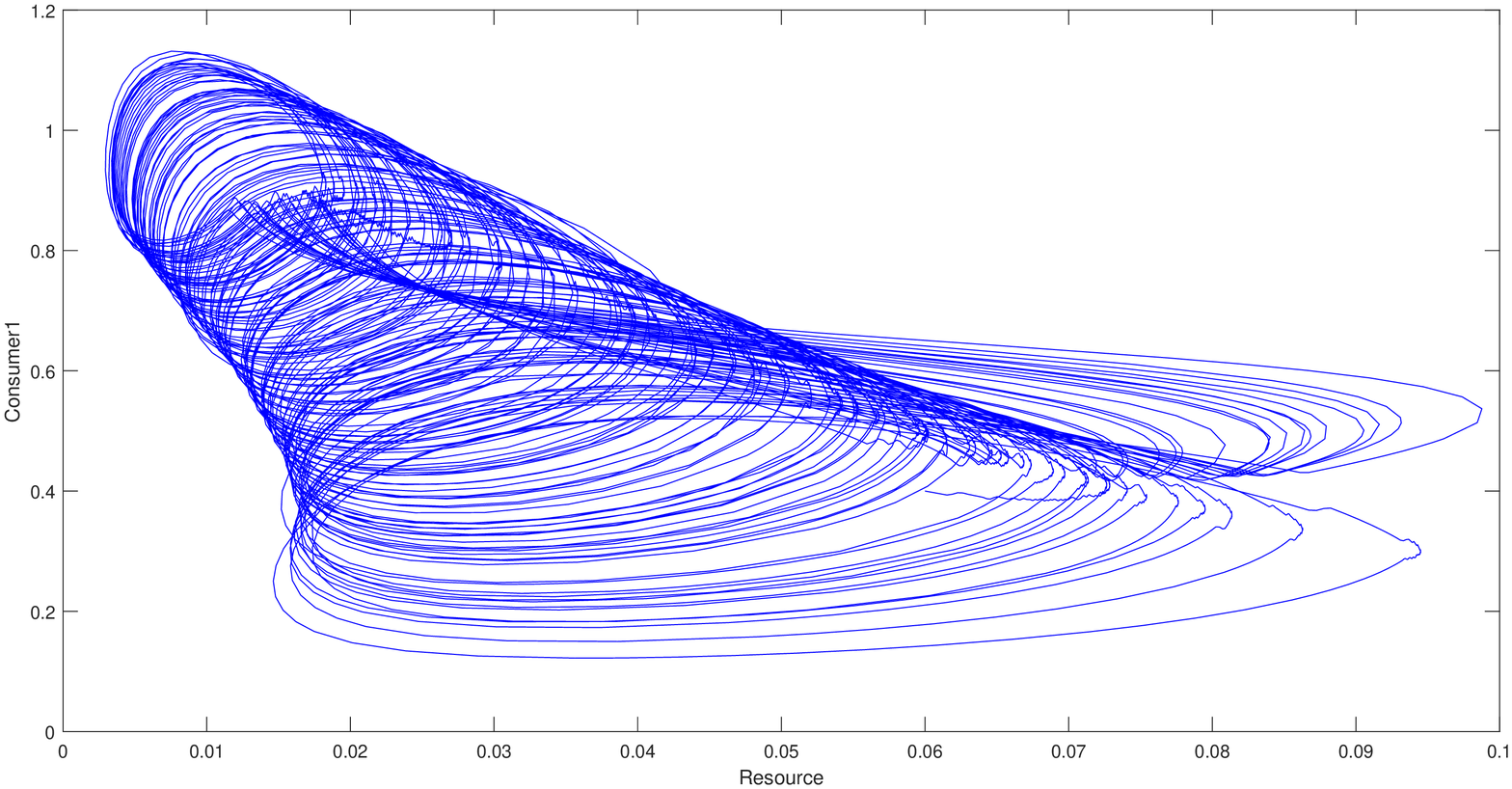}
		\caption{Resource-Consumer 1}
	\end{subfigure}
	\begin{subfigure}[b]{0.3\linewidth}
		\includegraphics[width=\linewidth]{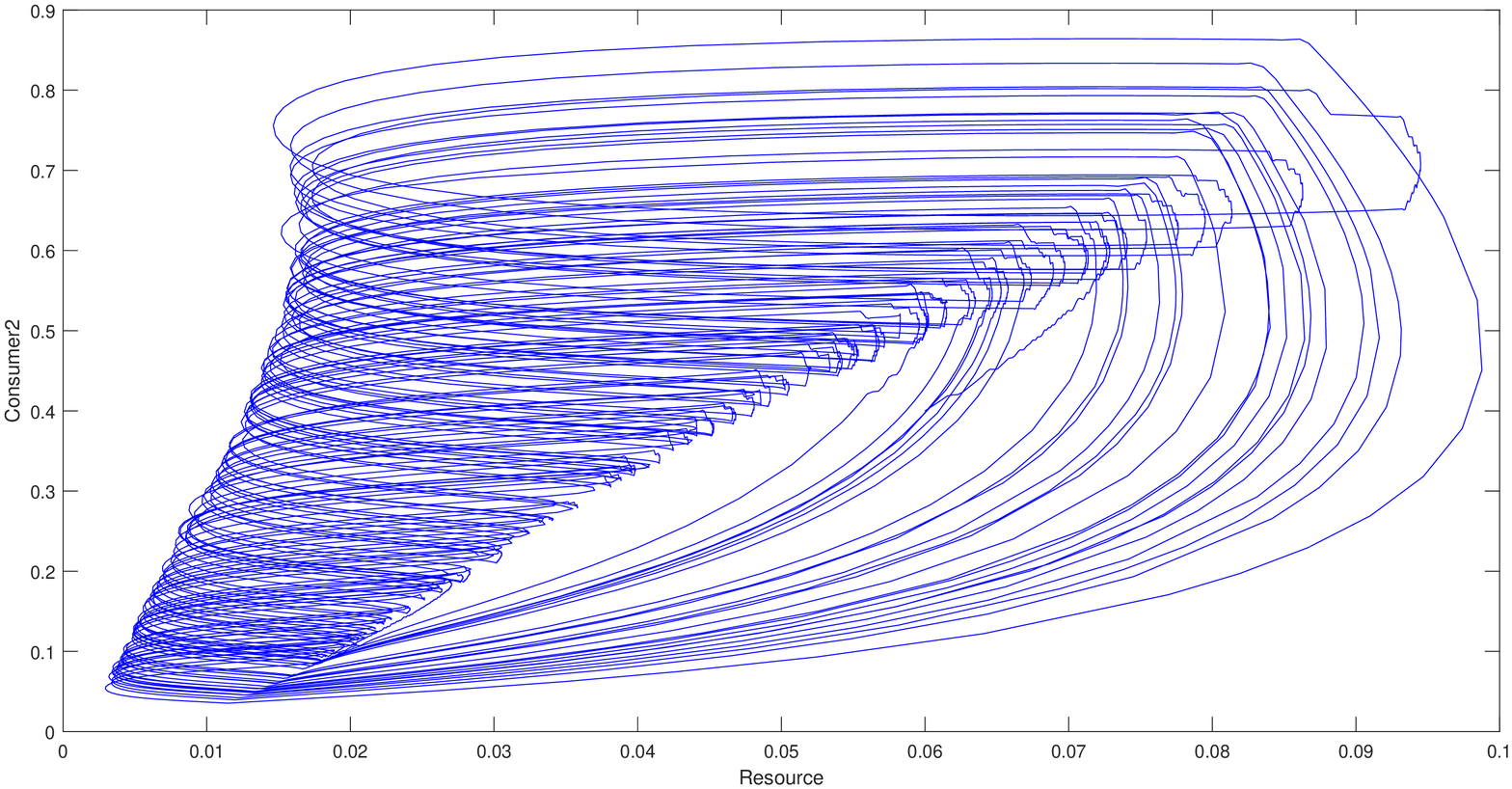}
		\caption{Resource-Consumer 2}
	\end{subfigure}
	\begin{subfigure}[b]{0.3\linewidth}
		\includegraphics[width=\linewidth]{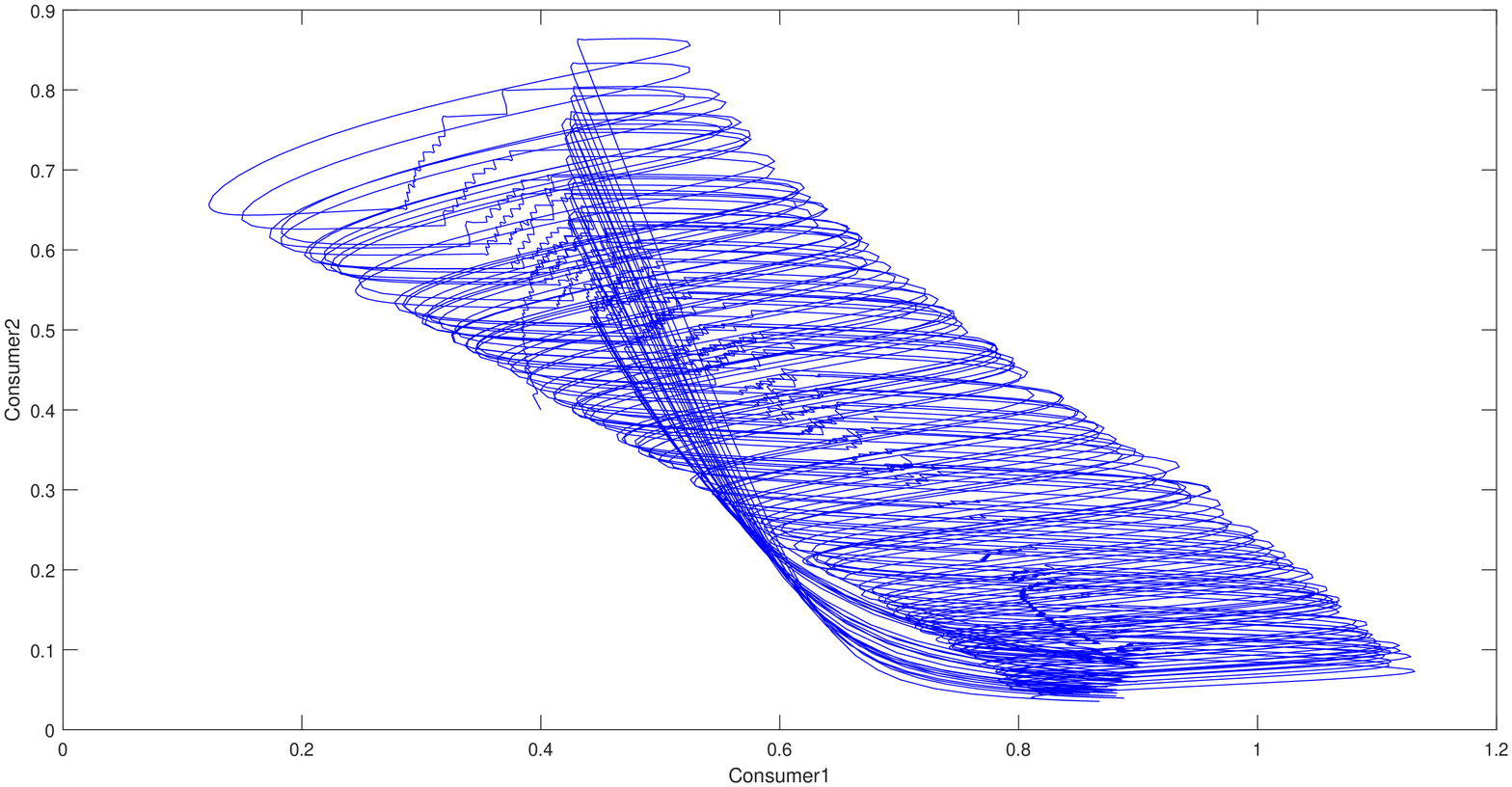}
		\caption{Consumer 1-Consumer 2}
	\end{subfigure}
	\caption{Two dimensional phase diagrams of the system \eqref{eq6} for the Hopf bifurcation case 2, $\lambda_3>0, \sigma<0$. (a,b,c) before Hopf bifurcation. (d,e,f) after Hopf bifurcation.}
	\label{fig6}
\end{figure}

\begin{figure}[h!]
	\centering
	\begin{subfigure}[b]{0.4\linewidth}
		\includegraphics[width=\linewidth]{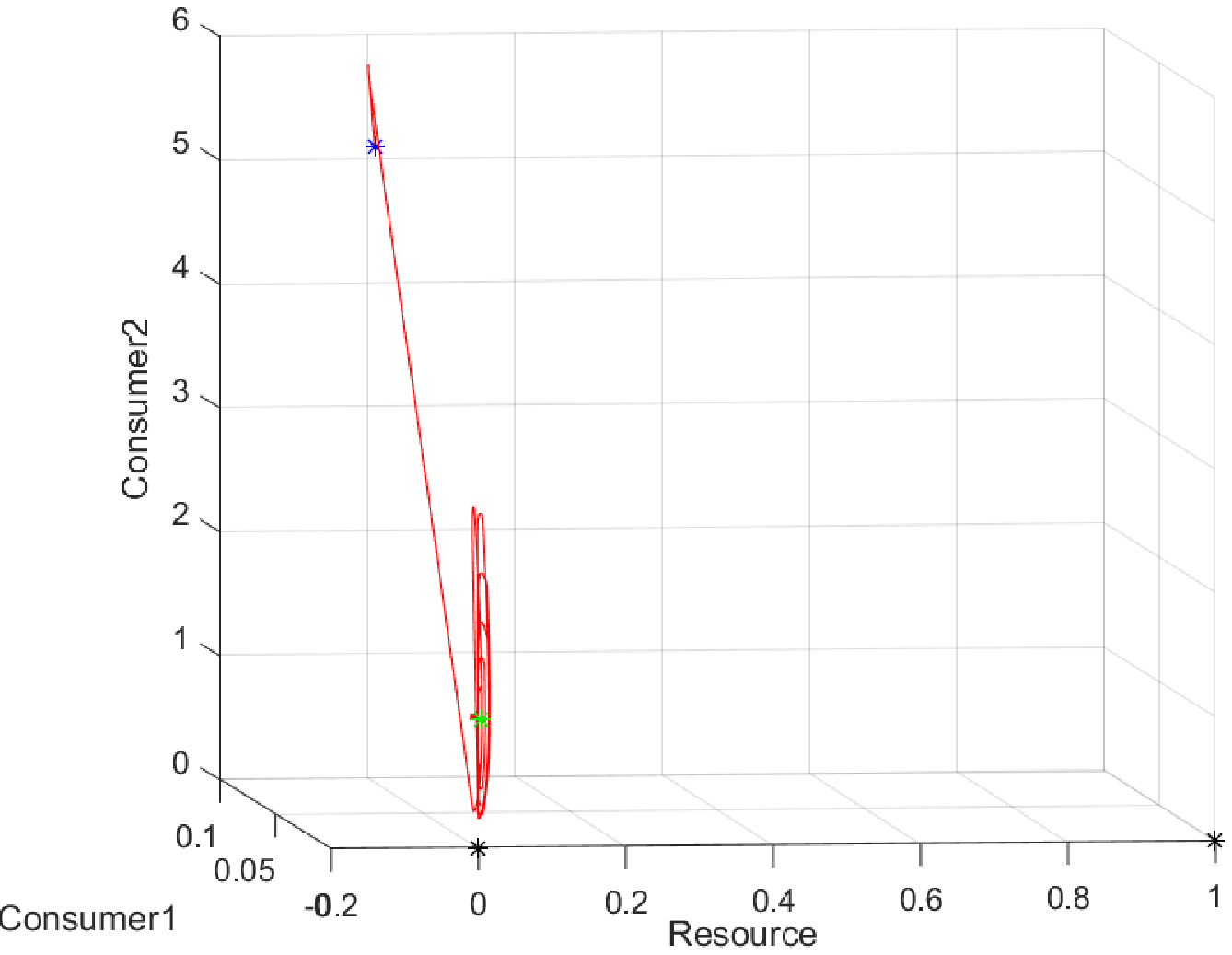}
		\caption{Before Hopf bifurcation}
	\end{subfigure}
	\begin{subfigure}[b]{0.5\linewidth}
		\includegraphics[width=\linewidth]{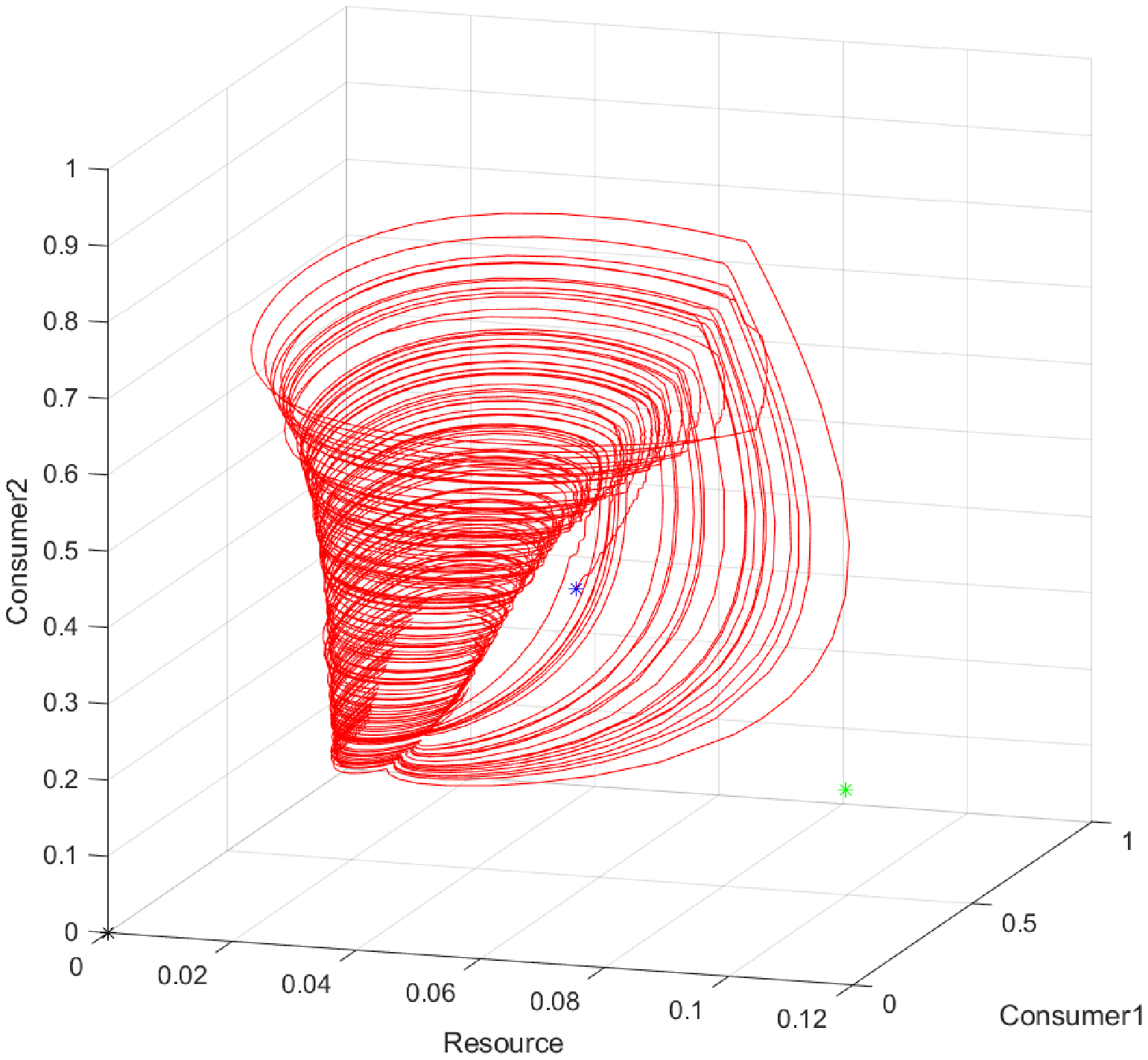}
		\caption{After Hopf bifurcation}
	\end{subfigure}
	\caption{Three dimensional phase diagram of the system \eqref{eq6} for the Hopf bifurcation in case 2 with $\lambda_3>0$ and $\sigma<0$. (a) before Hopf bifurcation. (b) after Hopf bifurcation. Blue, green and black stars are initial values, positive, and boundary equilibrium, respectively.}
	\label{fig7} 
\end{figure}

Competitive coefficients undergoing positive real root, $\lambda_3 >0$, negative first Lyapunov coefficient $\sigma < 0$ and supercritical Hopf bifurcation is used to simulate the figures~\ref{fig5} to~\ref{fig7}. Figure~\ref{fig5} (a,b,c) shows temporal dynamics of the system~\eqref{eq6} before supercritical Hopf bifurcation ($\mu=-0.0179$) for competitive coefficients $\alpha =0.9 \times 10^{-8}$ and $\beta =0.00041$, undergoing $\sigma = -0.08384$ and $\lambda_3=9.31 \times 10^{-9}$. Initial values for this simulation is $(0.01,0.1,5.1)$ and the $\mathrm{FRR}$' threshold value is 0.00001 with +60 improvement of attack rate after relaxation. In this case, consumer 1 is more effective in competition than consumer 2 with $\eta_{zy}^-=45555.6$. Symmetric parameters values are $a=b=10.5$, $c=d=13.11$, and $\mu = \nu = 0.05$ in this simulation. Figure~\ref{fig5} (d,e,f) shows temporal dynamics of the system components~\eqref{eq6} after supercritical Hopf bifurcation ($\mu=+0.0032$) for competitive coefficients $\alpha =0.3$ and $\beta =0.007$, undergoing negative first Lyapunov coefficient, $\sigma = -53.55$ and positive real root $\lambda_3=0.007$. Initial values for this simulation is $(0.06,0.4,0.4)$ and the $\mathrm{FRR}$' threshold value is 0.01 with +6 improvement of attack rate after relaxation. In this case, consumer 2 is more effective in competition than consumer 1 with $\eta_{yz}^+=42.86$. Symmetric parameters values are $a=b=1.2$, $c=d=0.5$, and $\mu = \nu = 0.03$ in this simulation. All phase portraits after this supercritical Hopf bifurcation exhibit R\"ossler-like attractors, \cite{rossler76}, similar to those reported for the first time by Hastings and Powell~\cite{hasting91}. For this kind of chaotic dynamics, change of stability leads to $\eta_{zy}^-/\eta_{yz}^+ \approx 1000$ times decreased competitive asymmetricity. 

Figure~\ref{fig6} shows two dimensional phase diagrams of the system~\eqref{eq6} for the supercritical Hopf bifurcation.

Accordingly, figure~\ref{fig7} shows the three dimensional phase portrait of the system~\eqref{eq6} for supercritical Hopf bifurcation.

\subsubsection{Hopf case 3, $\lambda_3>0, \sigma>0$}
\begin{figure}[h!]
	\centering
	\begin{subfigure}[b]{0.3\linewidth}
		\includegraphics[width=\linewidth]{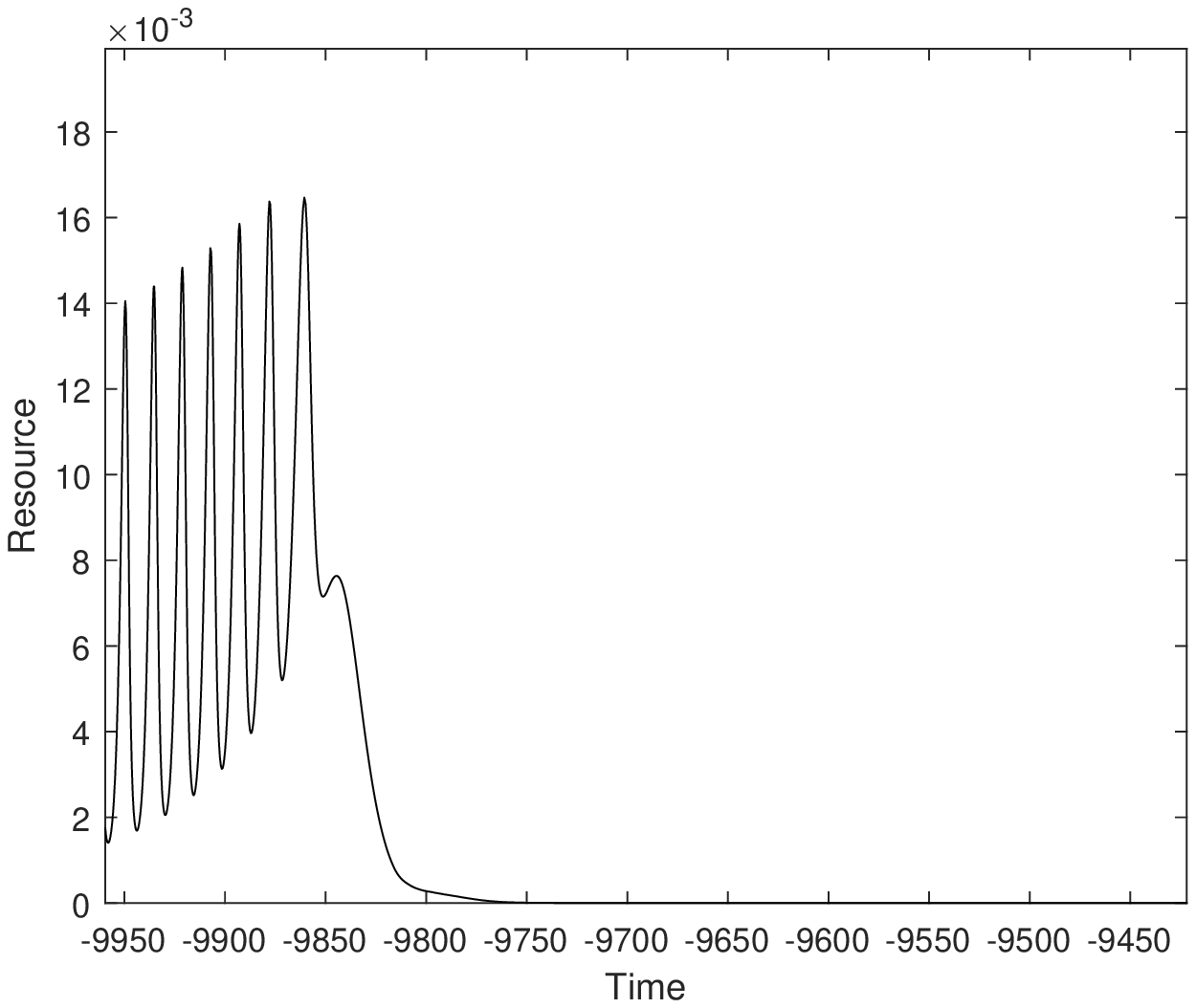}
		\caption{Resource time series}
	\end{subfigure}
	\begin{subfigure}[b]{0.3\linewidth}
		\includegraphics[width=\linewidth]{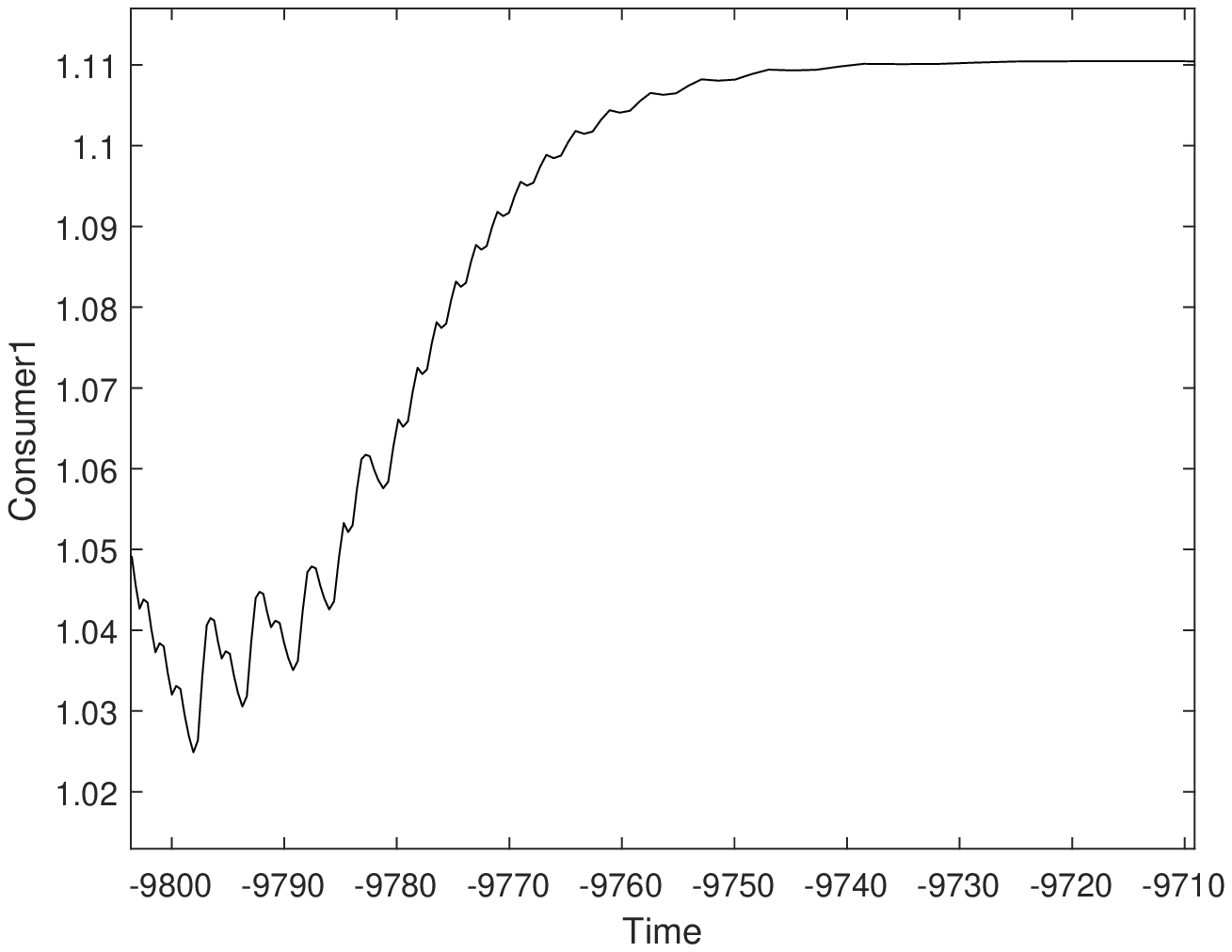}
		\caption{Consumer 1 time series}
	\end{subfigure}
	\begin{subfigure}[b]{0.3\linewidth}
		\includegraphics[width=\linewidth]{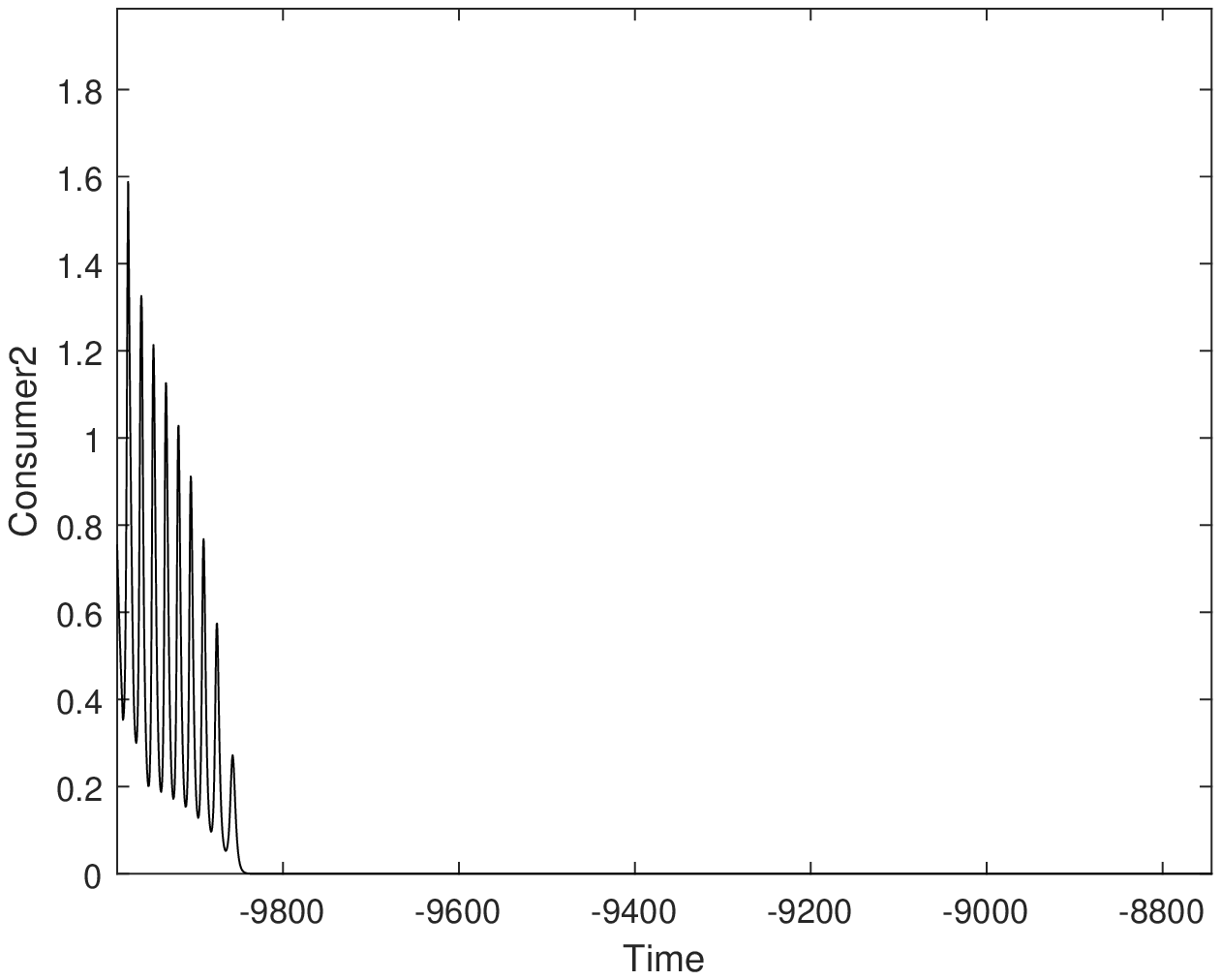}
		\caption{Consumer 2 time series}
	\end{subfigure}
	\begin{subfigure}[b]{0.3\linewidth}
		\includegraphics[width=\linewidth]{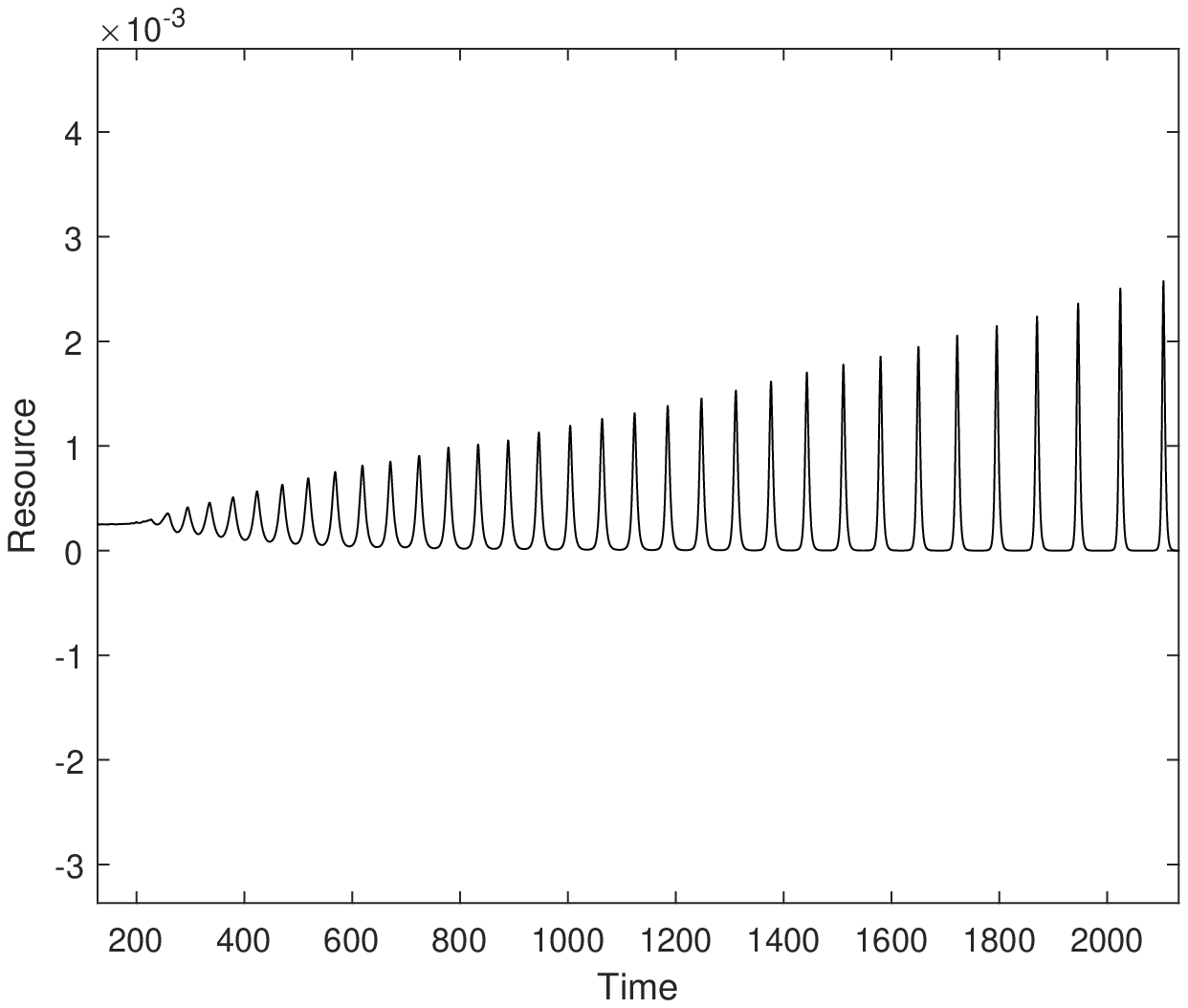}
		\caption{Consumer 2 time series}
	\end{subfigure}
	\begin{subfigure}[b]{0.3\linewidth}
		\includegraphics[width=\linewidth]{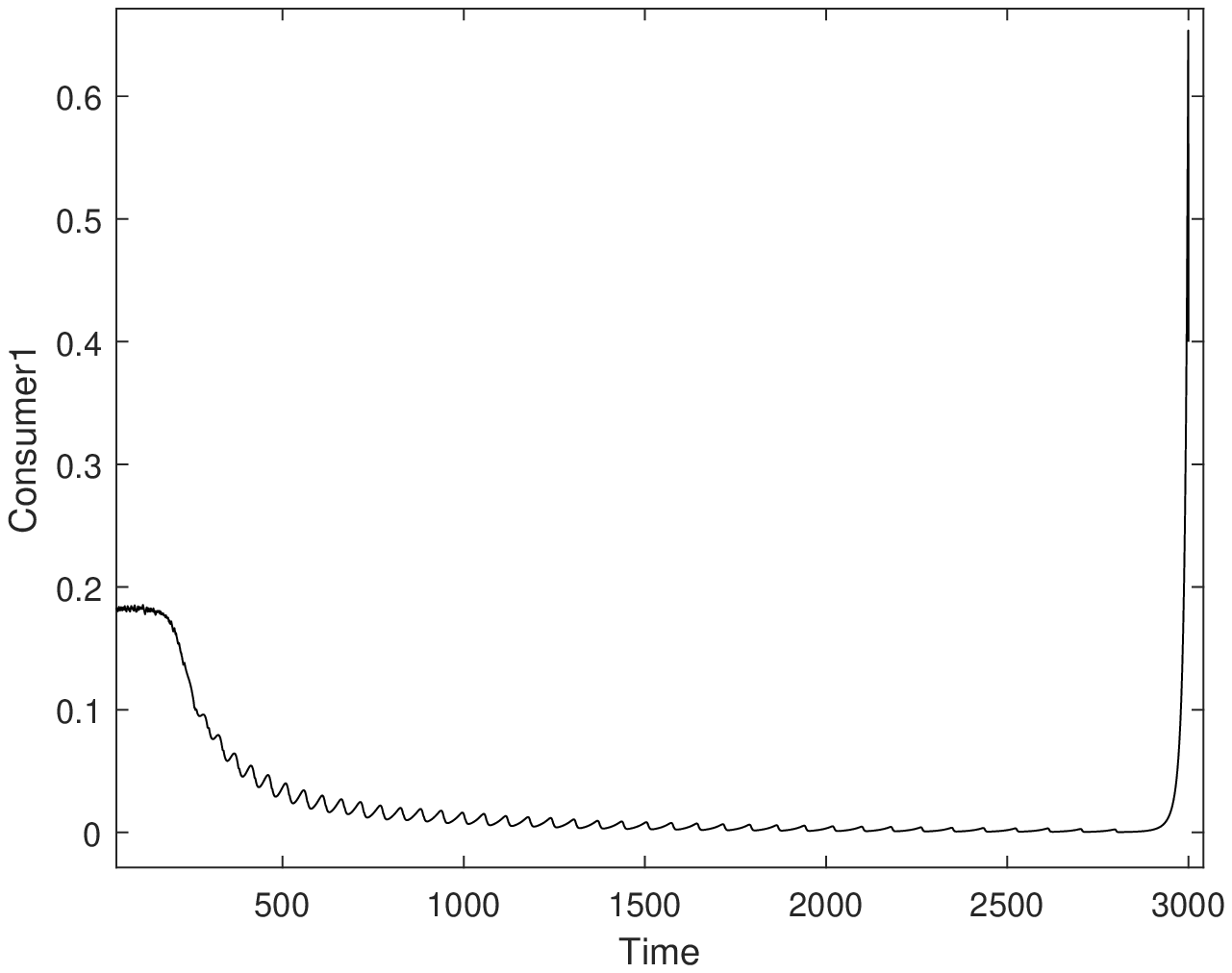}
		\caption{Consumer 2 time series}
	\end{subfigure}
	\begin{subfigure}[b]{0.3\linewidth}
		\includegraphics[width=\linewidth]{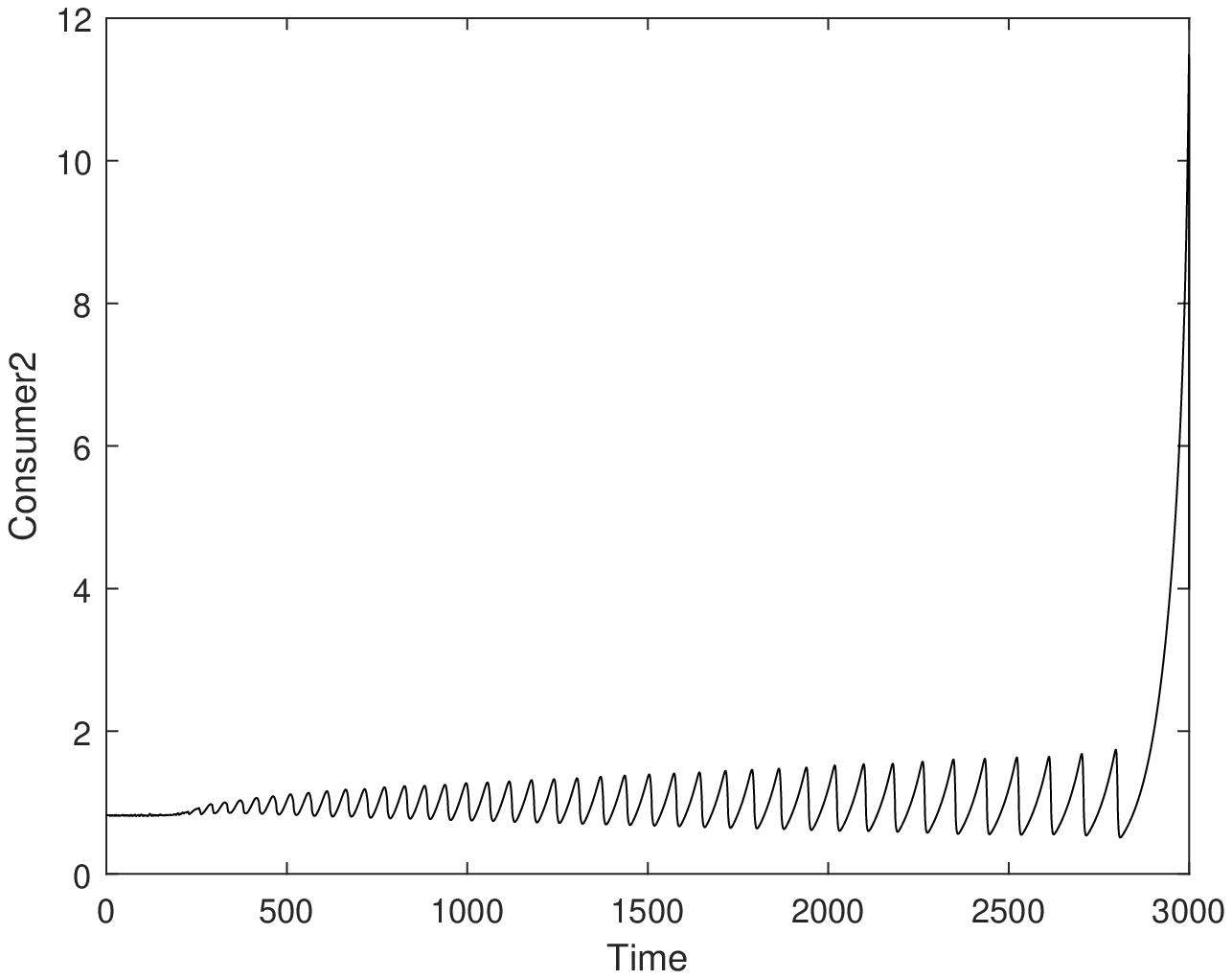}
		\caption{Consumer 2 time series}
	\end{subfigure}
	\caption{Temporal dynamics of system~\eqref{eq6} for the Hopf bifurcation case 3, $\lambda_3>0, \sigma>0$. (a,b,c) before Hopf bifurcation. (d,e,f) after Hopf bifurcation.}
	\label{fig8}
\end{figure}

\begin{figure}[h!]
	\centering
	\begin{subfigure}[b]{0.3\linewidth}
		\includegraphics[width=\linewidth]{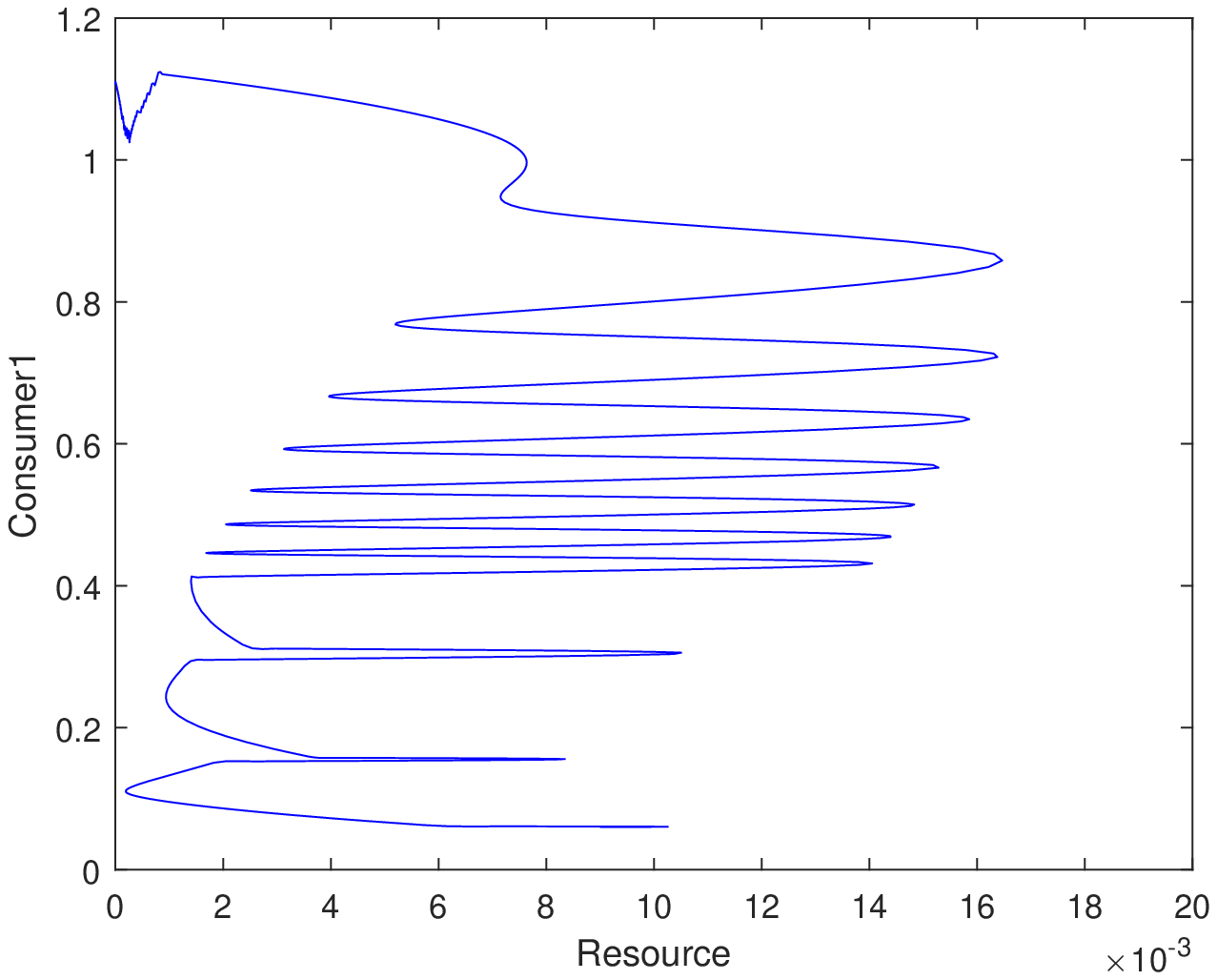}
		\caption{Resource-Consumer 1}
	\end{subfigure}
	\begin{subfigure}[b]{0.3\linewidth}
		\includegraphics[width=\linewidth]{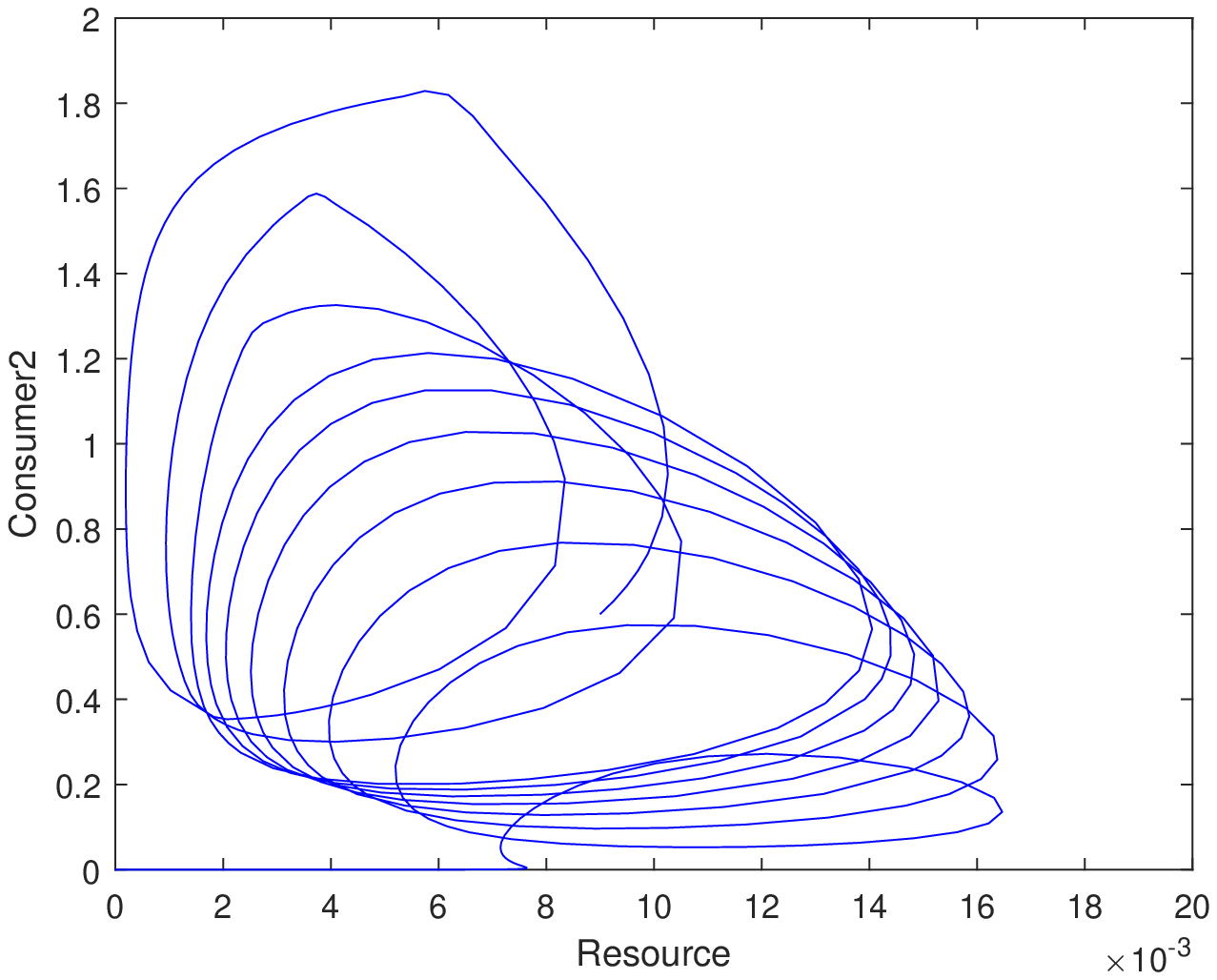}
		\caption{Resource- Consumer 2}
	\end{subfigure}
	\begin{subfigure}[b]{0.3\linewidth}
		\includegraphics[width=\linewidth]{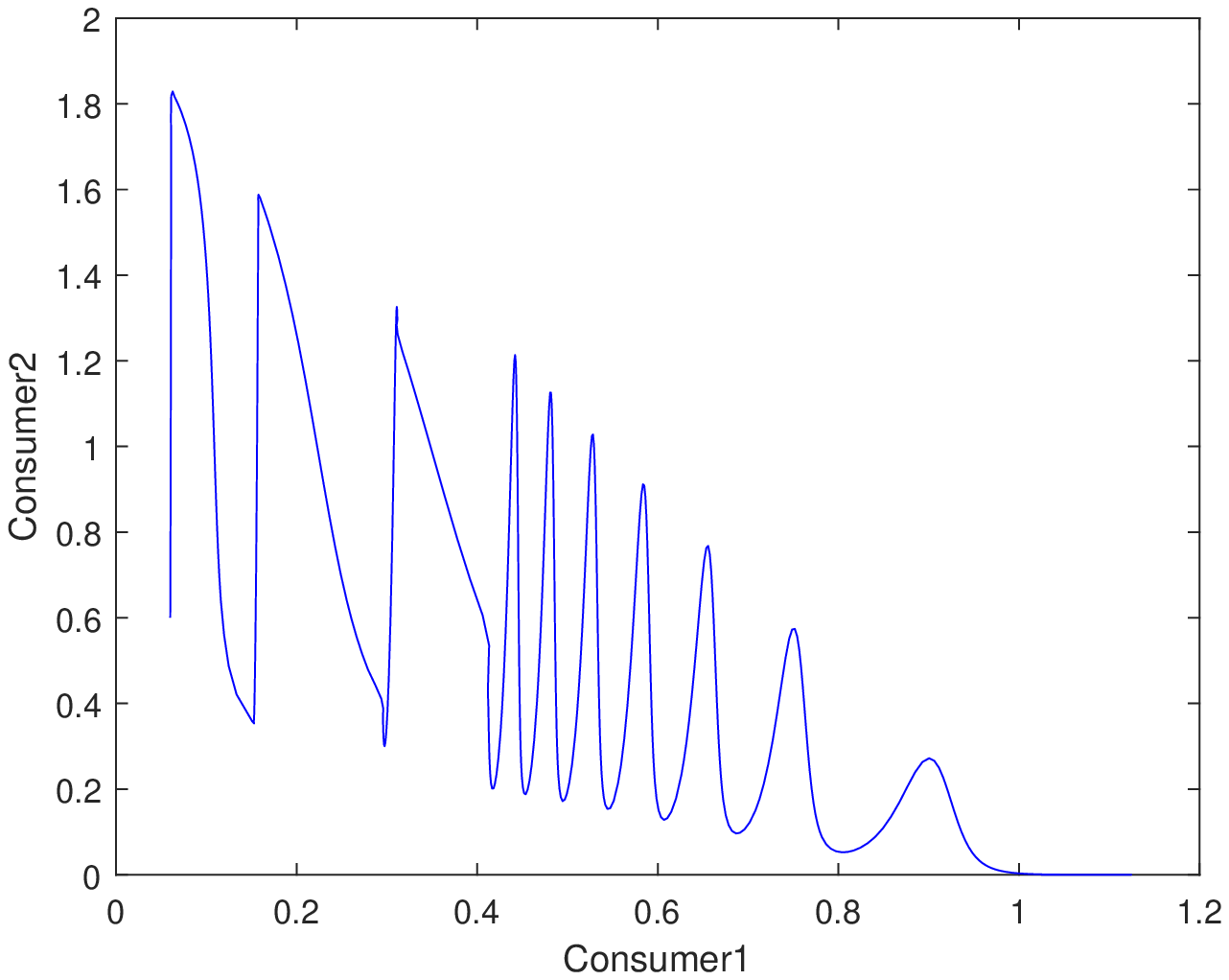}
		\caption{Consumer 1- Consumer 2}
	\end{subfigure}
	\begin{subfigure}[b]{0.3\linewidth}
		\includegraphics[width=\linewidth]{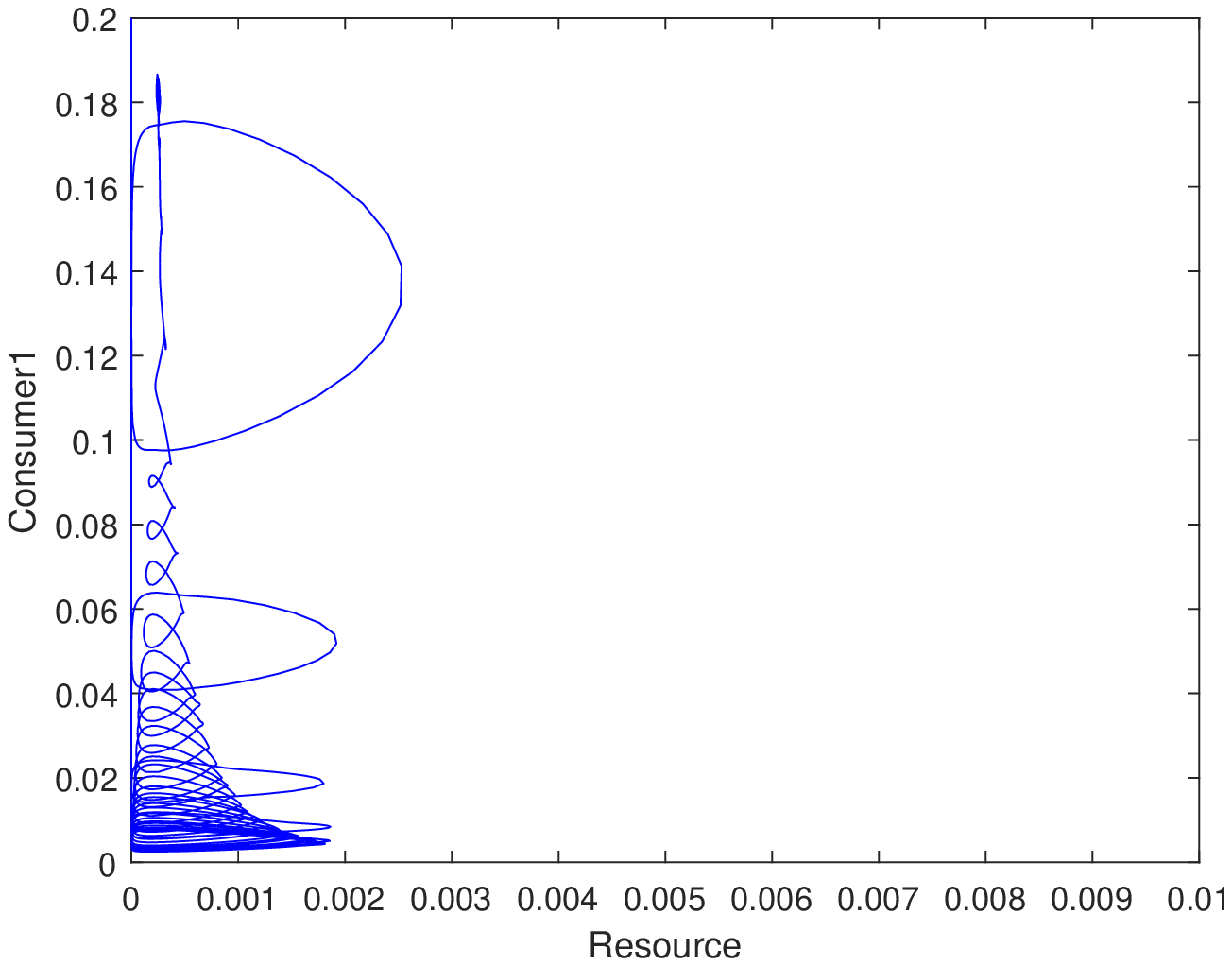}
		\caption{Resource- Consumer 1}
	\end{subfigure}
	\begin{subfigure}[b]{0.3\linewidth}
		\includegraphics[width=\linewidth]{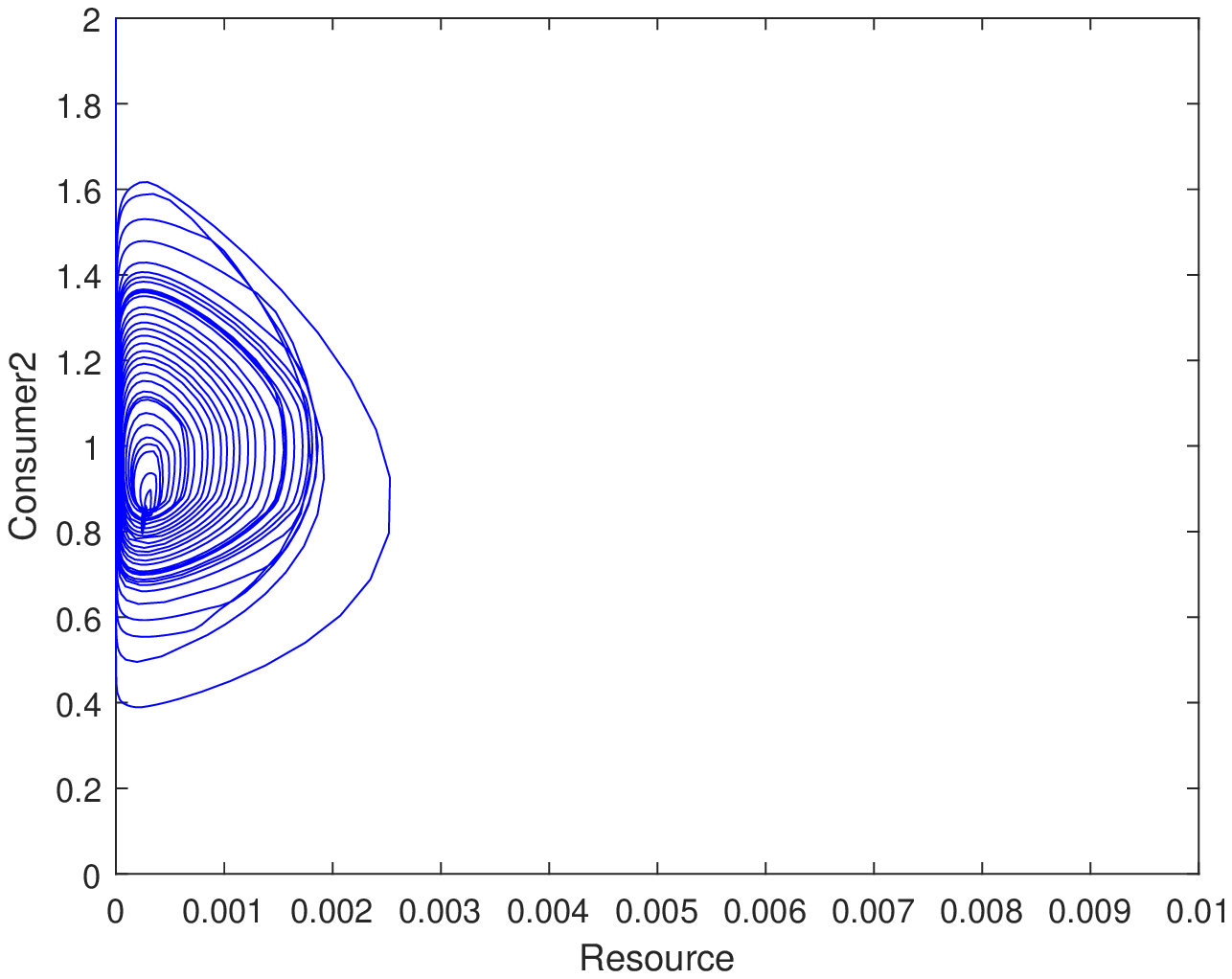}
		\caption{Resource- Consumer 2}
	\end{subfigure}
	\begin{subfigure}[b]{0.3\linewidth}
		\includegraphics[width=\linewidth]{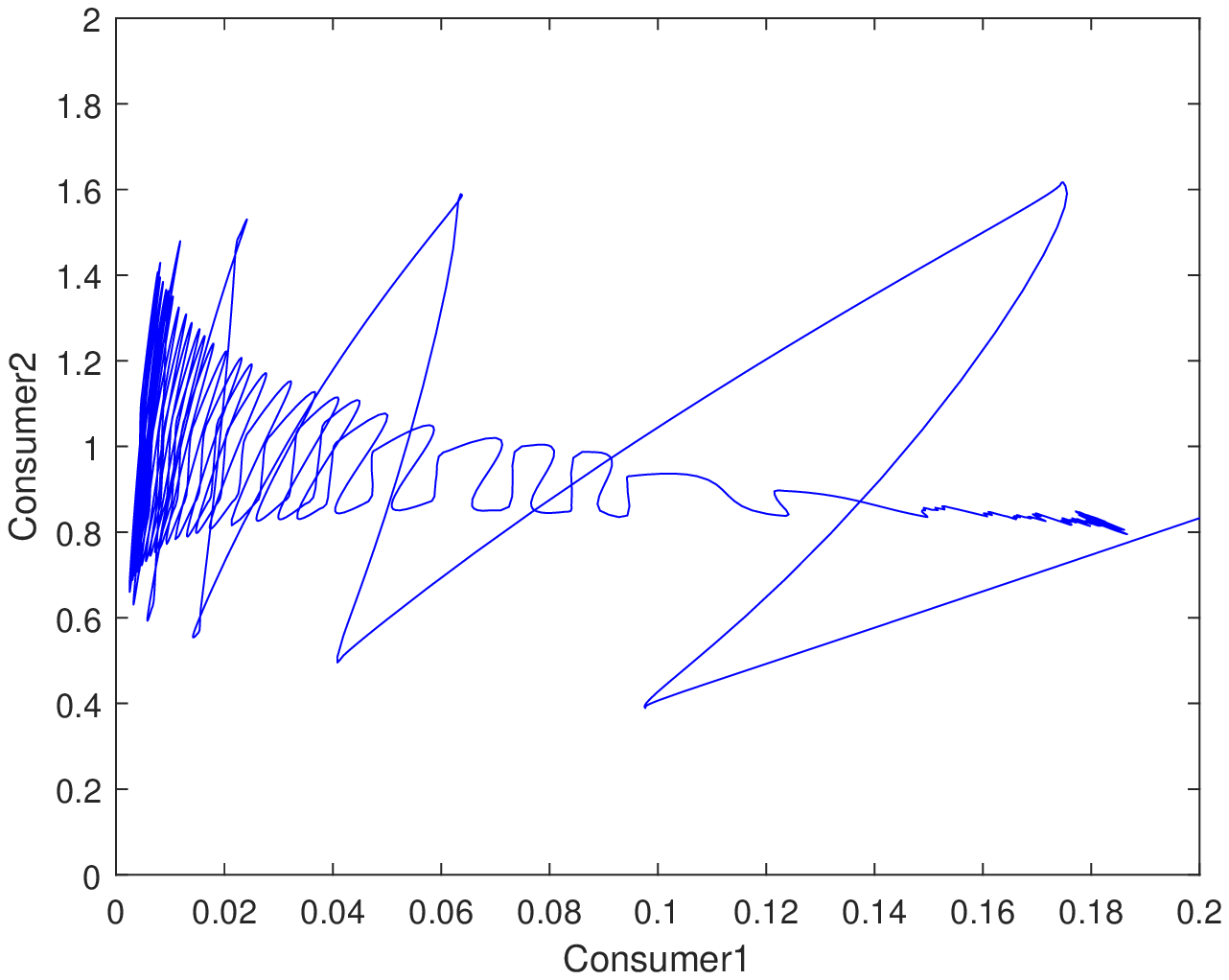}
		\caption{Consumer 1- Consumer 2}
	\end{subfigure}
	\caption{Two dimensional phase diagram of the system~\eqref{eq6} subcritical Hopf bifurcation, case 3, with $\lambda_3>0$ and $\sigma>0$. (a,b,c) before Hopf bifurcation. (d,e,f) after Hopf bifurcation.}
	\label{fig9}
\end{figure}

\begin{figure}[h!]
	\centering
	\begin{subfigure}[b]{0.45\linewidth}
		\includegraphics[width=\linewidth]{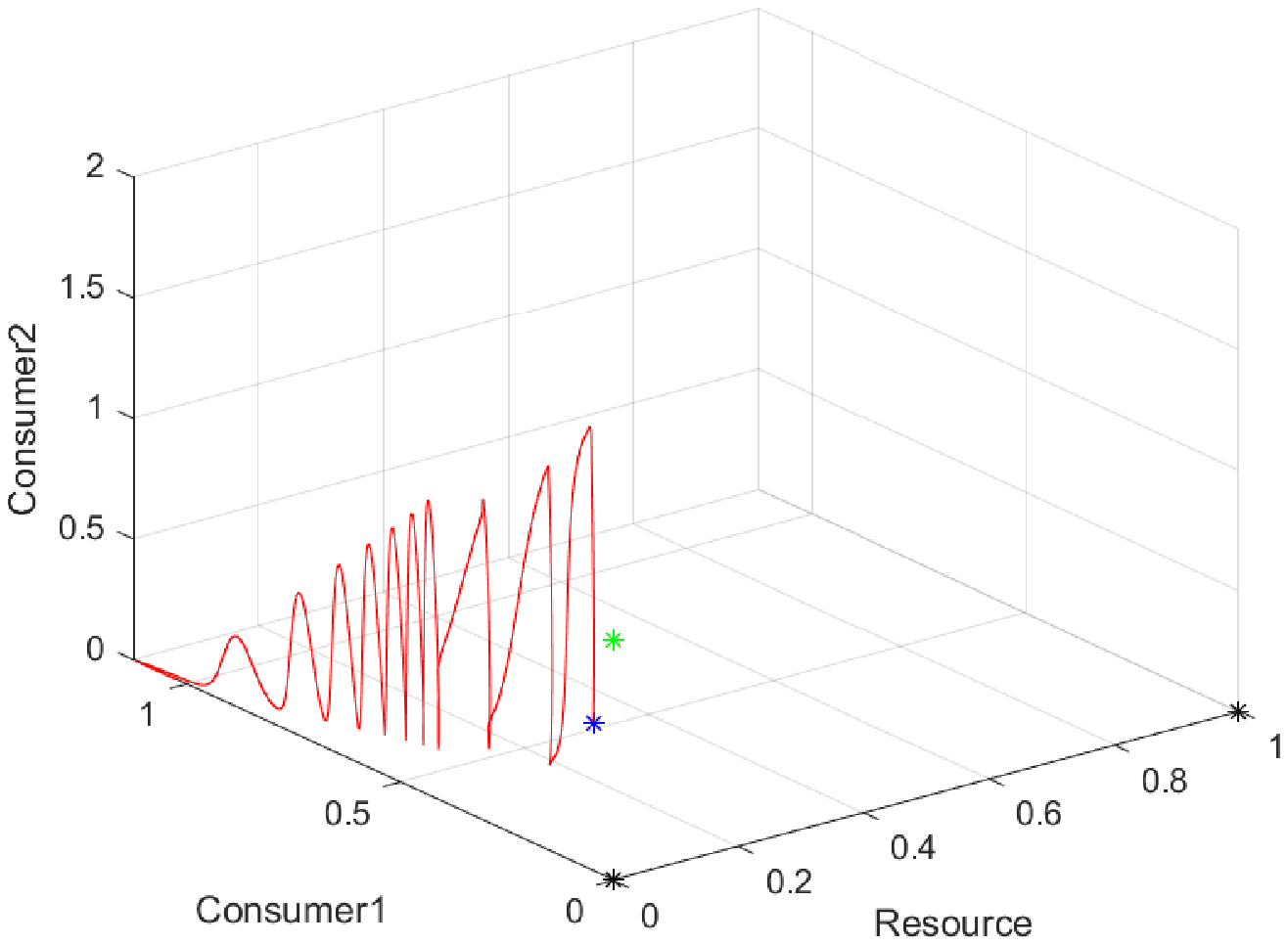}
		\caption{Before Hopf bifurcation}
	\end{subfigure}
	\begin{subfigure}[b]{0.45\linewidth}
		\includegraphics[width=\linewidth]{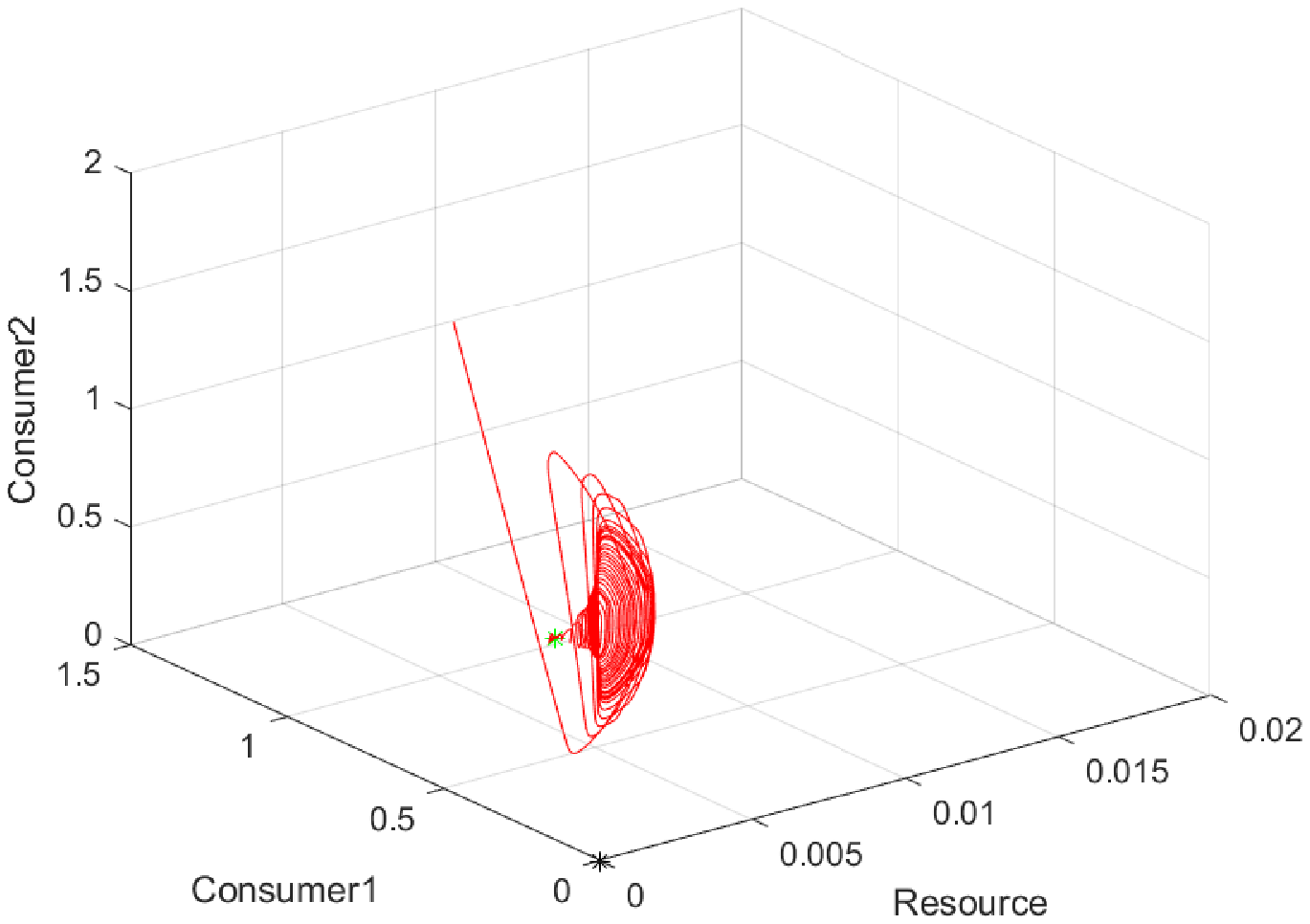}
		\caption{After Hopf bifurcation}
	\end{subfigure}
	\caption{Three dimensional phase diagram of system~\eqref{eq6} for subcritical Hopf bifurcation, case 3, with $\lambda_3>0$ and  $\sigma>0$ before (a) and after (b) bifurcation. Blue, green and black stars are initial values, positive, and boundary equilibrium, respectively.}
	\label{fig10}
\end{figure}

Third case of Hopf bifurcation leading to both positive $\lambda_3$ and $\sigma$ is showed in figures~\ref{fig8} to~\ref{fig10}. Figure~\ref{fig8} (a,b,c) shows temporal dynamics of system~\eqref{eq6} for competitive coefficients $\alpha =0.3369 \times 10^{-8}$ and $\beta =0.94$, undergoing positive first Lyapunov coefficient, $\sigma = 3.46 \times 10^{-6}$ and positive real root $\lambda_3=3.35 \times 10^{-9}$ and before subcritical Hopf bifurcation ($\mu=-7.816 \times 10^{-12}$). Initial values for this simulation is $(0.9,0.06,0.6)$, and the $\mathrm{FRR}$' threshold value is $0.0001$ with $+70$ improvement of attack rate after relaxation. In this case, consumer 1 is more effective in competition than consumer 2 with $\eta_{zy}^-=2.79 \times 10^{8}$. Symmetric parameter values are $a=b=1.01112$, $c=d=1.0075$, and $\mu = \nu = 0.0000003$ in this simulation. Figure~\ref{fig8} (d,e,f) shows temporal dynamics of system~\eqref{eq6} for competitive coefficients $\alpha =0.009$ and $\beta =0.041$, undergoing positive first Lyapunov coefficient, $\sigma = 0.059$ and positive real root $\lambda_3=0.0074$ after subcritical Hopf bifurcation with $\mu=+0.0035$.  Initial values for this simulation is $(0.06,0.4,1)$, and the $\mathrm{FRR}$' threshold value is $0.0001$ with $+60$ improvement of attack rate after relaxation. In this case, again consumer 1 is more effective in competition than consumer 2  with $\eta_{zy}^+=4.65$. Symmetric parameter values are $a=b=2.0143$, $c=d=60.11$, and $\mu = \nu = 0.015$ in this simulation. Figure~\ref{fig8} shows that resource and consumers move towards unstable limit cycles before Hopf bifurcation. After bifurcation, resource and consumers exhibit spirals. For this kind of dynamics change from unstable limit cycles to spirals, bifurcation leads to $\eta_{zy}^-/\eta_{zy}^+\approx 6 \times 10^{7}$ times decreased competitive asymmetricity.

Figure~\ref{fig9} (a,b,c) and (d,e,f) show the unstable limit cycles and spirals, respectively, in subcritical Hopf bifurcation.

Figure~\ref{fig10} shows the three dimensional phase portraits of the subcritical Hopf bifurcation of the system~\eqref{eq6}.

\subsubsection{Hopf case 4, $\lambda_3<0, \sigma>0$}
\begin{figure}[h!]
	\centering
	\begin{subfigure}[b]{0.3\linewidth}
		\includegraphics[width=\linewidth]{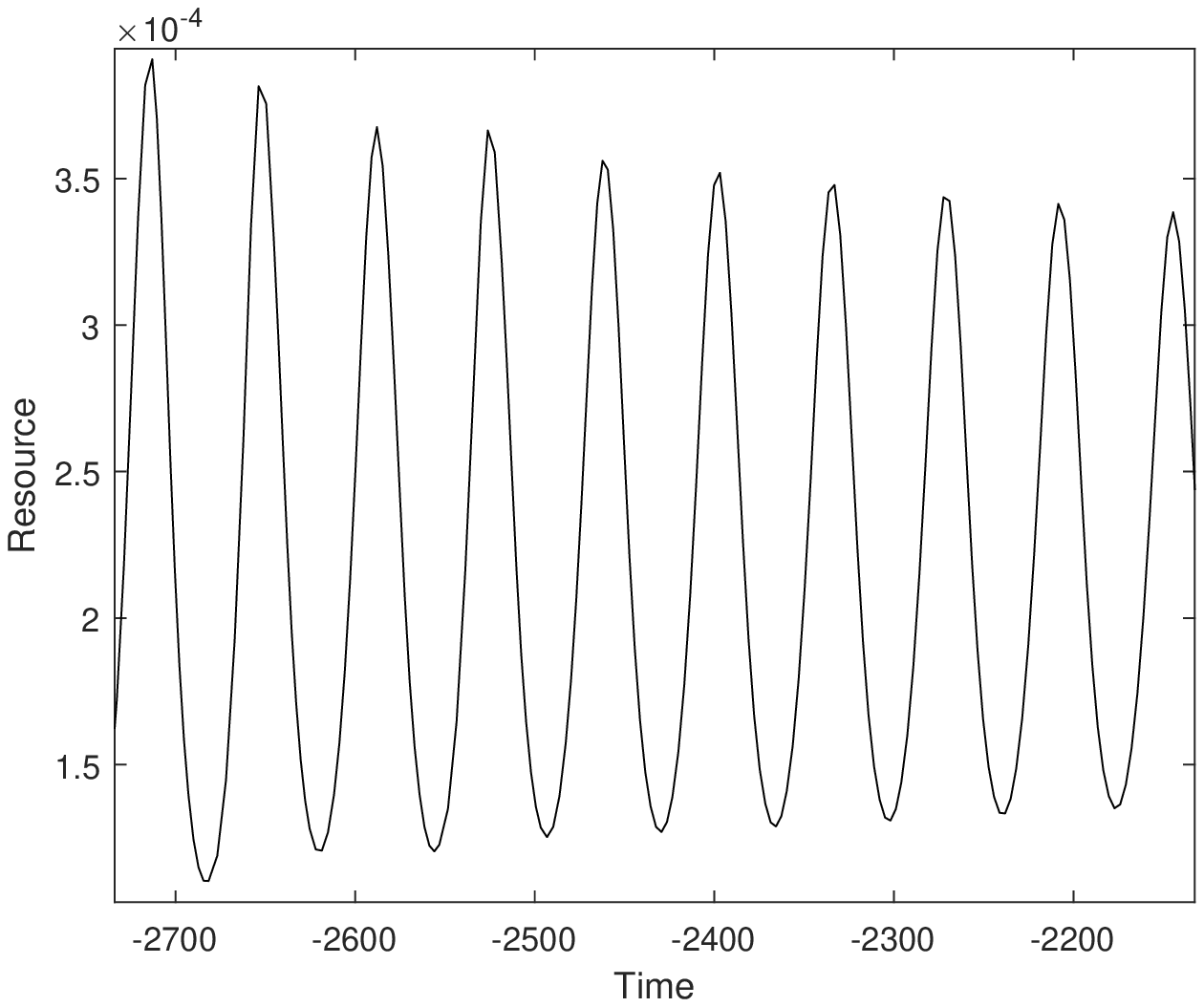}
		\caption{Resource success}
	\end{subfigure}
	\begin{subfigure}[b]{0.3\linewidth}
		\includegraphics[width=\linewidth]{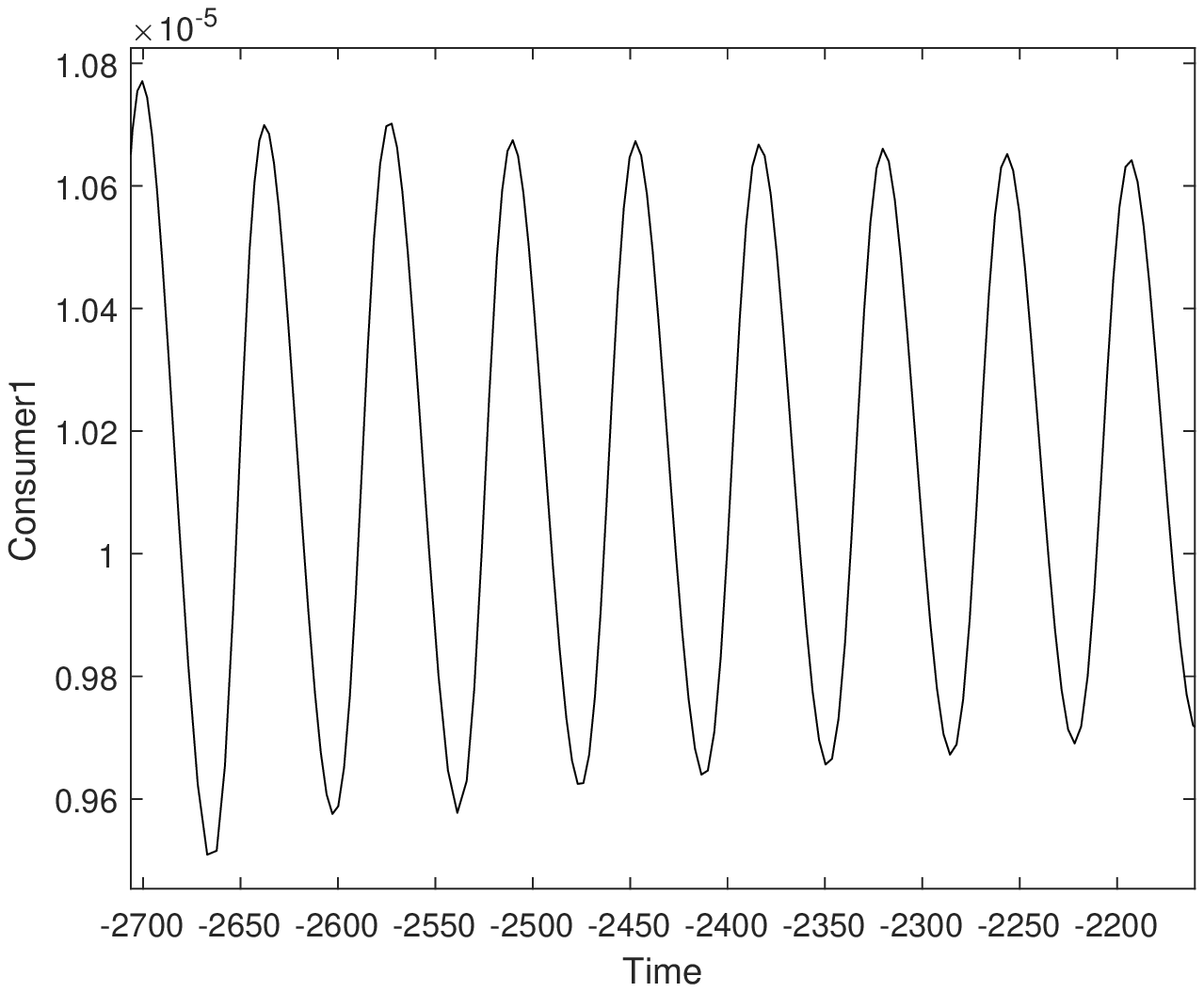}
		\caption{Consumer 1 success}
	\end{subfigure}
	\begin{subfigure}[b]{0.3\linewidth}
		\includegraphics[width=\linewidth]{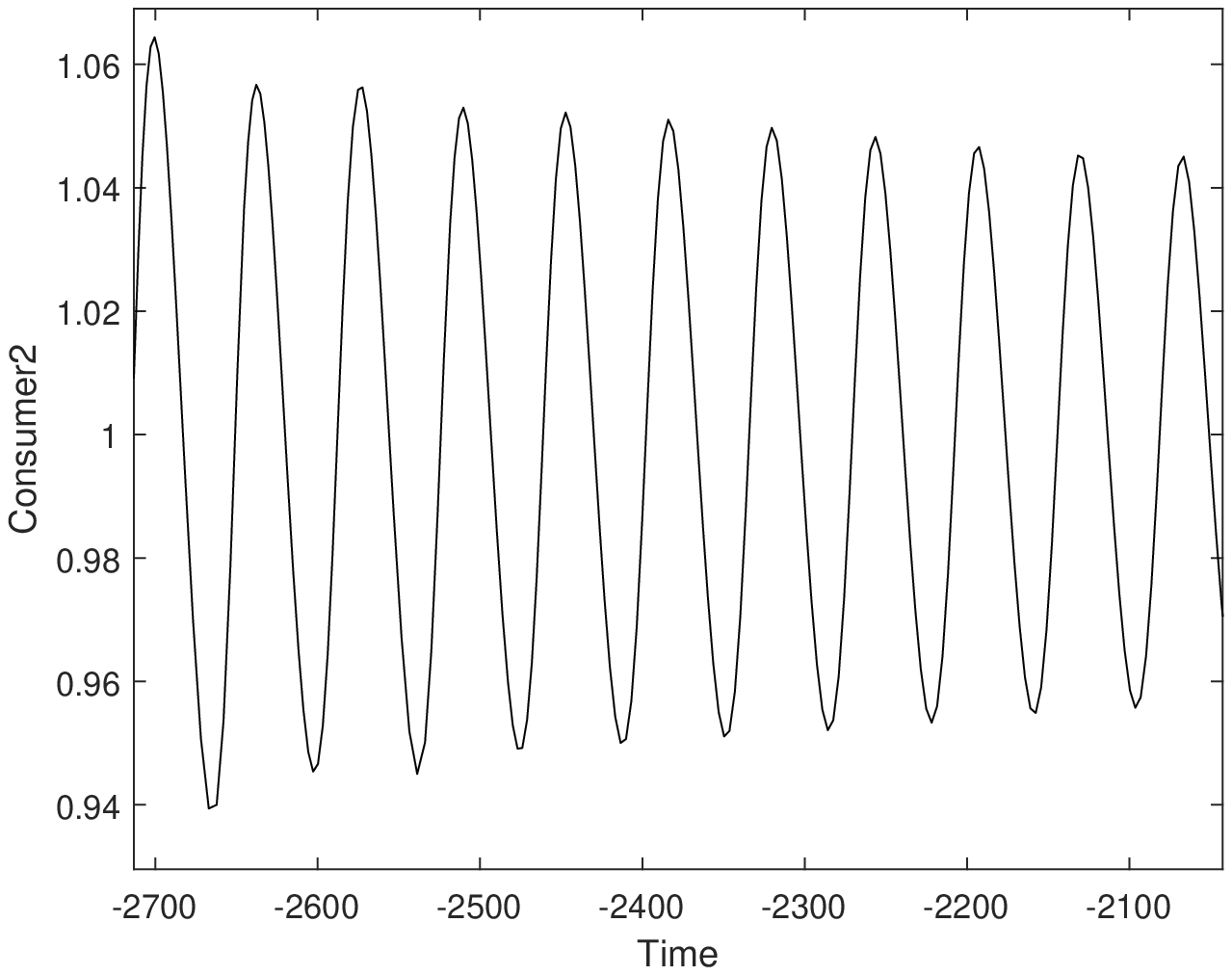}
		\caption{Consumer 2 success}
	\end{subfigure}
	\begin{subfigure}[b]{0.3\linewidth}
		\includegraphics[width=\linewidth]{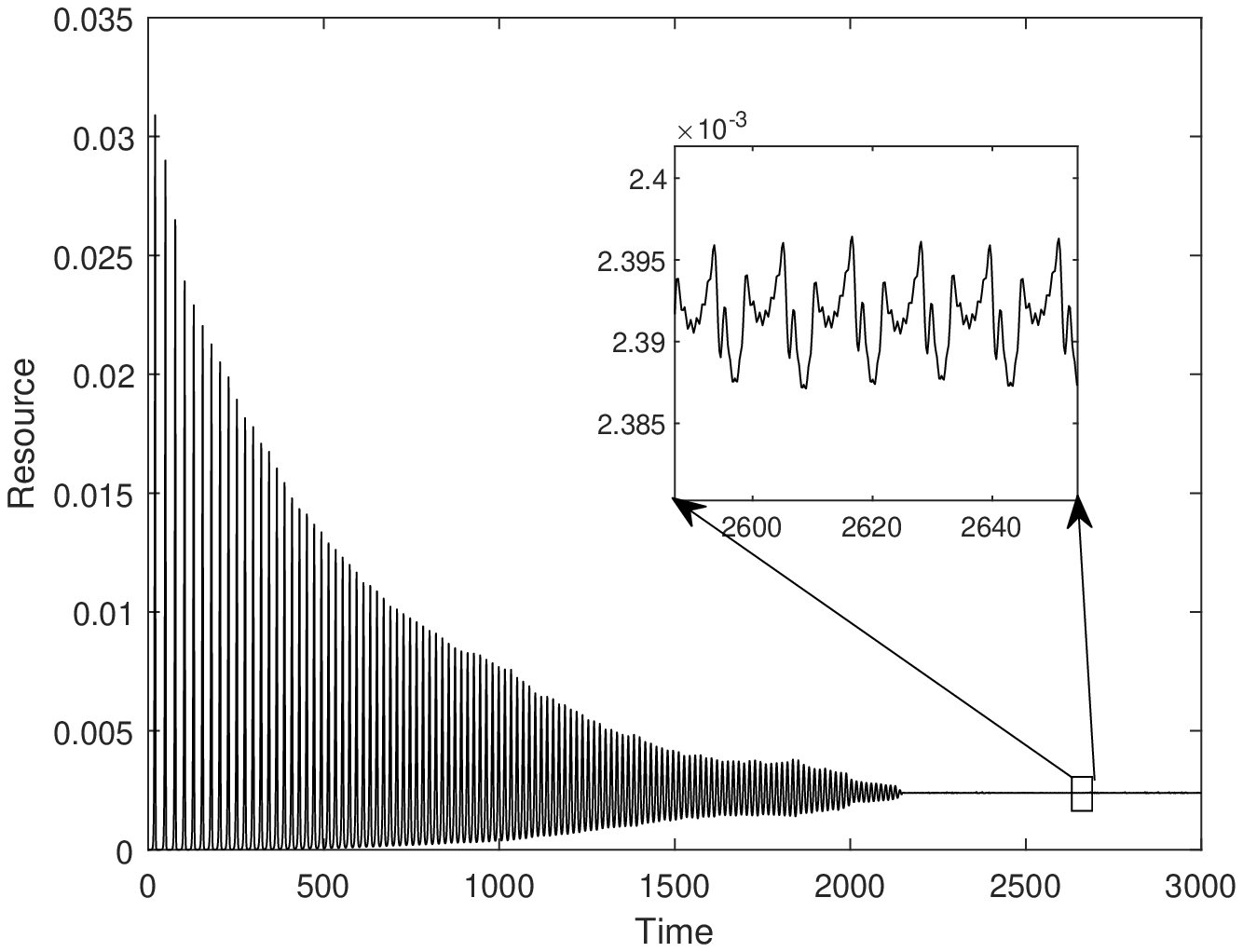}
		\caption{Resource success}
	\end{subfigure}
	\begin{subfigure}[b]{0.3\linewidth}
		\includegraphics[width=\linewidth]{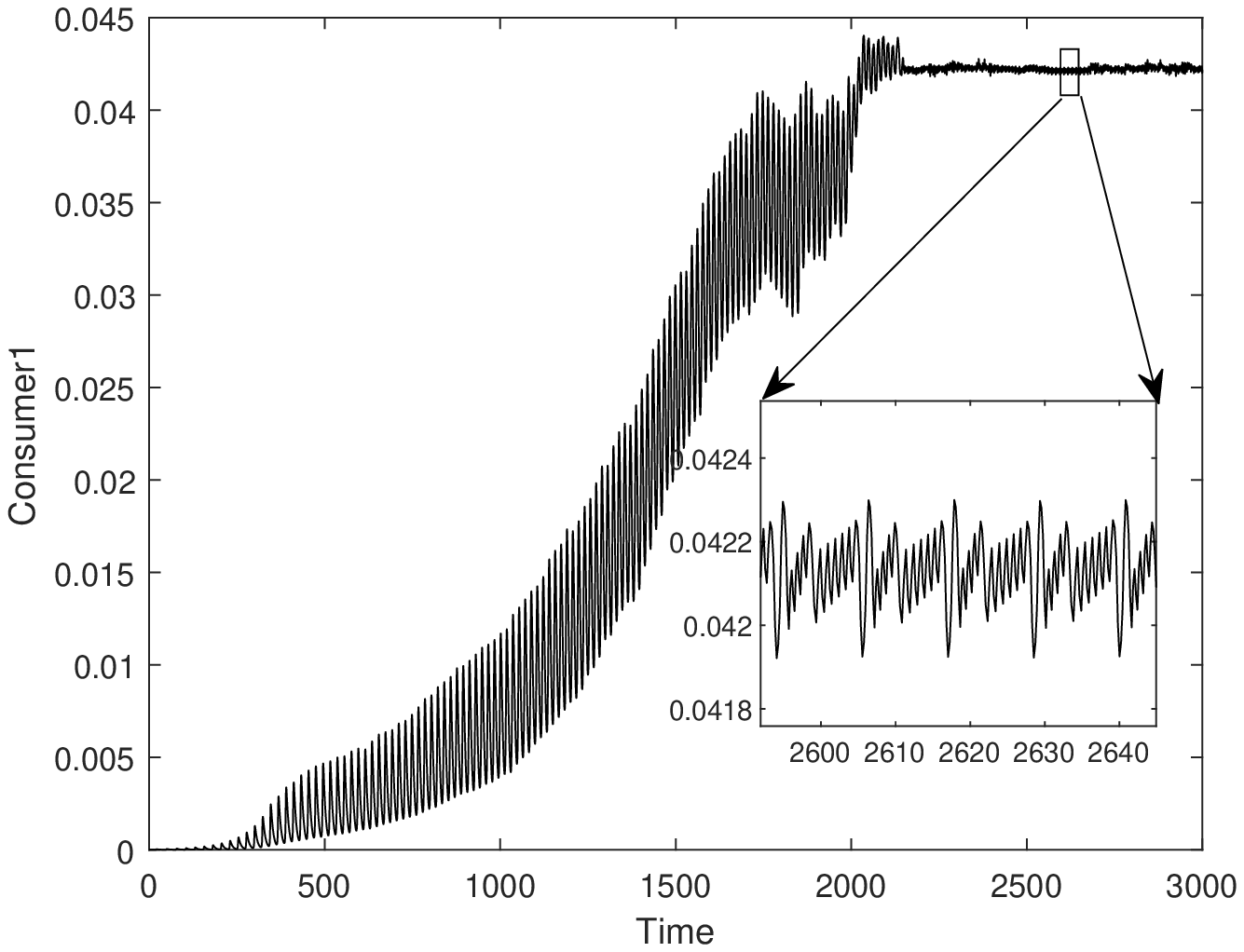}
		\caption{Consumer 1 success}
	\end{subfigure}
	\begin{subfigure}[b]{0.3\linewidth}
		\includegraphics[width=\linewidth]{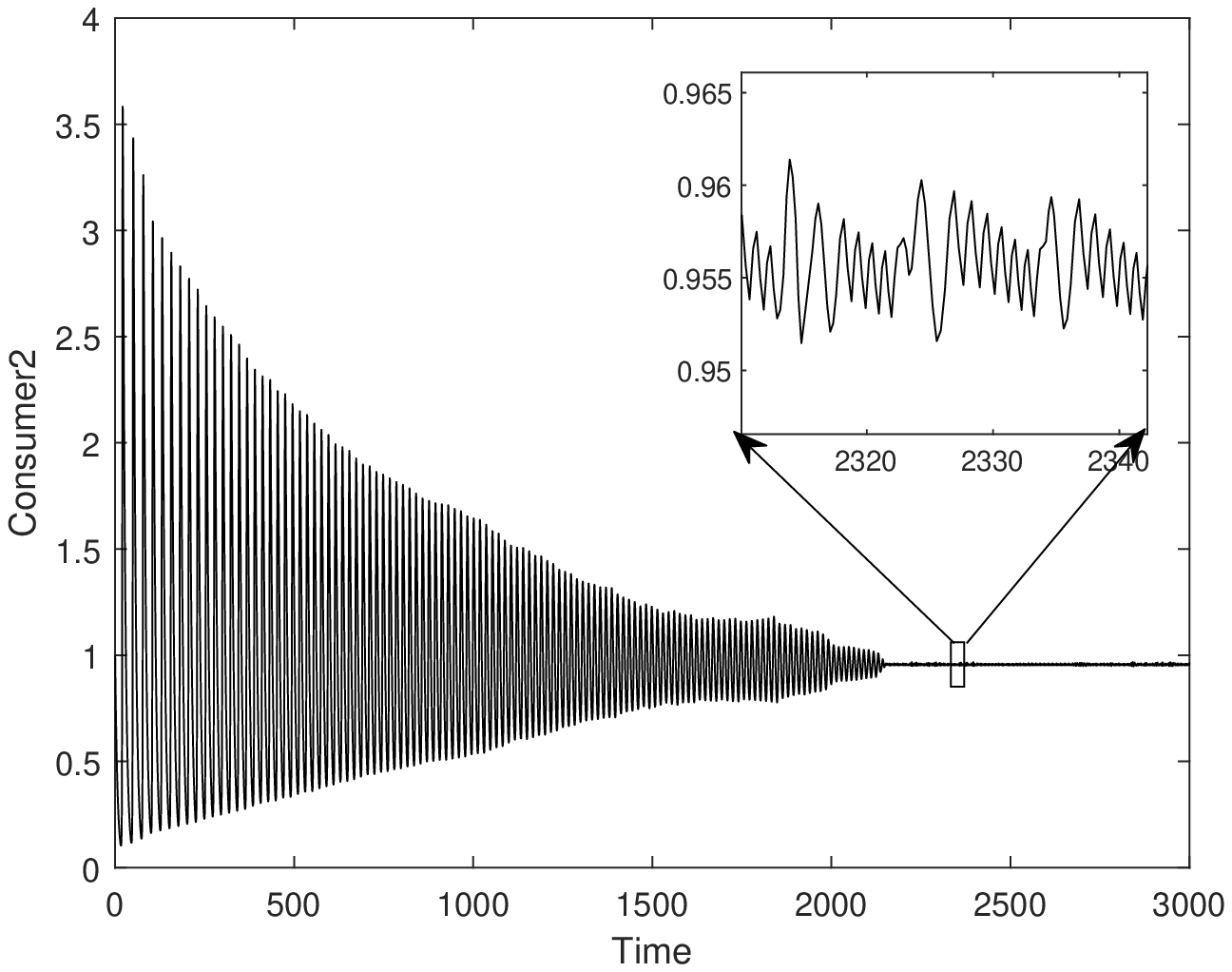}
		\caption{Consumer 2 success}
	\end{subfigure}
	\caption{Temporal dynamics of the system \eqref{eq6} for the Hopf bifurcation, case 4, with $\lambda_3<0$ and $\sigma>0$. (a,b,c) before Hopf bifurcation and (d,e,f) after Hopf bifurcation.}
	\label{fig11}
\end{figure}

\begin{figure}[h!]
	\centering
	\begin{subfigure}[b]{0.3\linewidth}
		\includegraphics[width=\linewidth]{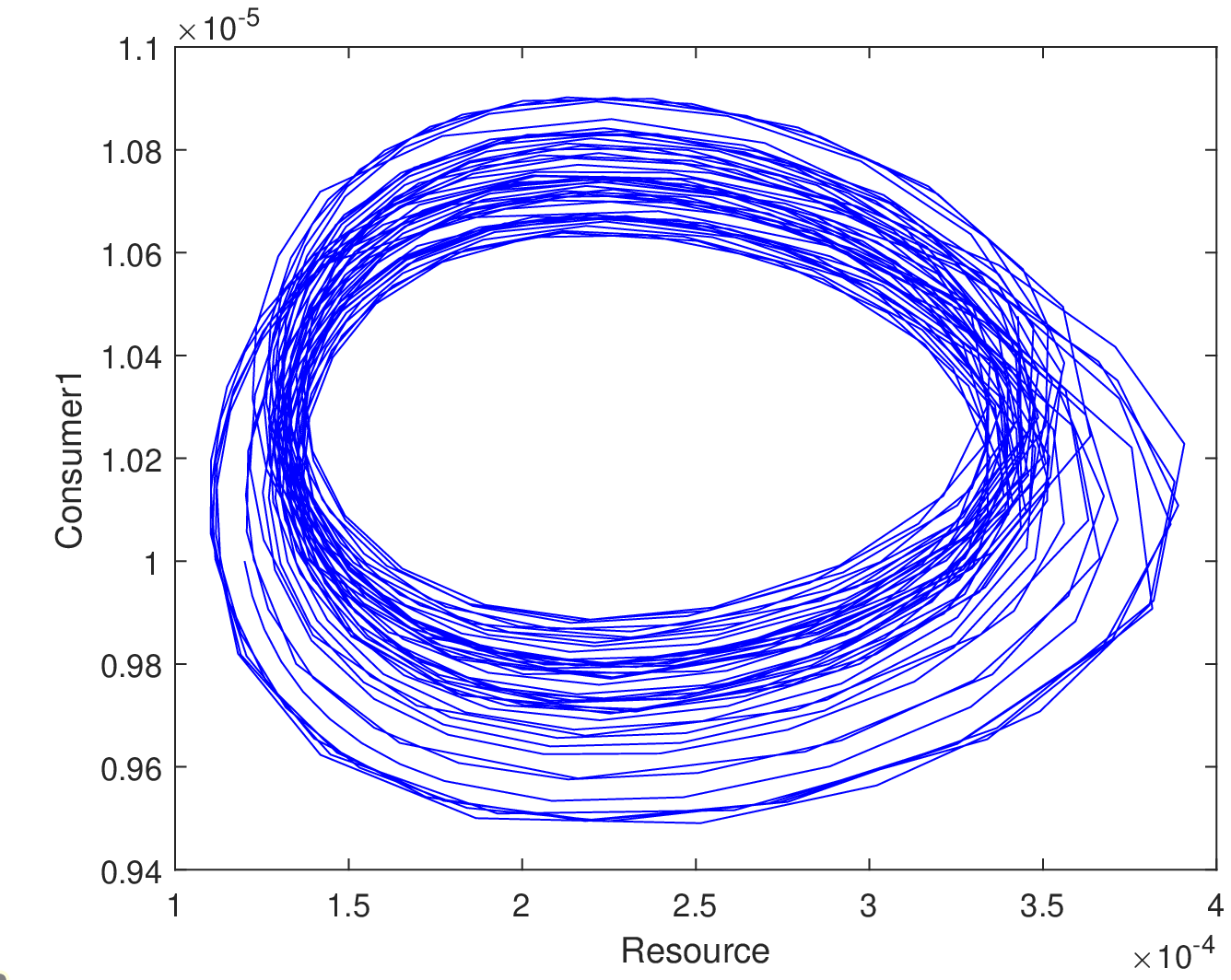}
		\caption{Resource- Consumer 1}
	\end{subfigure}
	\begin{subfigure}[b]{0.3\linewidth}
		\includegraphics[width=\linewidth]{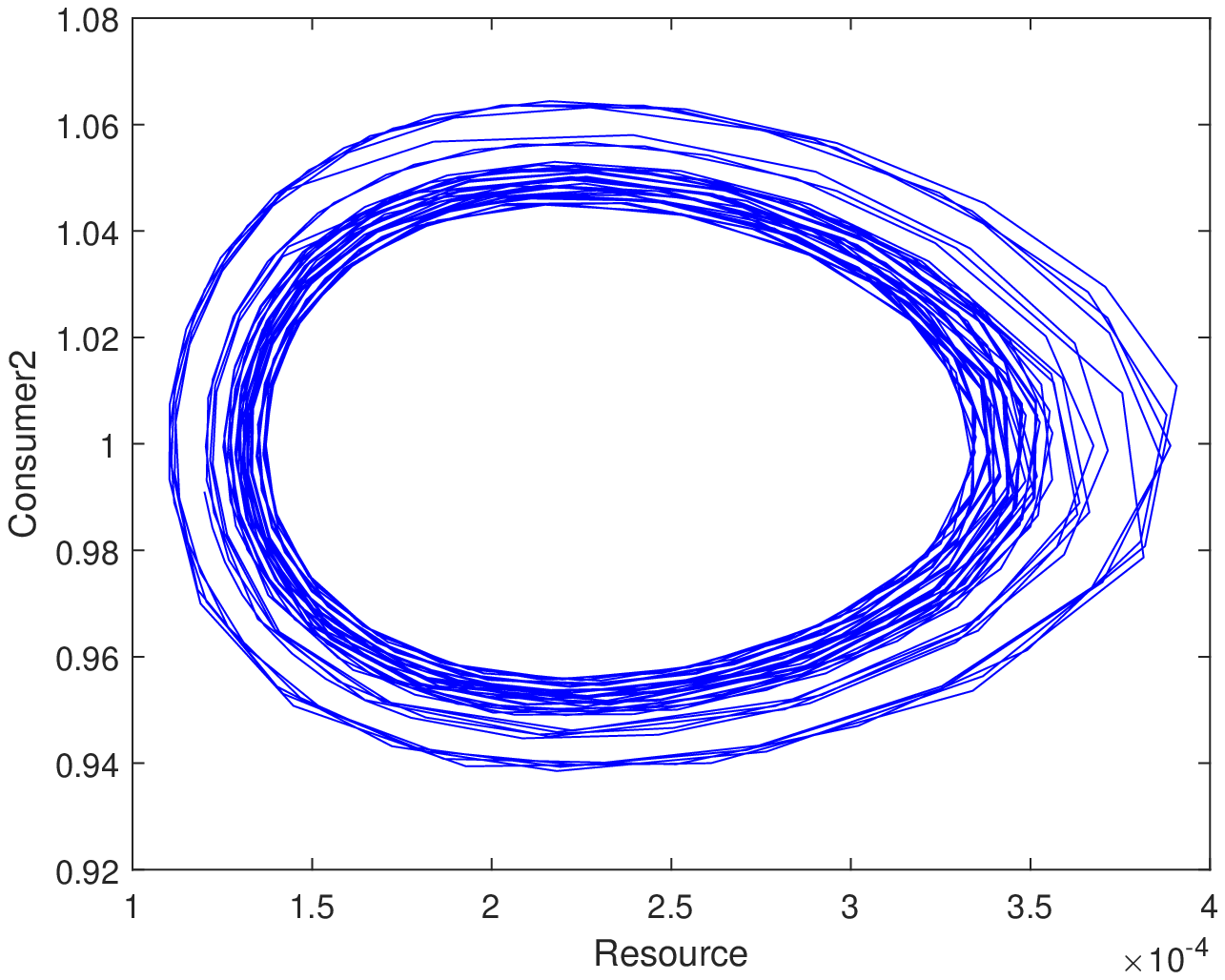}
		\caption{Resource- Consumer 2}
	\end{subfigure}
	\begin{subfigure}[b]{0.3\linewidth}
		\includegraphics[width=\linewidth]{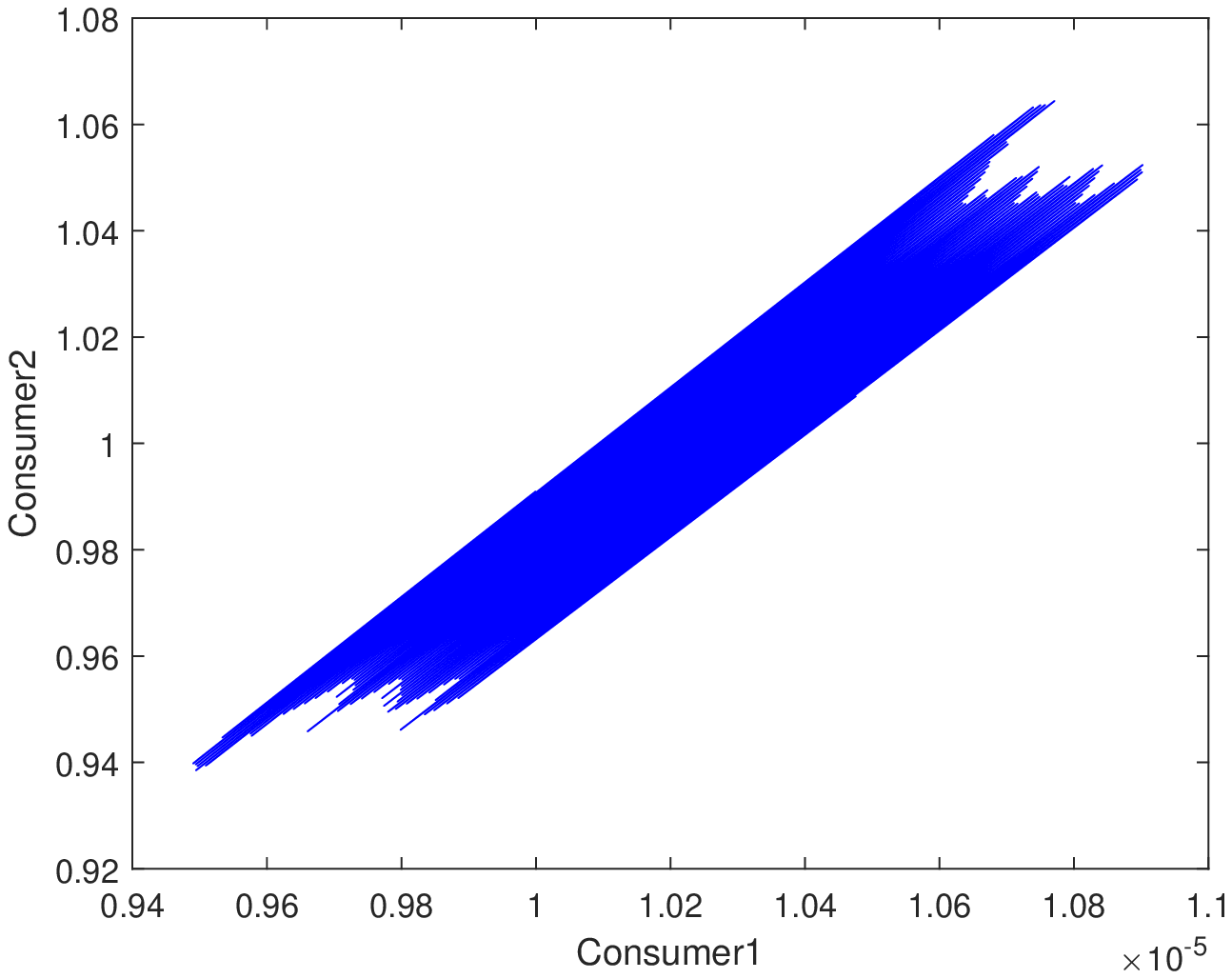}
		\caption{Consumer 1- Consumer 2}
	\end{subfigure}
	\begin{subfigure}[b]{0.3\linewidth}
		\includegraphics[width=\linewidth]{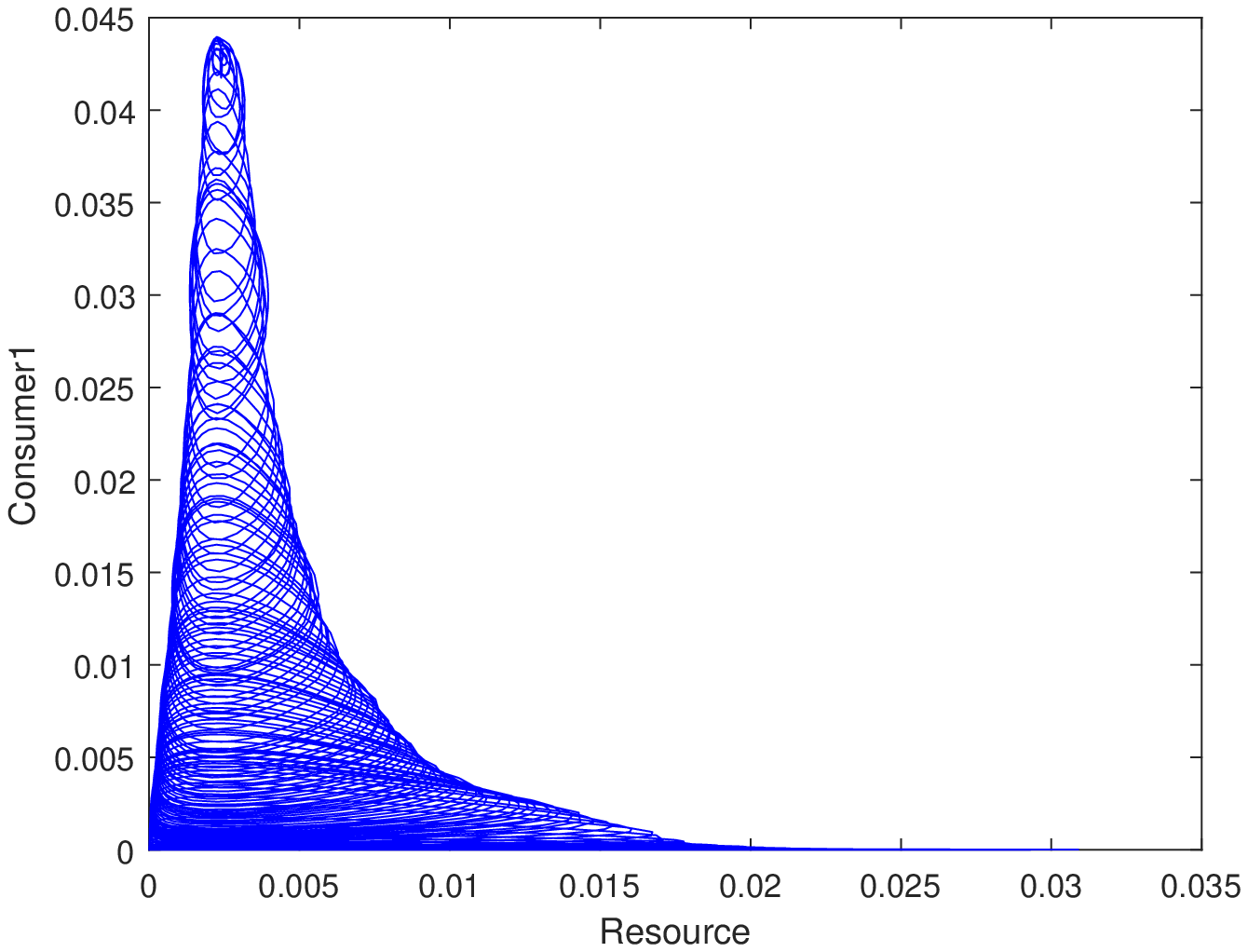}
		\caption{Resource- Consumer 1}
	\end{subfigure}
	\begin{subfigure}[b]{0.3\linewidth}
		\includegraphics[width=\linewidth]{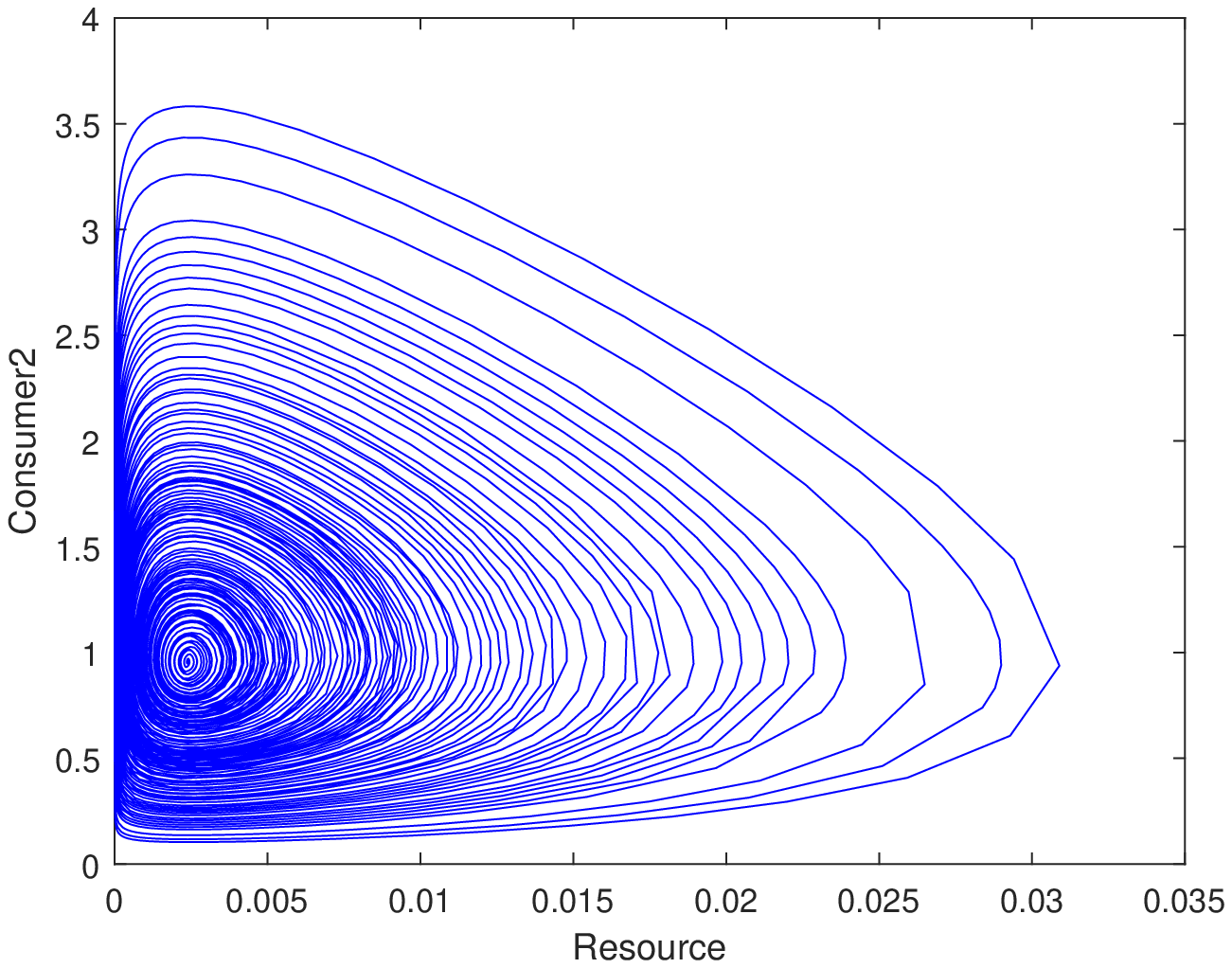}
		\caption{Resource- Consumer 2}
	\end{subfigure}
	\begin{subfigure}[b]{0.3\linewidth}
		\includegraphics[width=\linewidth]{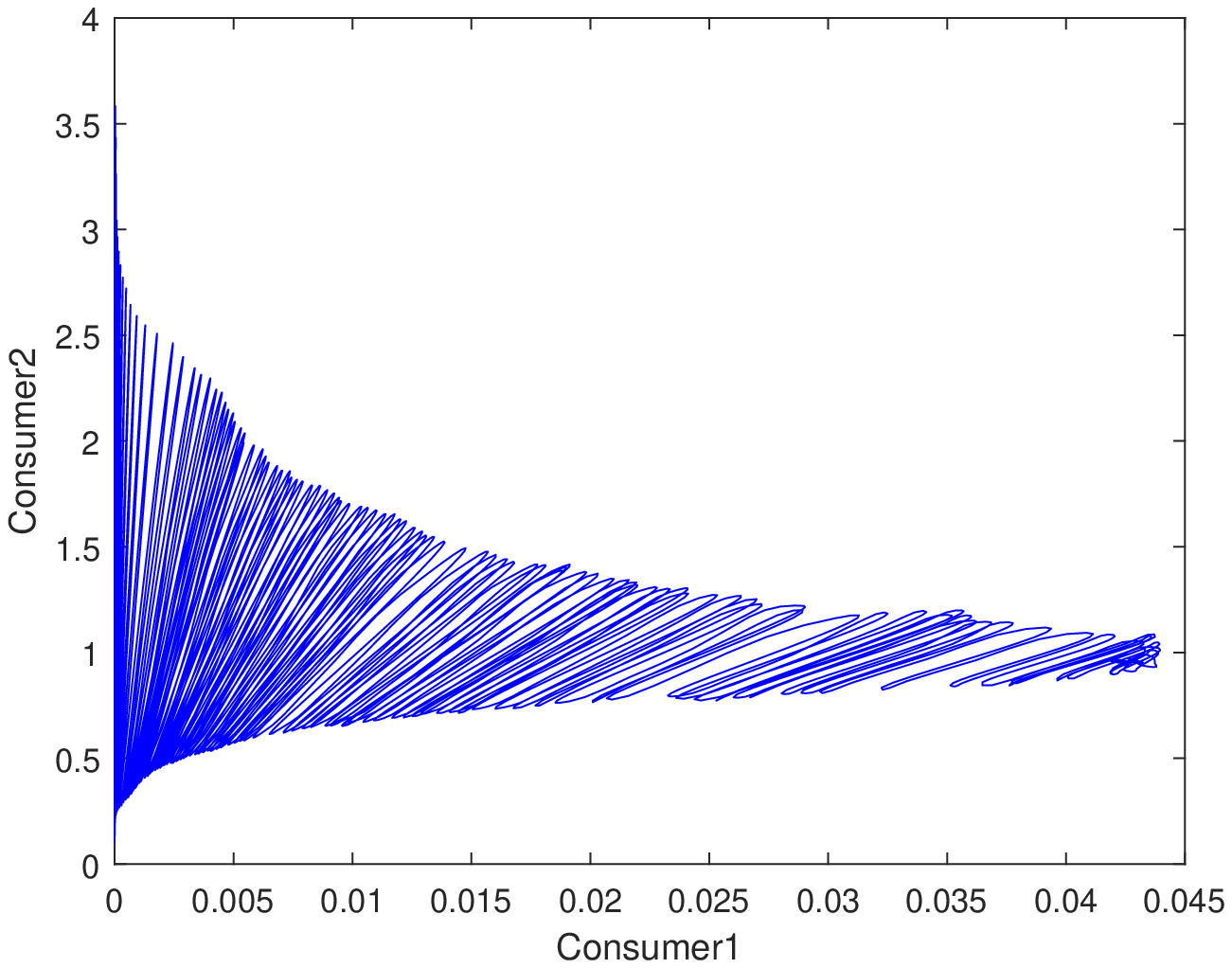}
		\caption{Consumer 1- Consumer 2}
	\end{subfigure}
	\caption{Two dimensional phase diagram of the system~\eqref{eq6} for the Hopf bifurcation case 4, $\lambda_3<0, \sigma>0$. (a,b,c) before Hopf bifurcation. (d,e,f) after Hopf bifurcation.}
	\label{fig12}
\end{figure}

\begin{figure}[h!]
	\centering
	\begin{subfigure}[b]{0.4\linewidth}
		\includegraphics[width=\linewidth]{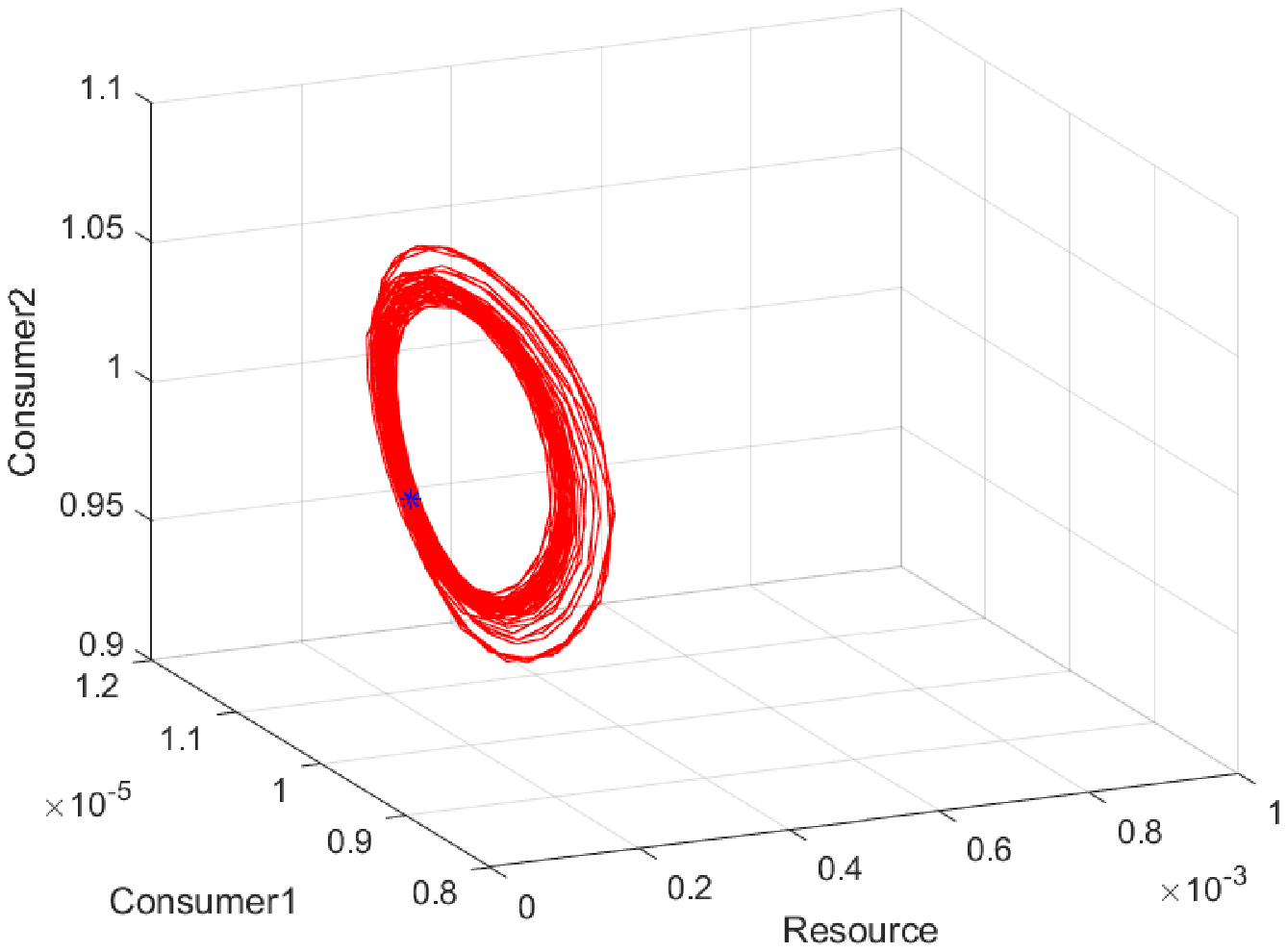}
		\caption{Before Hopf bifurcation}
	\end{subfigure}
	\begin{subfigure}[b]{0.4\linewidth}
		\includegraphics[width=\linewidth]{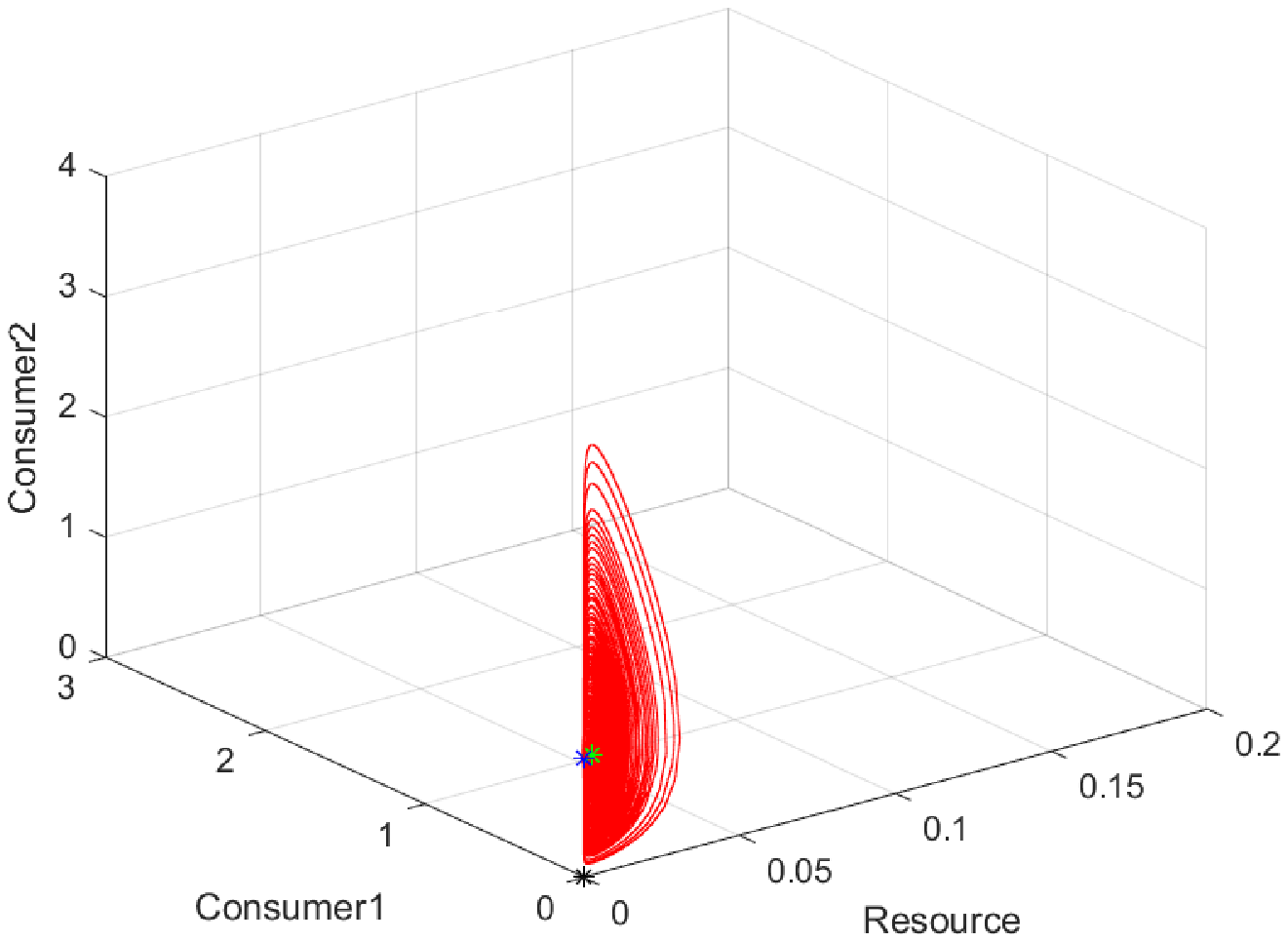}
		\caption{After Hopf bifurcation}
	\end{subfigure}
	\caption{Three dimensional phase diagram of the system~\eqref{eq6} for the Hopf bifurcation case 4 with $\lambda_3<0$ and $\sigma>0$. (a) shows the dynamics before Hopf bifurcation, and (b) shows the dynamics after Hopf bifurcation. Blue, green and black stars are initial values, positive, and boundary equilibrium, respectively.}
	\label{fig13}
\end{figure}

The fourth case of the Hopf bifurcation with  $\lambda_3<0$ and $\sigma>0$ is presented in figures~\ref{fig11} to~\ref{fig13}. Figure~\ref{fig11} (a,b,c) shows temporal dynamics of the system~\eqref{eq6} before subcritical Hopf bifurcation with $\mu=-0.000135$, competitive coefficients $\alpha =0.009\times 10^{-12}$, $\beta =0.91$, and the positive first Lyapunov coefficient $\sigma =0.000725$. Here, the third real root is negative, $\lambda_3=-0.0001635$. The initial point for this simulation is $(0.0012,0.00001,1)$, and the $\mathrm{FRR}$' threshold value is $0.01$ with $+0.1$ improvement of attack rate after relaxation. In this case, consumer 1 is more effective in competition than consumer 2 with $\eta_{zy}^-=1.01 \times 10^{14}$. In this simulation, Symmetric parameter values are $a=b=0.02$, $c=d=45.11$, and $\mu = \nu = 0.01$. Figure~\ref{fig11} (d,e,f) shows temporal dynamics of the system~\eqref{eq6} after subcritical Hopf bifurcation  with $\mu=+0.00104$, the competitive coefficients $\alpha =0.0009\times 10^{-10}$, $\beta =0.041$, the positive first Lyapunov coefficient $\sigma =0.00176$. Here, the third real root is negative, $\lambda_3=-1.27 \times 10^{-5}$. The initial point for this simulation is $(0.0000001,0.000005,1)$, and the $\mathrm{FRR}$' threshold value is $0.001$ with $+6$ improvement of attack rate after relaxation. In this case, again consumer 1 is more effective in competition than consumer 2 with $\eta_{zy}^+=4.56 \times 10^{11}$. In this simulation, the symmetric parameter values are $a=b=0.15$, $c=d=60.11$, and $\mu = \nu = 0.15$. Resource and consumers 1 and 2 exhibit unstable relaxation oscillations in their time series before the bifurcation. While after the Hopf bifurcation, stable relaxation-oscillations were dominant dynamics in the time series. For this kind of change in dynamics, from unstable limit cycles to stable relaxation$-$oscillations, bifurcation leads to $\eta_{zy}^-/\eta_{zy}^+\approx 221$ times decreased competitive asymmetricity.

\subsection{The case of zero-Hopf Bifurcation}
We have the possibility of zero-Hopf bifurcation to the system~\eqref{eq9} when $\lambda_{1,2}=\pm i\omega $ and $\lambda_3=0$. First we bring the zero$-$Hopf type equilibrium definition and after some algebraic calculations, we obtain the following lemma.

\begin{definition}
	A zero$-$Hopf type equilibrium for a first$-$order three dimensional autonomous differential system is an isolated equilibrium point of the system for which its linear part has one zero and two purely imaginary eigenvalues \cite{Guckenheimer}
\end{definition}

\begin{lemma}\label{lem2}
  The jacobian matrix $J_1$ of the linearization of the system~\eqref{eq6} at the positive equilibrium $E^*$ has a pair of imaginary eigenvalues $\lambda_{1,2}=\pm i\omega_0$ and a zero eigenvalue  $\lambda_3 =0$ if and only if $\delta=\delta^0\in\Omega^0$ where

  \begin{align}\label{eq102}
    \Omega^0=\{\delta:\; \mathrm{a}(\delta)= 0,\; \mathrm{b}(\delta)>0,\; \mathrm{c}(\delta)=0\}
	\end{align}
 
\end{lemma}
\proof
 We look for the region of parameters for which $P(\lambda)$ has a pair of complex conjugate eigenvalues $\lambda_{1,2}=\pm i\omega_0$ and a real zero eigenvalue $\lambda_3=0$. This setting requires that
\begin{align}
  P(\lambda)=\lambda(\lambda^2+\omega^2)+(\mathrm{b}-\omega^2) \lambda.
\end{align}

By some simple algebra for $P(\lambda)=0$ we have $\lambda_3=0$, and $\lambda_{1,2}=\pm i\omega_0=\pm i\sqrt{\mathrm{b}^0}$ when $\mathrm{a}^0=0),\; \mathrm{b}^0=\mathrm{b}(\delta^0),\; \mathrm{c}^0=0)$, and $\delta^0\in\Omega^0$.
\hfill $\square$\\

We can choose $\phi^0$ such that $\delta^0=\mathrm{P}\phi^0\in\Omega^0$.
Suppose $\delta^0\in\Omega^0$ and $J_1(0)=J_1(\delta^0)$. For the given $\epsilon-$ perturbation $\delta(\epsilon)=\delta^0+O(\epsilon)$, an $\epsilon-$perturbation of matrix $J_1(0)$ is the matrix $J_1^{\epsilon}=J_1\left(\delta^0+O(\epsilon)\right)$. Then $\mathrm{a}(\epsilon)=\mathrm{a}^0+O(\epsilon)=O(\epsilon)$, $\mathrm{b}(\epsilon)=\mathrm{b}^0+O(\epsilon)$, and $\mathrm{c}(\epsilon)=\mathrm{c}^0+O(\epsilon)=O(\epsilon)$. 

\begin{theorem}\label{th2}
  For small enough real number $\epsilon$, with $|\epsilon|<<1$ the jacobian matrix $J_1^{\epsilon}$ of the linearization of the system~\eqref{eq6} at the positive equilibrium $E^*$ has a pair of complex conjugate eigenvalues $\lambda_{1,2}=\mu(\epsilon)\pm i\omega(\epsilon)$ and a real eigenvalue  $\lambda_3(\epsilon) \in \mathbb{R}$. For $\epsilon=0$, we have $\lambda_3(0)=0$, $\mu(0)=0$, $\mu'(0)\neq 0$, and $\omega(0)=\omega_0>0$. 
\end{theorem}
\proof

We look for a perturbation $\delta(\epsilon)=\delta_0+O(\epsilon)$ such that $P(\lambda)$ has a pair of complex conjugate eigenvalues $\lambda_{1,2}=\mu(\epsilon)\pm i\omega(\epsilon)$ and a real eigenvalue $\lambda_3(\epsilon)$.
Based on the continuity of the roots of a polynomial with respect to its coefficients, this setting requires that
\begin{align}
  P(\lambda)=\left(\lambda+2\mu\right)\left(\lambda^2-2\mu\lambda+(\mu^2+\omega^2)\right),
\end{align}
subject to the following conditions for $|\epsilon|<<1$.
	\begin{align}\label{eq102-2}
          & \mathrm{a}=\mathrm{a}^0+O(\epsilon)=O(\epsilon),\; \mathrm{b}=\mathrm{b}^0+O(\epsilon),\; \mathrm{c}=\mathrm{c}^0+O(\epsilon)=O(\epsilon),\nonumber\\
          &   \omega^2(\epsilon)=\mathrm{b}^0+O(\epsilon),\; \mu=O(\epsilon). \nonumber
	\end{align}
These conditions fulfill the requirements. Moreover, we have $\mu'(0)\neq 0$. \hfill $\square$\\

Based on Theorem~\ref{th1}, the positive equilibrium $E^*$ of the system~\eqref{eq6} undergoes a generic zero-Hopf bifurcation at $\delta^0\in\Omega^0$ in equation (\ref{eq102}).
To study the qualitative behavior of the system~\eqref{eq6} near the bifurcation point, we employ the normal form of zero-Hopf bifurcation in~\cite{Guckenheimer}.
There, the near identity change of variable $Y=\xi+W(\phi,\chi)\in R^3$ and the adjoint operator $ad_L:H_k\to H_K$ for $k=2,3$ are employed to simplify the nonlinear vector field $\hat{F}_{\phi^0}(\chi)$.
Here, $H_k$ is the vector space of monomials of order $k$ in $R^3$ and $L=J_1(\delta)Y$. Then, the polar change of coordinates $(y_1,y_2)\to (r,\theta)$, truncation up to the third order, removing the azimuthal term, truncation up the second order, and rescaling $\bar{r}=-\sqrt{c_1c_3}\, r,\; \bar{z}=-c_3\, z$ reduces~
\footnote{The reduction procedure is well explained in~\cite{Guckenheimer} with details. Therefore, we do not duplicate them here.}
the system~\eqref{eq6} around $\delta^0\in\Omega^0$ as follows.
\begin{equation}\label{eq44}
	\begin{cases}
		\dot{\bar{r}}=-\frac{a_1}{c_3}\, \bar{z}\bar{r},\\
		\dot{\bar{z}}=-\frac{c3}{|c3|}\frac{c1}{|c1|} \bar{r}^2-\bar{z}^2.
	\end{cases}
\end{equation}

\begin{figure}[t]
  \begin{center}
\begin{picture}(420,150)(0,0)
  \put(0,20)  {\psfig{file=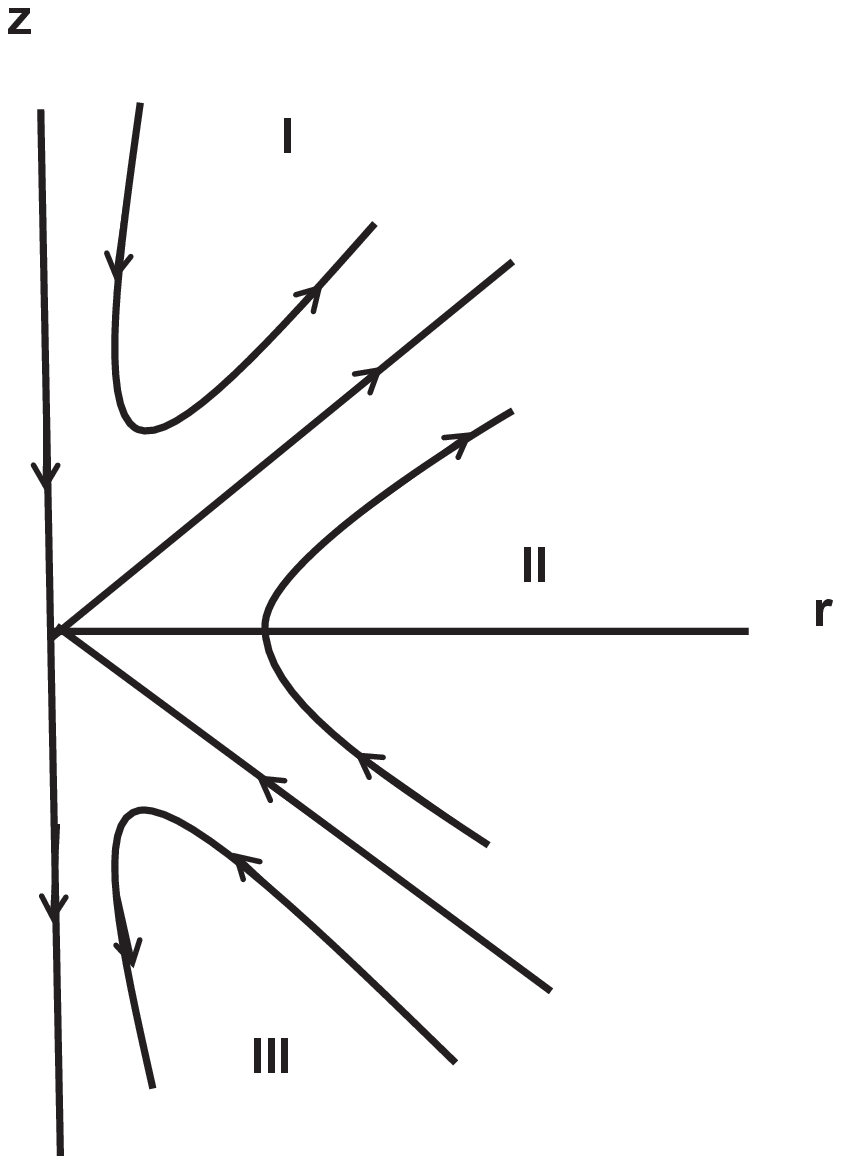, width=3cm}}  \put(30,10) {Form I  }
  \put(90,20) {\psfig{file=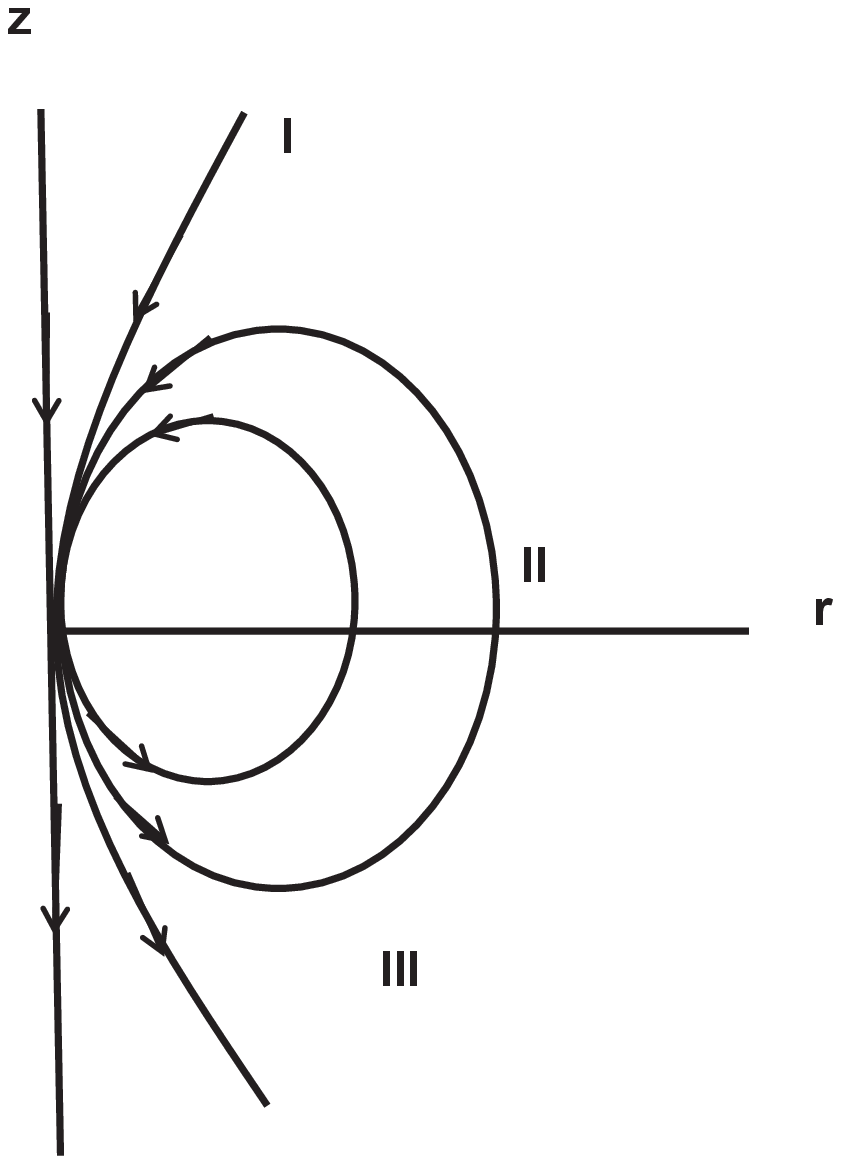,width=3cm}}  \put(115,10){Form IIb}
  \put(180,20){\psfig{file=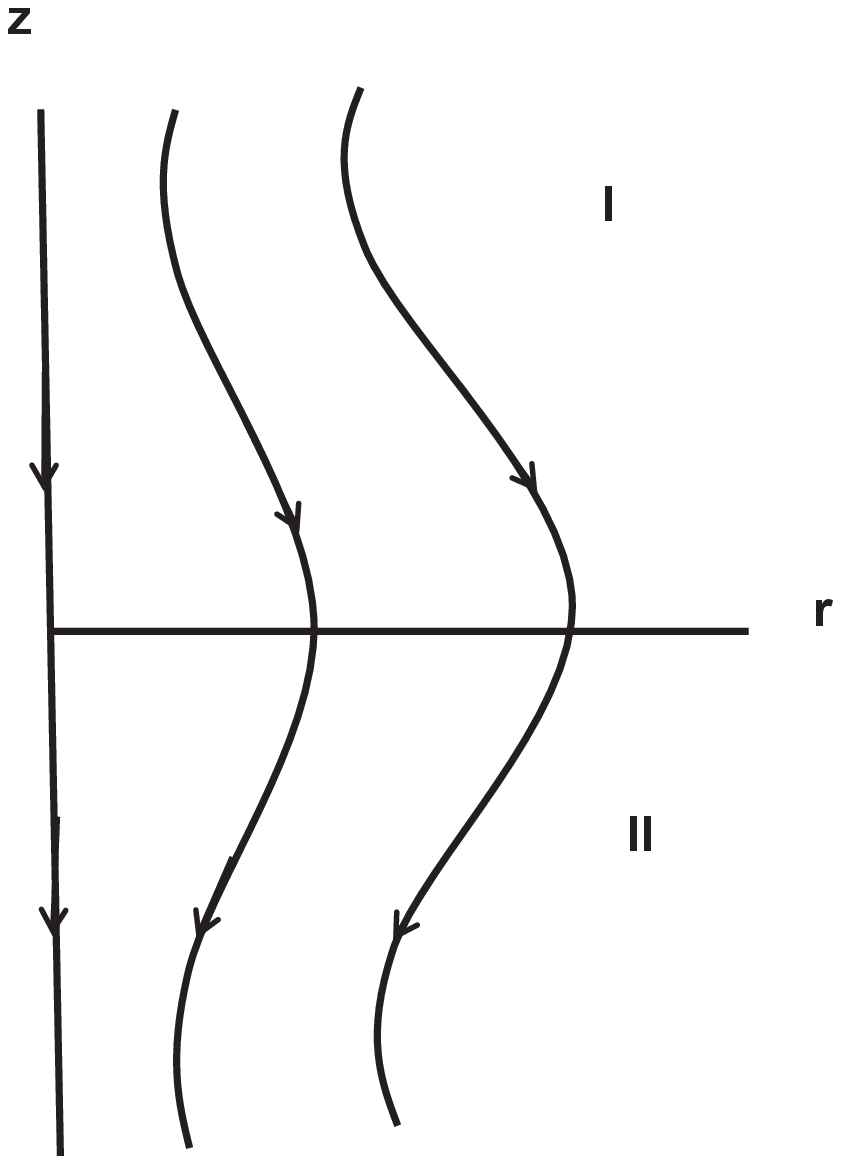, width=3cm}}  \put(205,10){Form III}
  \put(270,20){\psfig{file=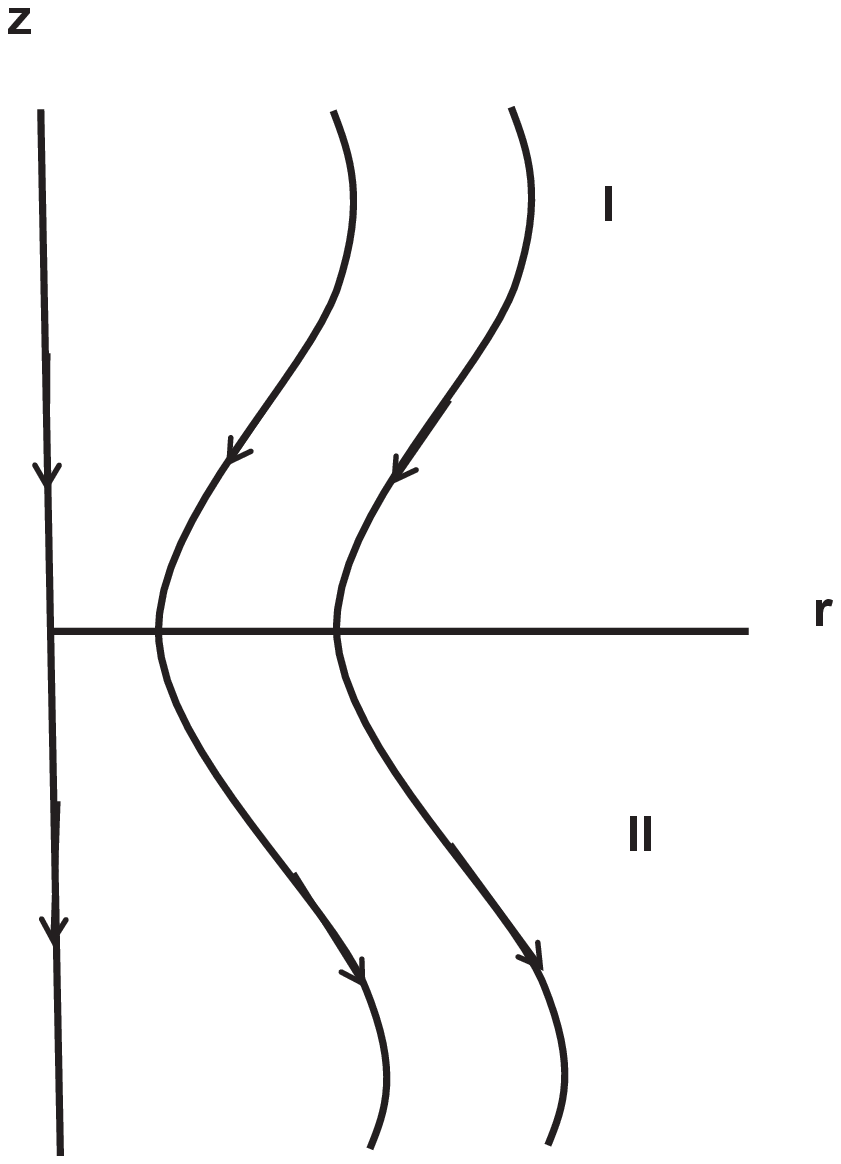,width=3cm}}  \put(295,10){Form IVa}
  \put(360,20){\psfig{file=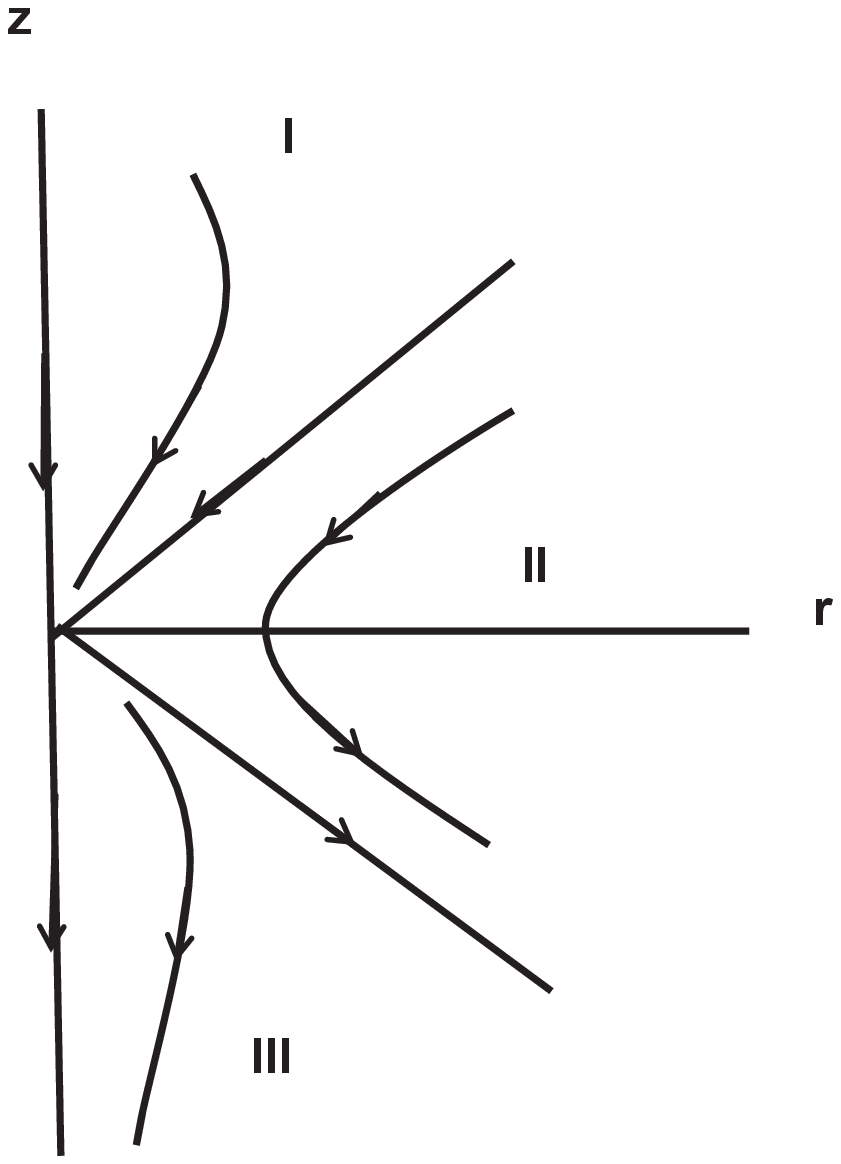,width=3cm}}  \put(380,10){Form IVb}
\end{picture}
\caption{\label{fig14}
  Five types of zero$-$Hopf normal form of the respective catalog presented in ~\cite[Sec. 7.4]{Guckenheimer}.
  Type IIa was not depicted in our simulations after the range of the parameters value checked. Regions I to III are used to simulate the system repective to each parametric domain.}
  \end{center}
\end{figure}

Typical patterns of the zero$-$Hopf bifurcation are depicted in Figure~\ref{fig14}; see~\cite{Guckenheimer} for more details.
By varying coefficients of competition, various initial values, $\mathrm{FRR}$ threshold values, and respective improved attack rates, different phase portraits of the system~\eqref{eq44} around the zero-Hopf bifurcation can be obtained.  

\subsubsection{Type I of zero$-$Hopf normal form catalog, $a>0,b=+1$}
\begin{figure}[h!]
	\centering
	\begin{subfigure}[b]{0.45\linewidth}
		\includegraphics[width=\linewidth]{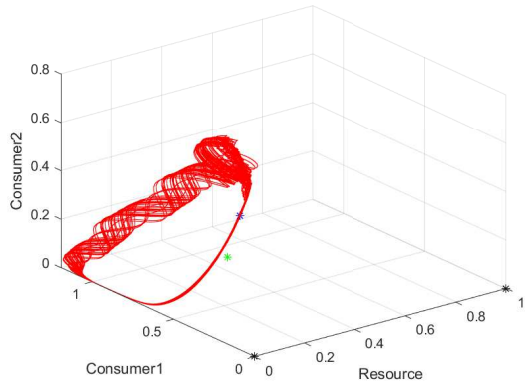}
		\caption{Dynamics in the region I in the form I}
	\end{subfigure}
	\begin{subfigure}[b]{0.45\linewidth}
		\includegraphics[width=\linewidth]{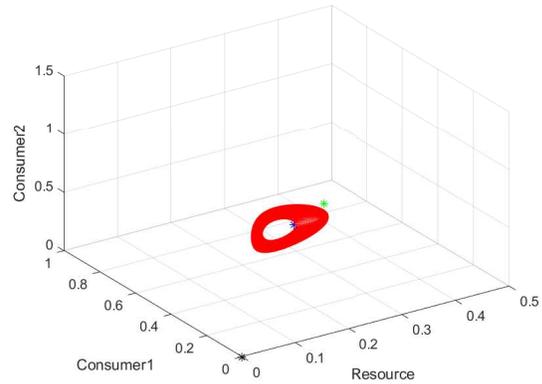}
		\caption{Dynamics in the region II in the form I}
	\end{subfigure}
	\caption{Three dimensional phase diagram of system \eqref{eq6} components for the type I of zero$-$Hopf normal form catalog, $a>0,b=+1$. The plots show the three-dimensional phase diagram of the system in regions I (a) and II (b) of the typical pattern. The dynamics of region III is symmetric to that of region I. Blue, green and black stars are initial values, positive, and boundary equilibrium, respectively.}
	\label{fig15}
\end{figure}

For this type of zero$-$Hopf bifurcation, Figure~\ref{fig14} (a) shows the typical dynamics pattern of type I. Figure~\ref{fig15} (a) shows the three dimensional phase diagram of the system considering competitive coefficients $\alpha =3.41$ and $\beta =0.13444002$, leading to normal form for $a= 0.378167$,  $b=+1$ and $\lambda_3=0$ corresponding to the region I in the typical dynamics pattern. The $\mathrm{FRR}$' threshold value was set as $0.16$ with a $+16$ improvement in attack rate after relaxation. In this case, consumer 2 is more effective in competition than consumer 1 with $\eta_{yz}=25.36$. Figure~\ref{fig15} (b) shows the related three dimensional phase diagram of the system for these competitive coefficients in region II in the typical dynamics pattern. Initial values for the simulation in the regions I and II are $(0.2,0.4,0.4)$ and $(0.1,0.01,1)$ respectively. Dynamics in regions I and  III are symmetric to each other. Symmetric parameters values are $a=b=1.2$, $c=d=0.5$, and $\mu = \nu = 0.03$ in this simulation. The tornado shape attractor is due to the commuted trajectories in region I with the typical dynamics pattern. Stable relaxation$-$oscillation is apparent in the time series of the resource and consumers in the region I of the form I catalog. Also, the dynamics in region II of form I, are stable relaxation-oscillation for all three variables.

\begin{figure}[h!]
	\centering
	\begin{subfigure}[b]{0.3\linewidth}
		\includegraphics[width=\linewidth]{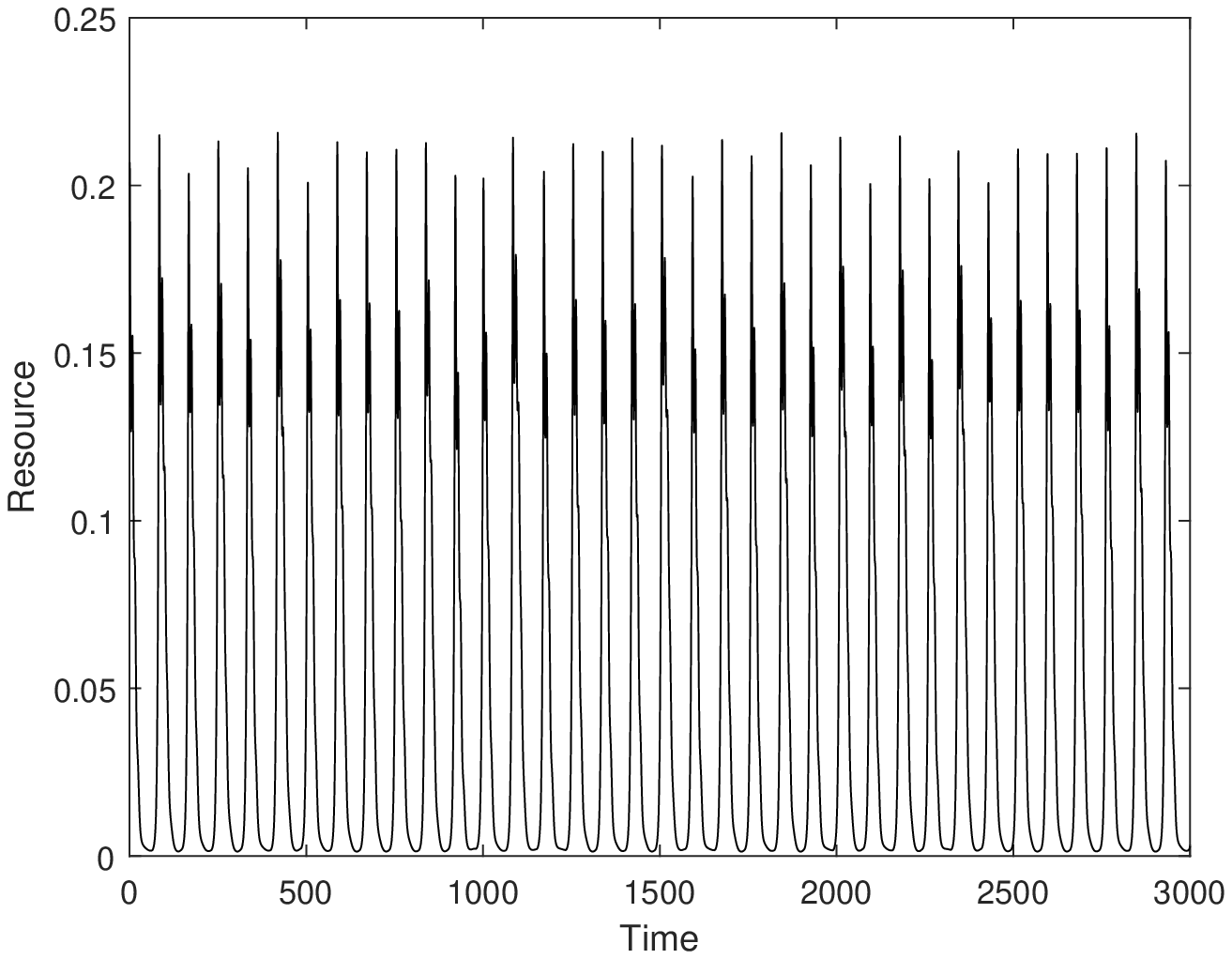}
		\caption{Resource success in region I}
	\end{subfigure}
	\begin{subfigure}[b]{0.3\linewidth}
		\includegraphics[width=\linewidth]{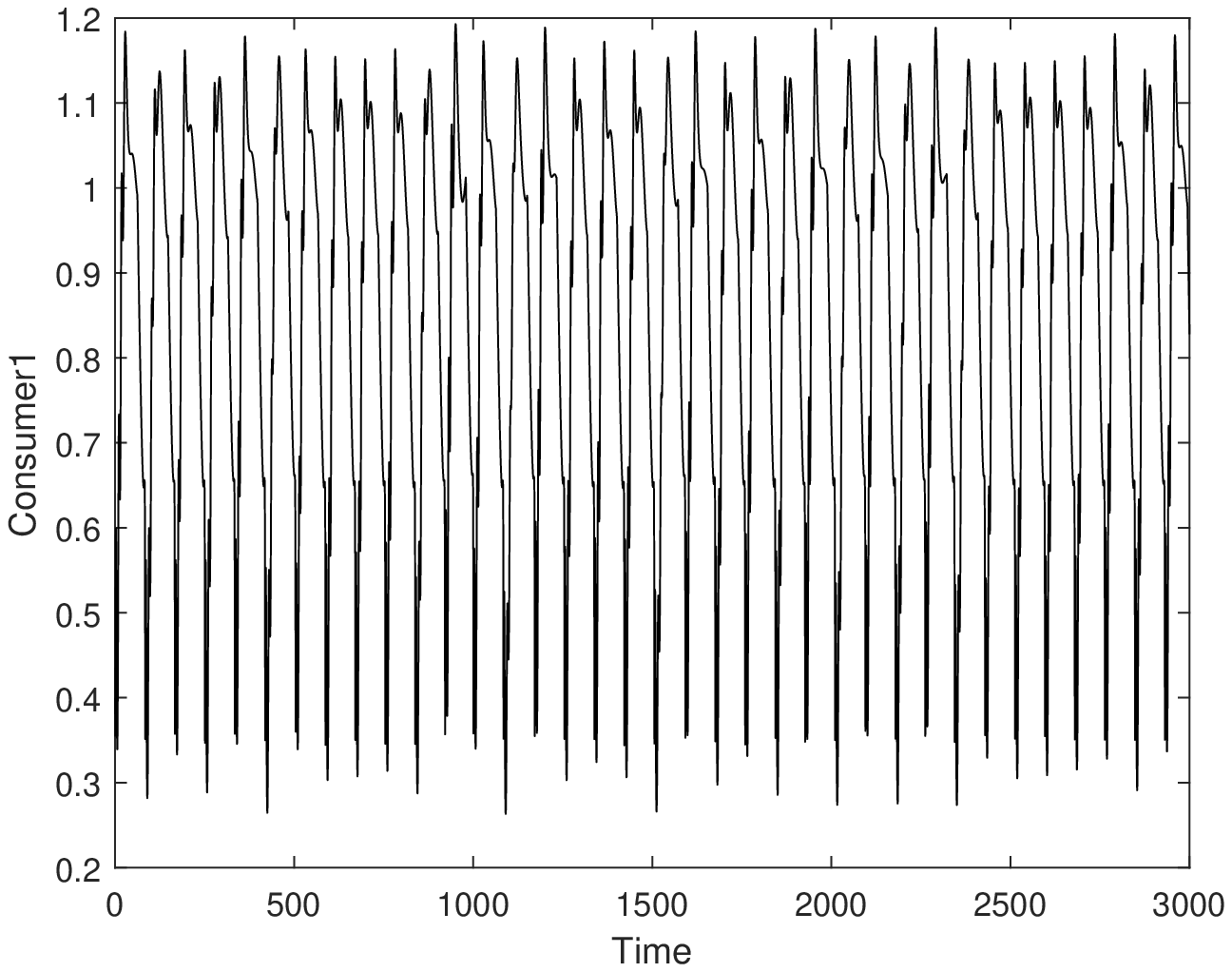}
		\caption{Consumer 1 success in region I}
	\end{subfigure}
	\begin{subfigure}[b]{0.3\linewidth}
		\includegraphics[width=\linewidth]{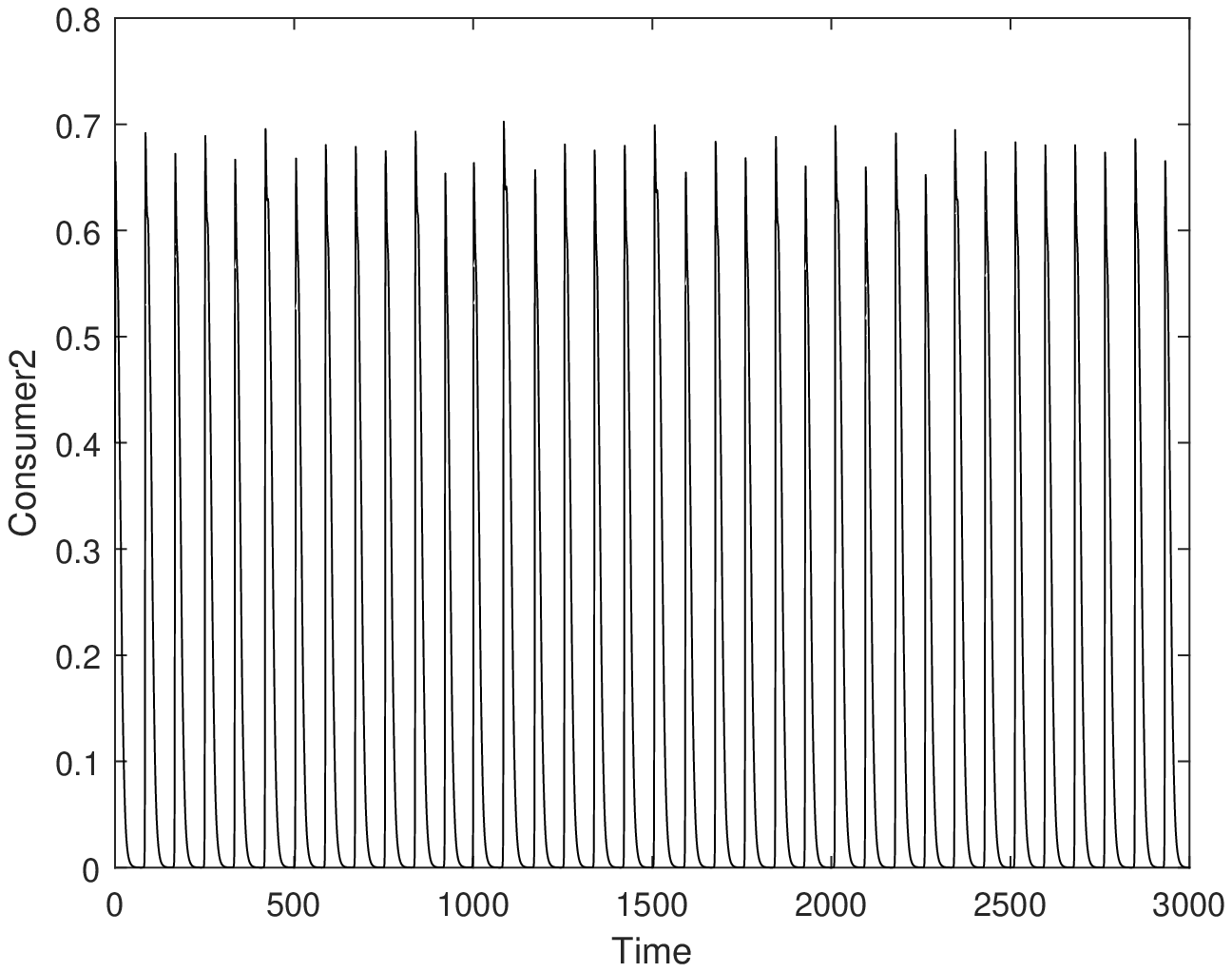}
		\caption{Consumer 2 success in region I}
	\end{subfigure}
	\begin{subfigure}[b]{0.3\linewidth}
		\includegraphics[width=\linewidth]{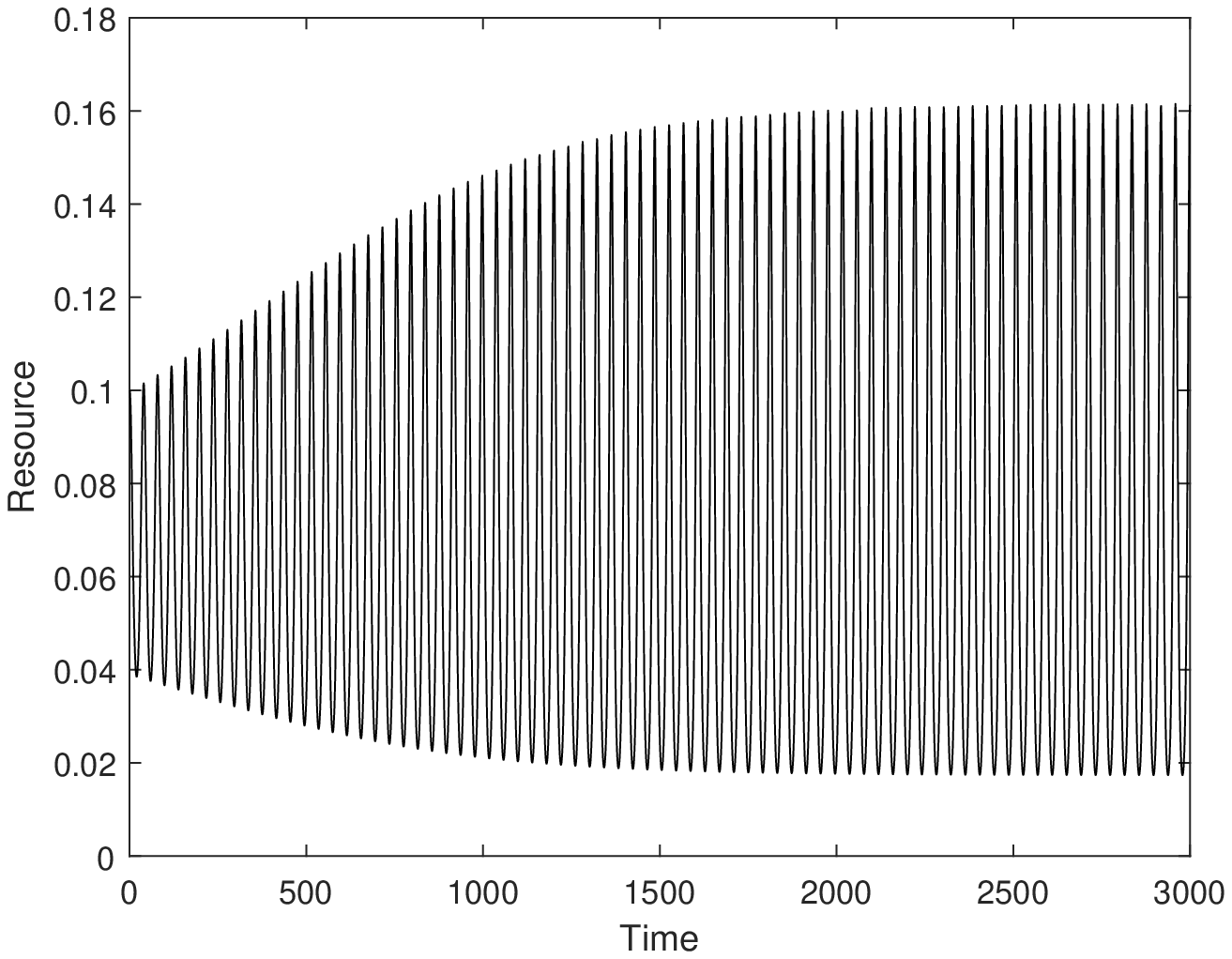}
		\caption{Resource success in region II}
	\end{subfigure}
	\begin{subfigure}[b]{0.3\linewidth}
		\includegraphics[width=\linewidth]{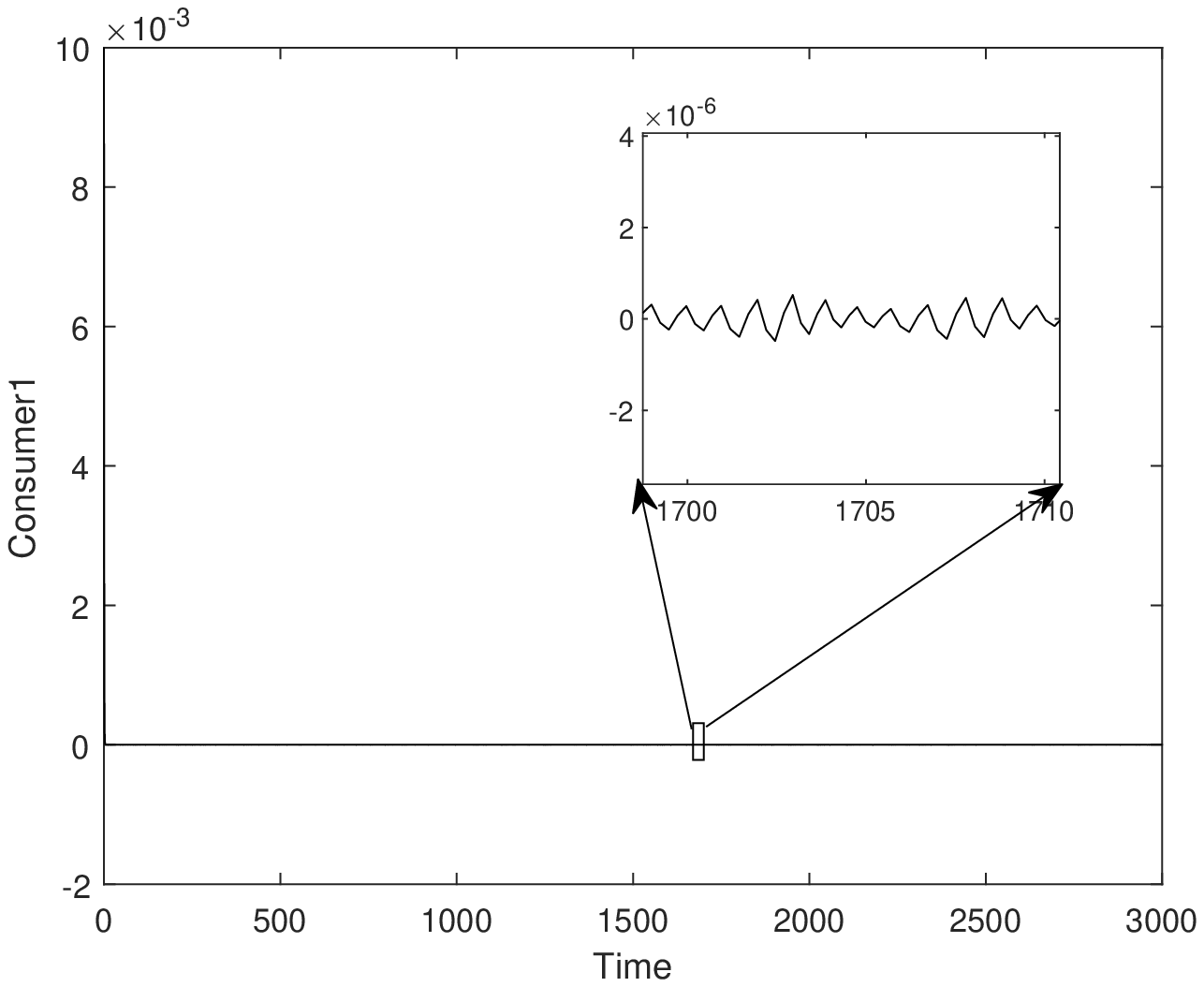}
		\caption{Consumer 1 success in region II}
	\end{subfigure}
	\begin{subfigure}[b]{0.3\linewidth}
		\includegraphics[width=\linewidth]{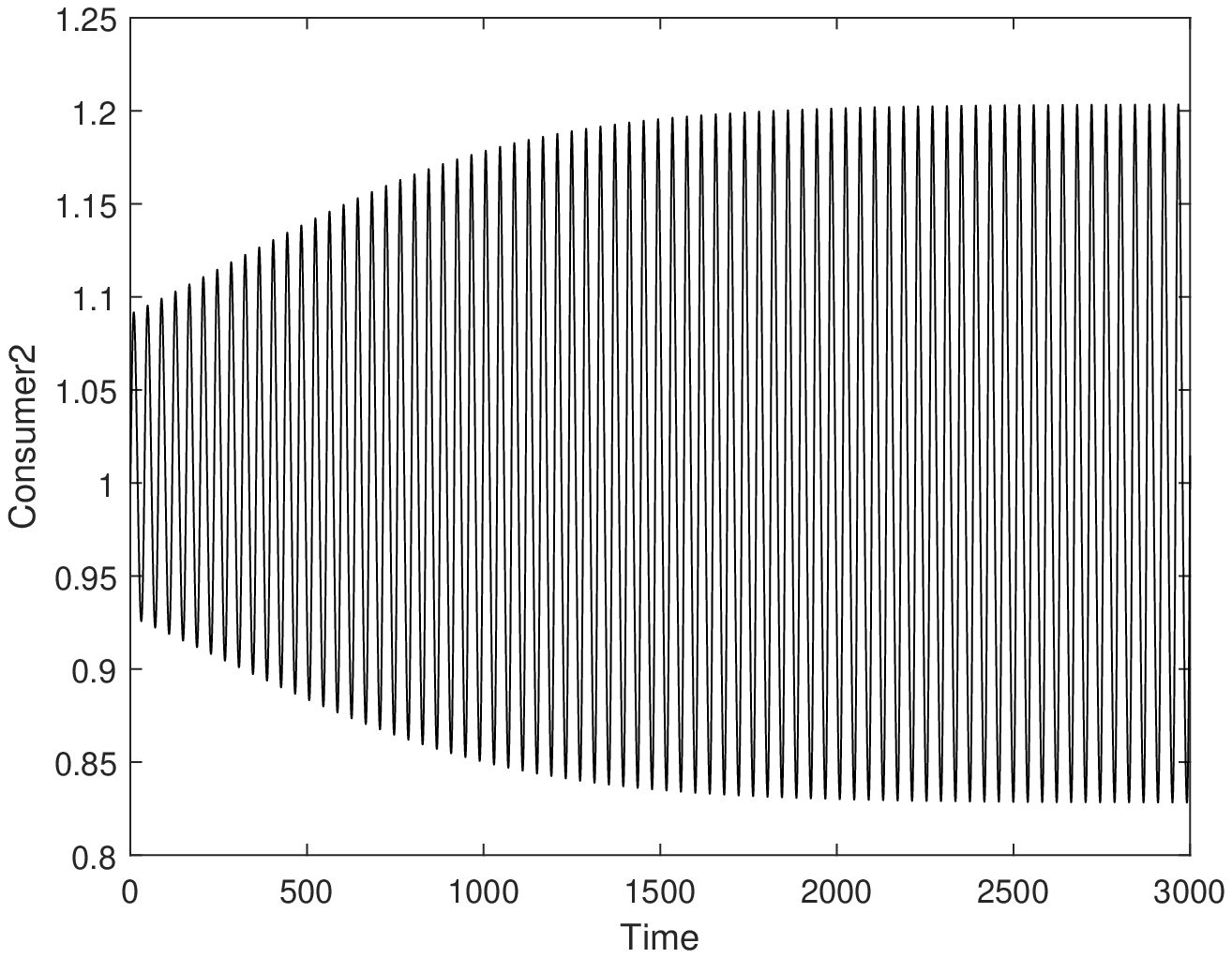}
		\caption{Consumer 2 success in region II}
	\end{subfigure}
	\caption{Time series of the system~\eqref{eq6} for the type I of zero$-$Hopf normal form catalog, $a>0,b=+1$. a,b, and c show time series in region$-$I and d, e, and f in region$-$II of the typical pattern. Temporal dynamics in region$-$III is symmetric to that of region$-$I.}
	\label{fig16}
\end{figure}

Time series of the system~\eqref{eq6}, undergoing type I zero$-$Hopf bifurcation shown in figure~\ref{fig16}. 

\subsubsection{Type IIb of zero$-$Hopf normal form catalog, $a<=-1,b=+1$}
\begin{figure}[h!]
	\centering
	\includegraphics[width=0.5\linewidth]{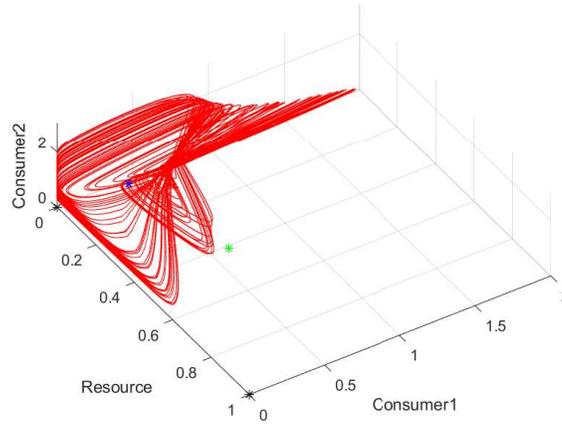}
	\caption{Three dimensional phase diagram of the system \eqref{eq6} for type IIb of zero$-$Hopf normal form catalog, $a<=-1,b=+1$. Dynamics in region$-$II and III is symmetric to those of region$-$I and II. Blue, green and black stars are initial values, positive, and boundary equilibrium, respectively.}
	\label{fig17}
\end{figure}

For this type of zero$-$Hopf bifurcation, figure~\ref{fig14} (b) shows the typical dynamics pattern of type IIb in which three regions are determined that show some dynamics according to respective initial values. Figure~\ref{fig17} shows three-dimensional phase diagram of the system~\eqref{eq6} in region$-$I and II of the typical pattern that visualizes the dynamics for competitive coefficients $\alpha =0.41$ and $\beta =0.19890678$, leading to normal form parameters $a= -9.63212$,  $b=+1$ with $\lambda_3=0$. The $\mathrm{FRR}$' threshold value is $0.058$ with a $+16$ improvement in attack rate after relaxation. In this case, consumer 2 is more effective in competition than consumer 1 with $\eta_{yz}=2.06$. Initial values for this simulation is $(0.06,0.4,0.4)$. Dynamics in region$-$II and III in the normal form IIb is symmetric to the other regions. Symmetric parameter values are $a=b=1.2$, $c=d=0.5$, and $\mu = \nu = 0.03$ in this simulation. The dynamics is similar to Gilpins' spiral chaos attractor which is refereed to Otto R\"ossler's continuous chaos in the Gilpins paper~\cite{gilpin79}. The least asymmetricity in competition is calculated in this bifurcation at which consumer 1 is $2.06$ times more effective than consumer 2.

\begin{figure}[h!]
	\centering
	\begin{subfigure}[b]{0.3\linewidth}
		\includegraphics[width=\linewidth]{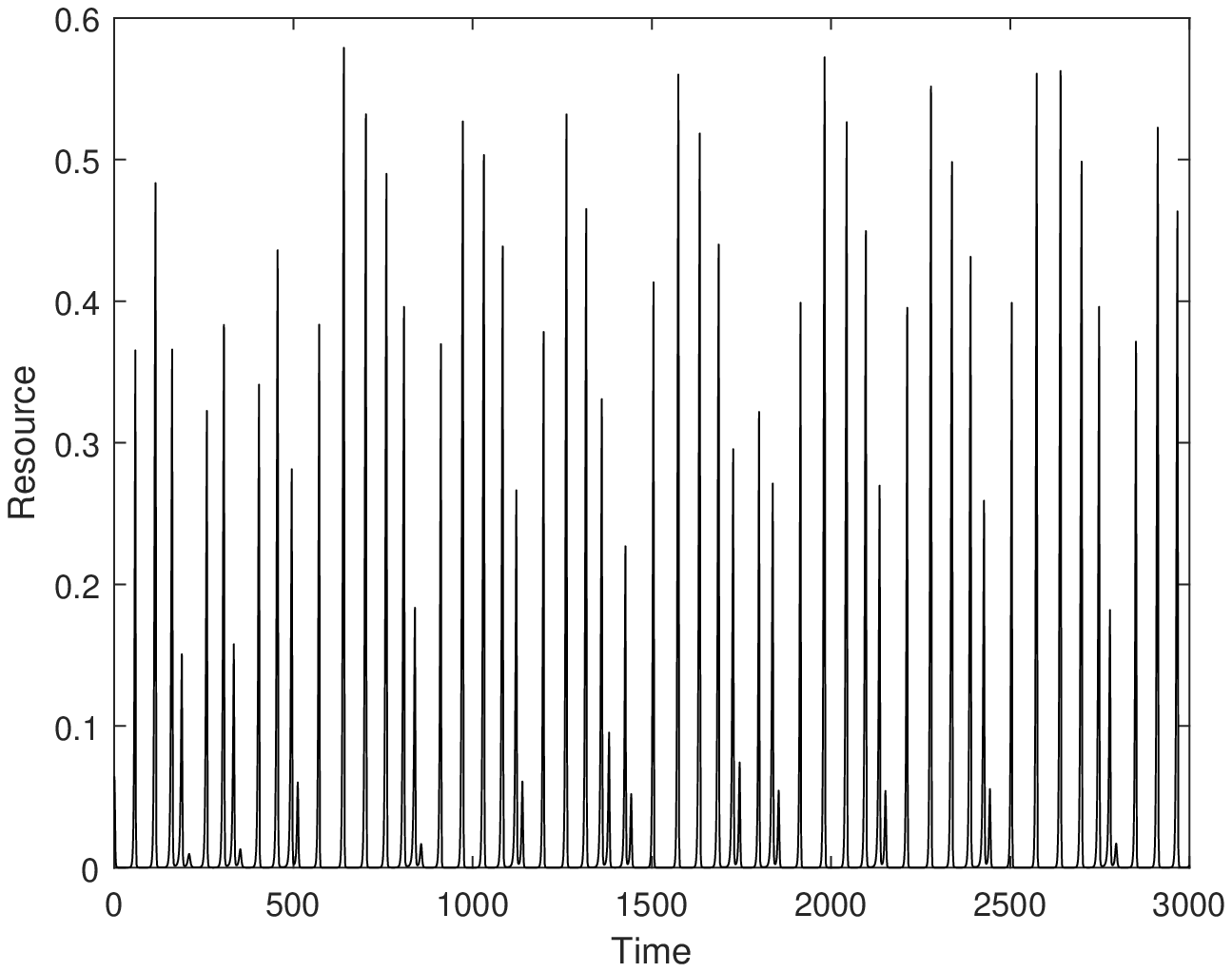}
		\caption{Resource success in region I+II}
	\end{subfigure}
	\begin{subfigure}[b]{0.3\linewidth}
		\includegraphics[width=\linewidth]{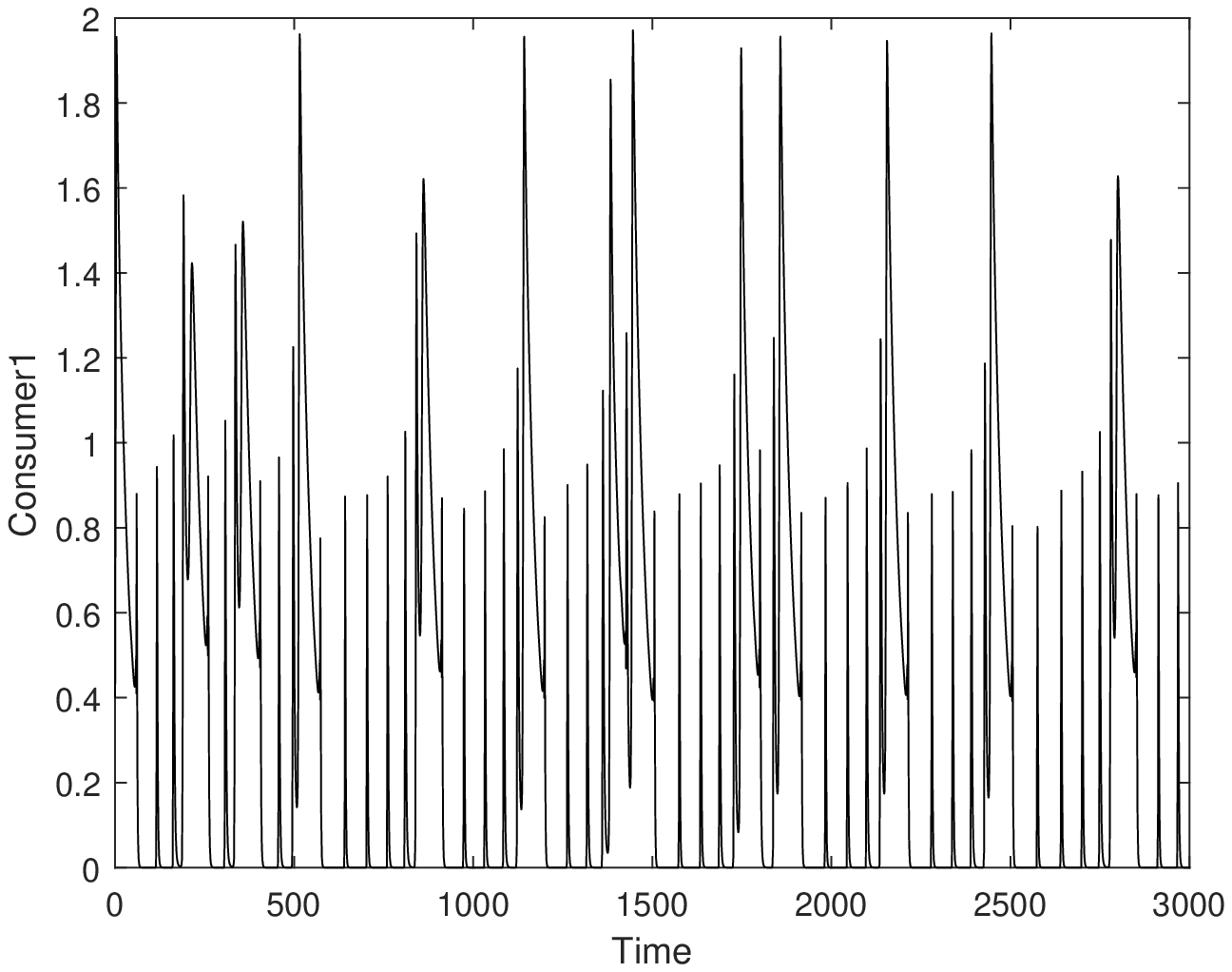}
		\caption{Consumer 1 success in region I+II}
	\end{subfigure}
	\begin{subfigure}[b]{0.3\linewidth}
		\includegraphics[width=\linewidth]{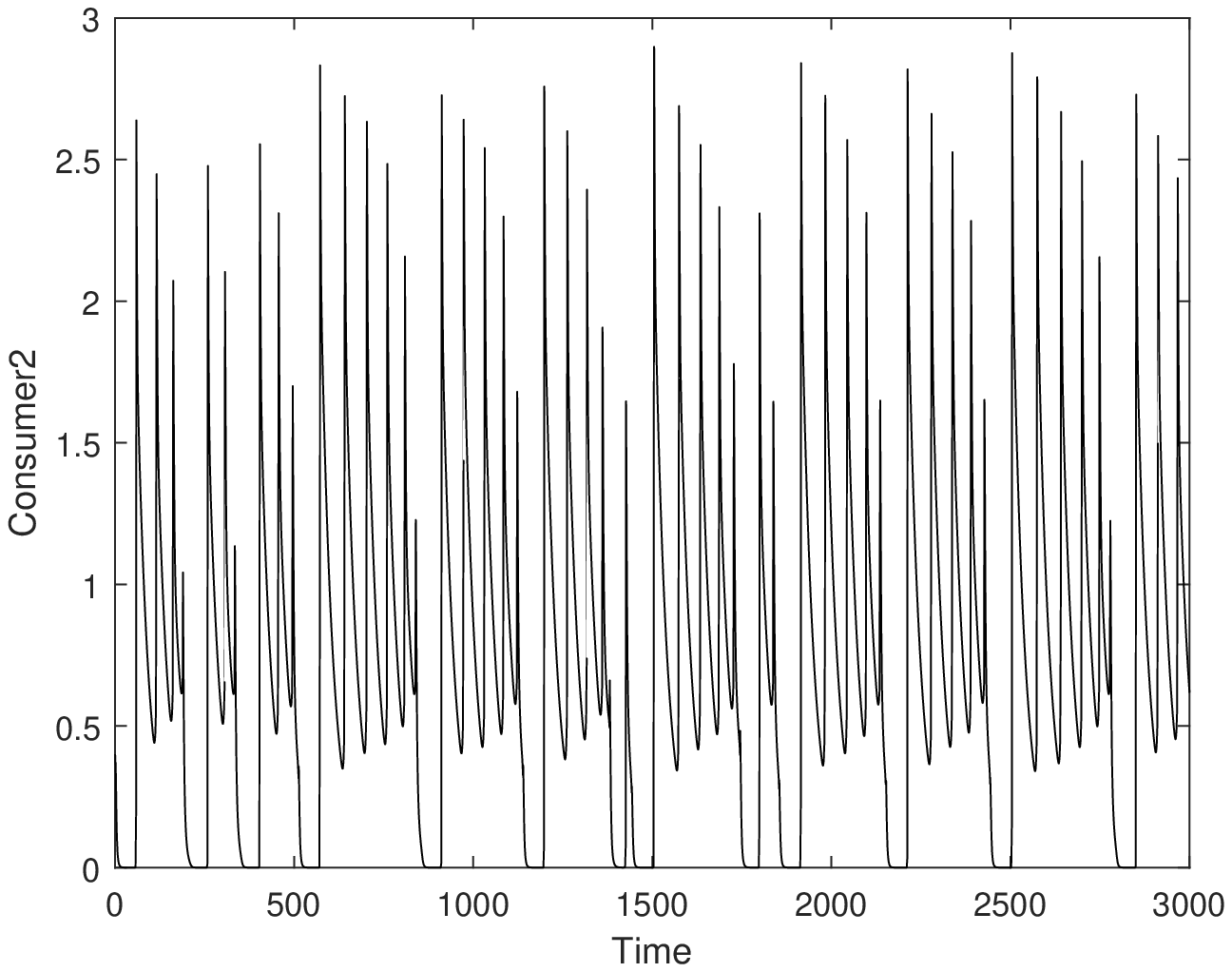}
		\caption{Consumer 2 success in region I+II}
	\end{subfigure}
	\caption{The time series of system~\eqref{eq6} of type IIb zero$-$Hopf bifurcation for $a<=-1 and b=+1$ in region$-$I and II shows the typical pattern of dynamics. Temporal dynamics in region$-$II and III is symmetric to that of region$-$I and II.}
	\label{fig18}
\end{figure}

Figure~\ref{fig18} illustrates the time series of the system~\eqref{eq6} undergoing type$-$IIb zero$-$Hopf bifurcation. 

\subsubsection{Type III of zero$-$Hopf normal form catalog, $a>0,b=-1$}
\begin{figure}[h!]
	\centering
	\includegraphics[width=0.5\linewidth]{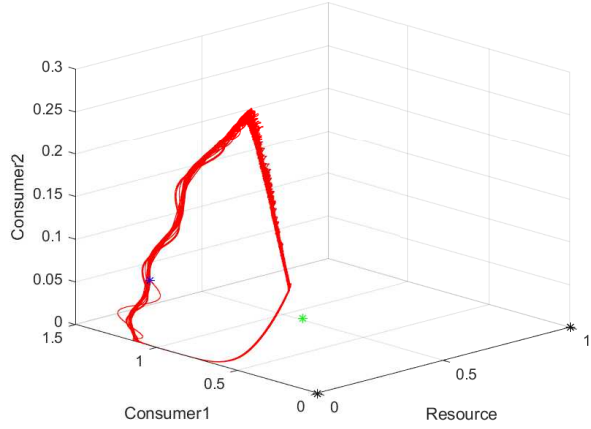}
	\caption{Three dimensional phase diagram of the system~\eqref{eq6} of type III zero$-$Hopf bifurcation with $a>0,b=-1$. The plot shows three dimensional phase diagram of the system in region$-$I and II. Blue, green and black stars are initial values, positive, and boundary equilibrium, respectively.}
	\label{fig19}
\end{figure}

For this type of zero$-$Hopf bifurcation, normal form stability, parametric values $\alpha=15.85492$ and $\beta=0.1454865$ are used to simulate the system~\eqref{eq6}. Figure~\ref{fig14} (c) shows the typical dynamics pattern of type III in which two regions are determined that show some dynamics according to respective initial values. Figure~\ref{fig19} shows three dimensional phase diagram of system~\eqref{eq6} visualizing the dynamics for the above competitive coefficients in the region$-$I and II of the typical pattern leading to normal form with $a= 0.063585$, $b=-1$ and the real root $\lambda_3=0$. Initial values for this simulation is $(0.04,1.1,0.07)$ and $\mathrm{FRR}$' threshold value, is $1$ with $+27$ improvement of attack rate after relaxation. In this case, consumer 2 is more effective in competition than consumer 1 with $\eta_{yz}=108.98$. In this simulation, symmetric parameter values are $a=b=1.2$, $c=d=0.5$, and $\mu = \nu = 0.03$. The jar shape attractor is depicted in three dimensional space. Stable relaxation oscillation is represented in the time series of resource and consumers 1 and 2 in region$-$I and II of the form III catalog \ref{fig14} (c). 

\begin{figure}[h!]
	\centering
	\begin{subfigure}[b]{0.3\linewidth}
		\includegraphics[width=\linewidth]{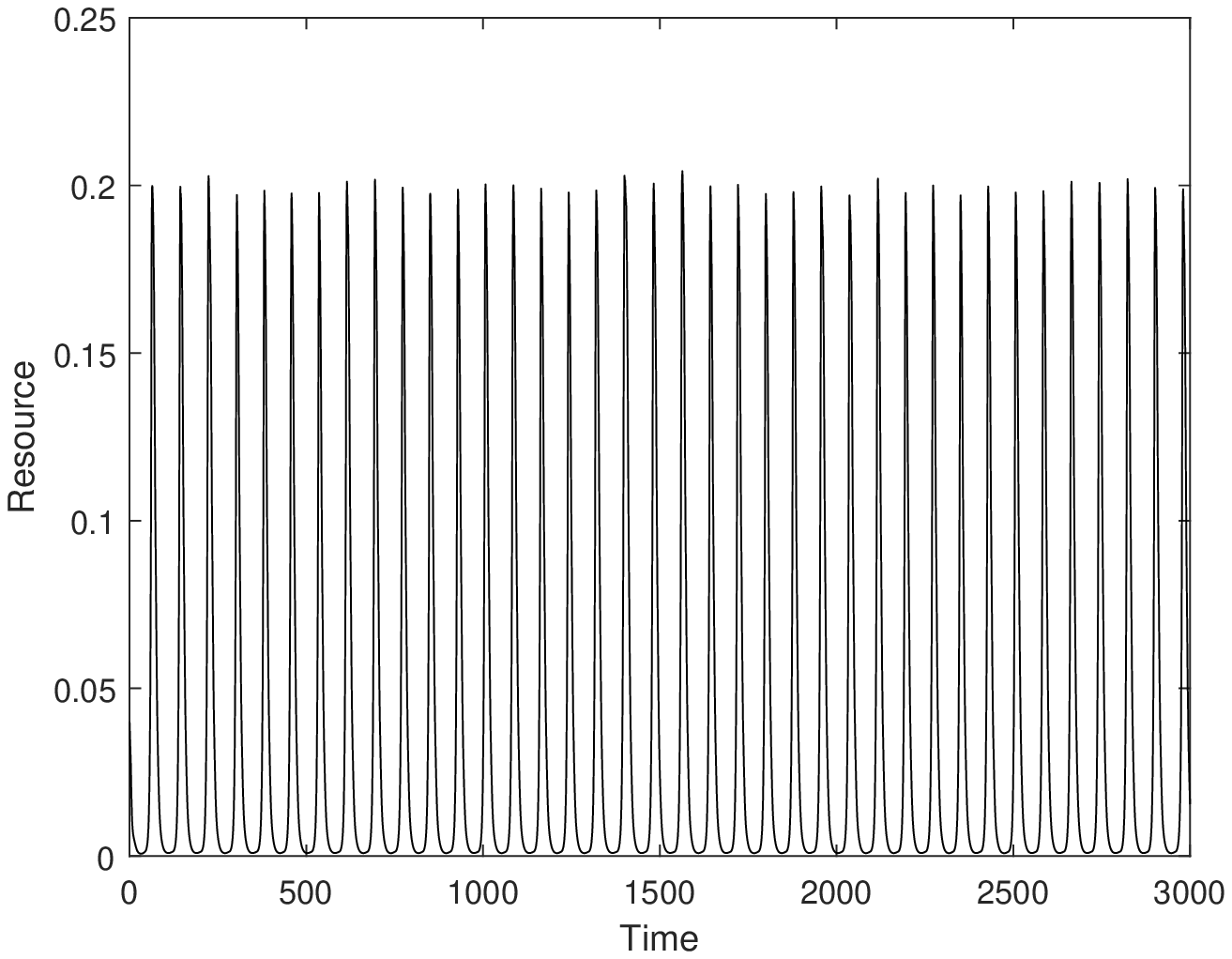}
		\caption{Resource success in region I+II}
	\end{subfigure}
	\begin{subfigure}[b]{0.3\linewidth}
		\includegraphics[width=\linewidth]{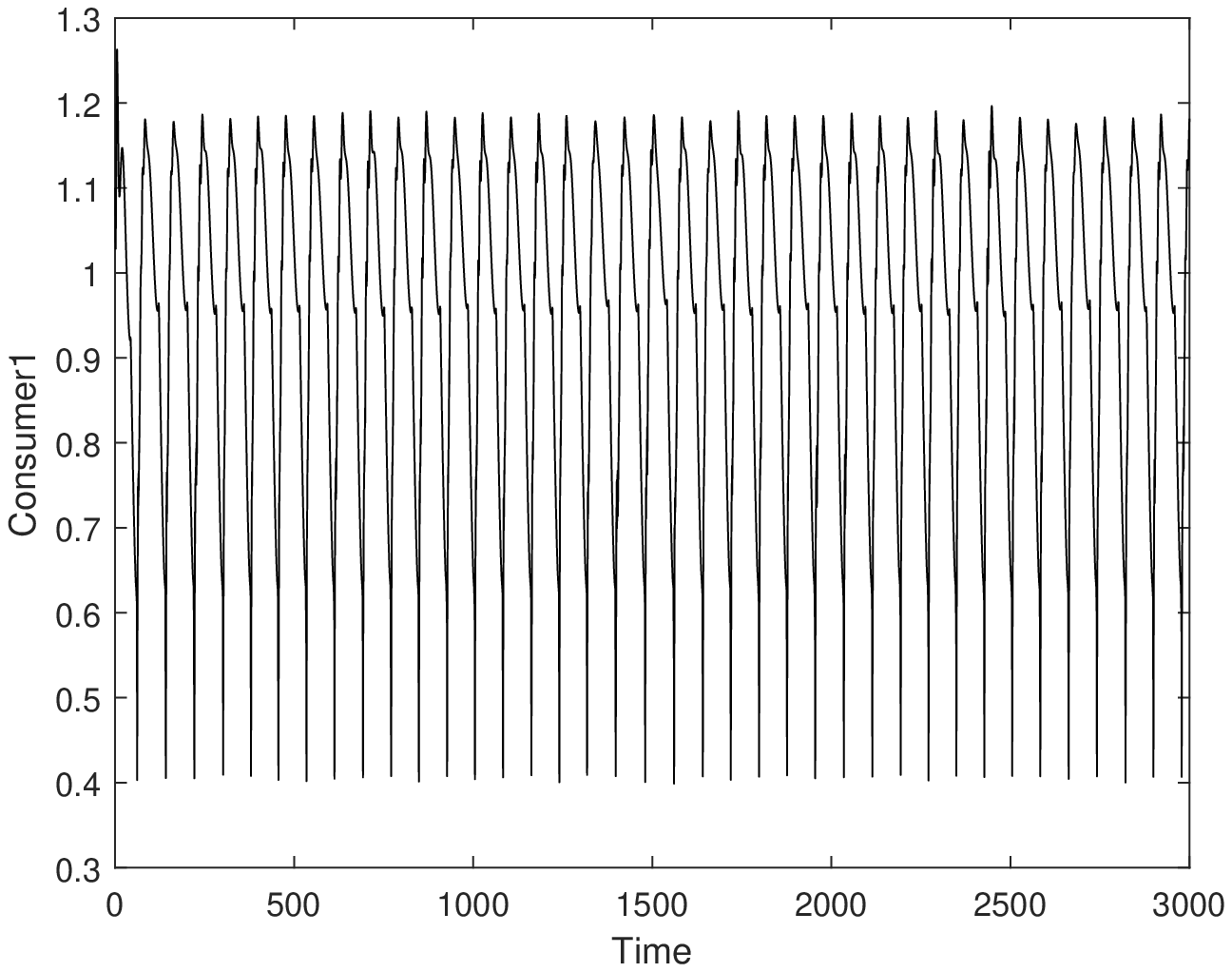}
		\caption{Consumer 1 success in region I+II}
	\end{subfigure}
	\begin{subfigure}[b]{0.3\linewidth}
		\includegraphics[width=\linewidth]{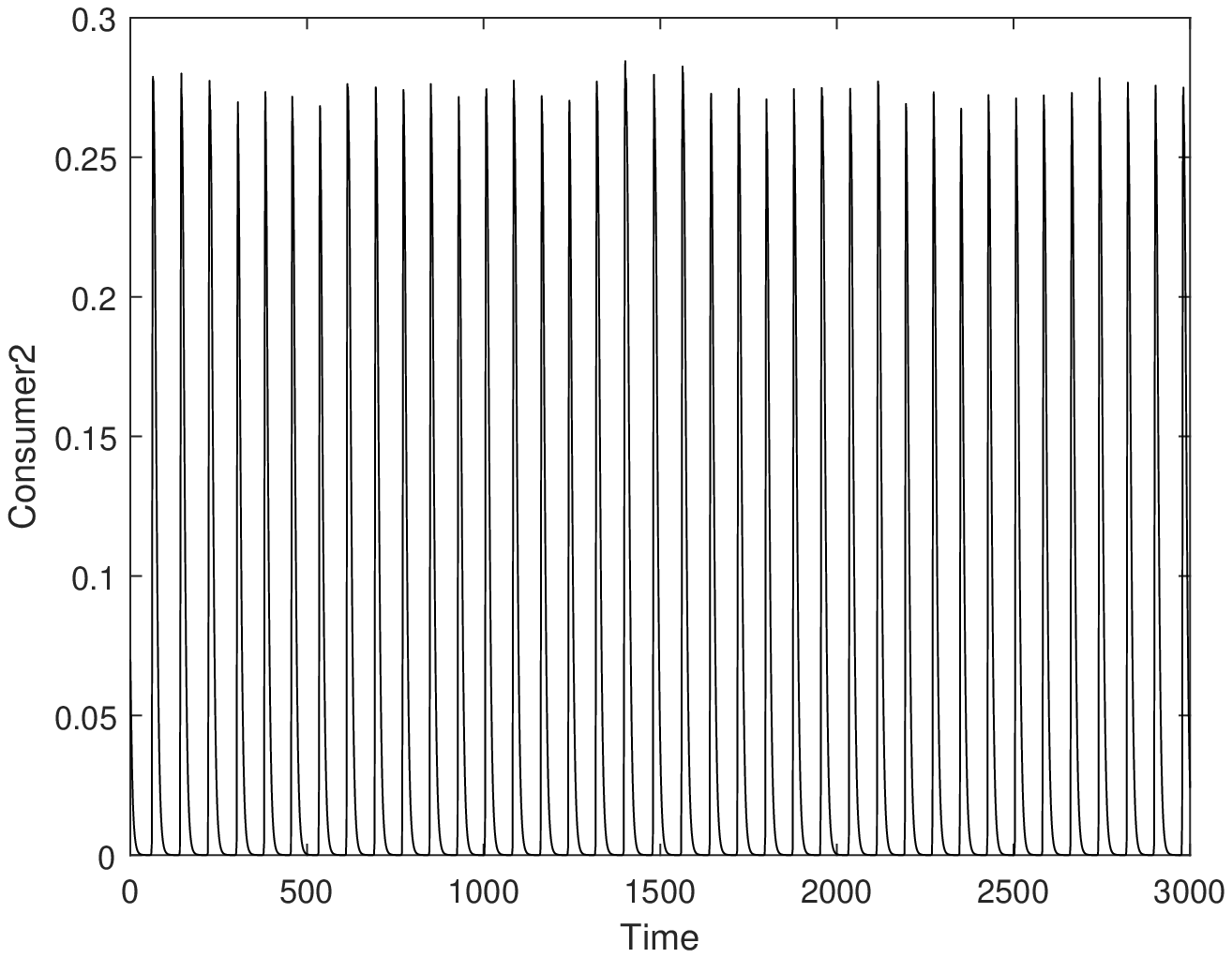}
		\caption{Consumer 2 success in region I+II}
	\end{subfigure}
	\caption{Time series of system~\eqref{eq6} for the type$-$III of zero$-$Hopf bifurcation with $a>0,b=-1$ in region$-$I and II with typical dynamics pattern.}
	\label{fig22}
\end{figure}

Time series of the system~\eqref{eq6} undergoing type$-$III of zero$-$Hopf bifurcation catalog are shown in figure~\ref{fig22}. 

\subsubsection{Type IVa of zero$-$Hopf normal form catalog, $-1<a<0,b=-1$}
\begin{figure}[h!]
	\centering
	\includegraphics[width=0.5\linewidth]{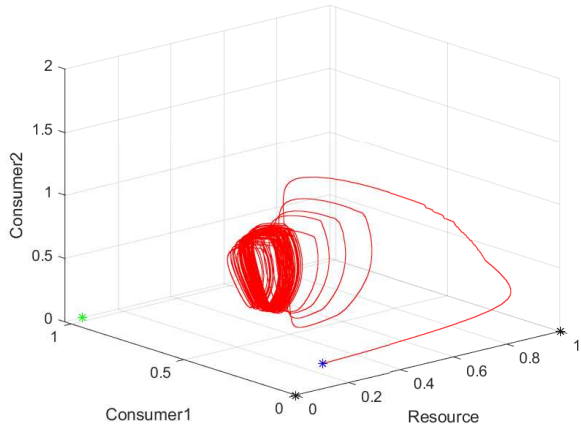}
	\caption{Three dimensional phase diagram of the system~\eqref{eq6} for type IVa of zero-Hopf bifurcation with  $-1<a<0, b=-1$. The plot shows three dimensional phase diagram of the system in region$-$I and II. Blue, green and black stars are initial values, positive, and boundary equilibrium, respectively.}
	\label{fig23}
\end{figure}

For this type of zero-Hopf bifurcation, parametric values $\alpha=0.000812293111191911$ and $\beta=0.00000000000099118319153423057911$ are used for simulation of the system~\eqref{eq6}. Figure~\ref{fig14} (d) shows the typical dynamics pattern of type IVa in which two regions are determined to show the dynamics separately in each one. Figure~\ref{fig23} (b) shows three-dimensional phase diagram of the system~\eqref{eq6} that visualizes the dynamics for the above competitive coefficient values in regions$-$I and II with $a= -0.271557$,  $b=-1$ and $\lambda_3=0$. Initial values for this simulation is $(0.1,0.0001,0.2)$ and $\mathrm{FRR}$' threshold value is $0.001$ with $+1$ improvement of attack rate after relaxation. In this case, consumer 2 is more effective in competition than consumer 1 with $\eta_{yz}=8.19 \times 10^{8}$. Symmetric parameters  values are $a=b=1.5$, $c=d=0.5$, and $\mu = \nu = 0.030008$. Rolling trajectories lead to an hourglass attractor. Stable relaxation oscillations, the time series of resource and consumers 1 and 2 in region$-$I and II of the form IVa in the catalog, are shown in Figure~\ref{fig23}. The most asymmetricity in this competition is that at which consumer 1 was more effective than consumer 2 with $\eta_{yz}=8.19 \times 10^{8}$. 

\begin{figure}[h!]
	\centering
	\begin{subfigure}[b]{0.3\linewidth}
		\includegraphics[width=\linewidth]{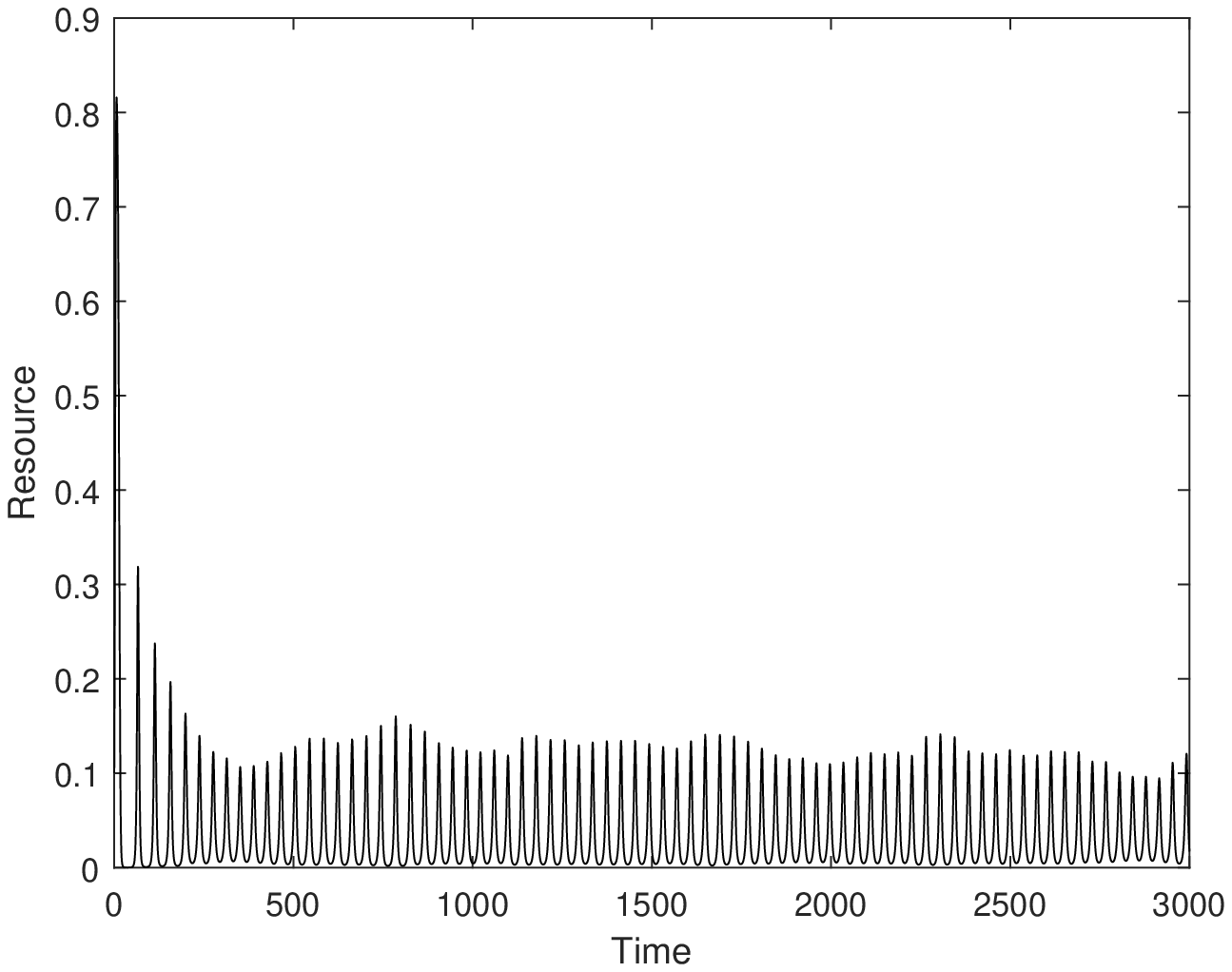}
		\caption{Resource success in region I+II}
	\end{subfigure}
	\begin{subfigure}[b]{0.3\linewidth}
		\includegraphics[width=\linewidth]{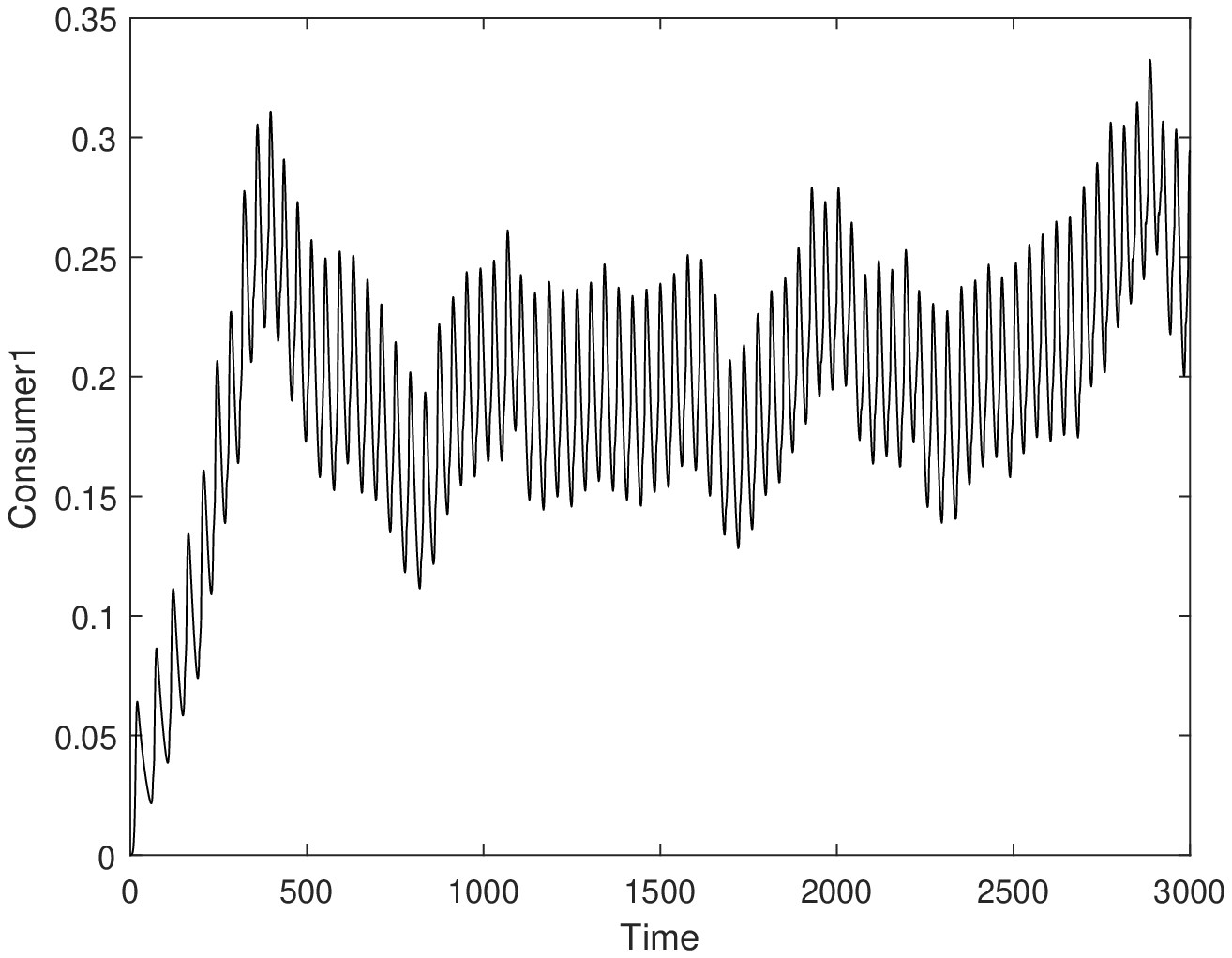}
		\caption{Consumer 1 success in region I+II}
	\end{subfigure}
	\begin{subfigure}[b]{0.3\linewidth}
		\includegraphics[width=\linewidth]{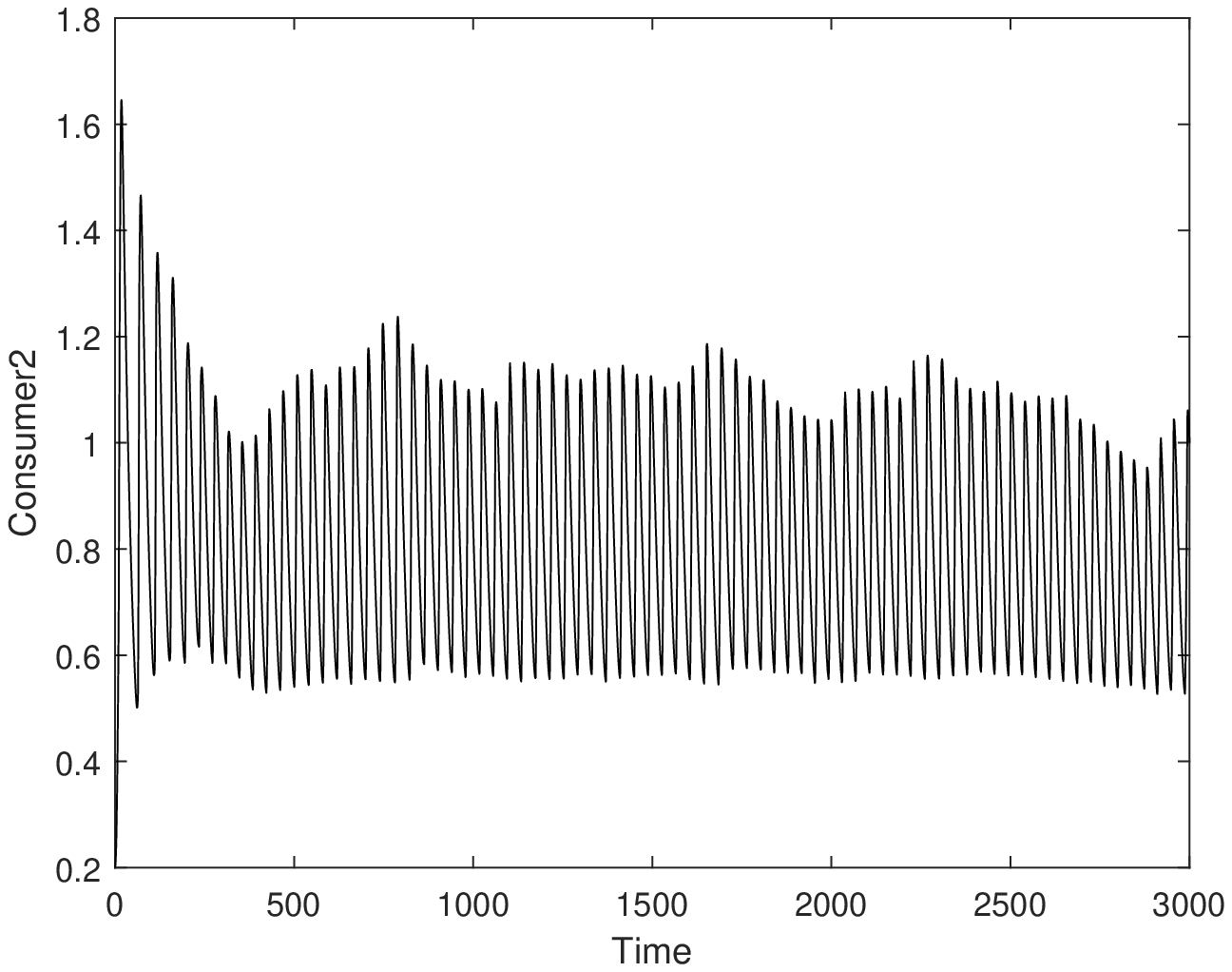}
		\caption{Consumer 2 success in region I+II}
	\end{subfigure}
	\caption{Time series of system~\eqref{eq6} success for the type IVa zero$-$Hopf bifurcation with $-1<a<0, b=-1$ in region$-$I and II.}
	\label{fig24}
\end{figure}

Time series of the system~\eqref{eq6}, undergoing type IVa zero$-$Hopf bifurcation, are shown in figure~\ref{fig24}. 

\subsubsection{Type IVb of zero$-$Hopf normal form catalog, $a<-1,b=-1$}
\begin{figure}[h!]
	\centering
	\begin{subfigure}[b]{0.45\linewidth}
		\includegraphics[width=\linewidth]{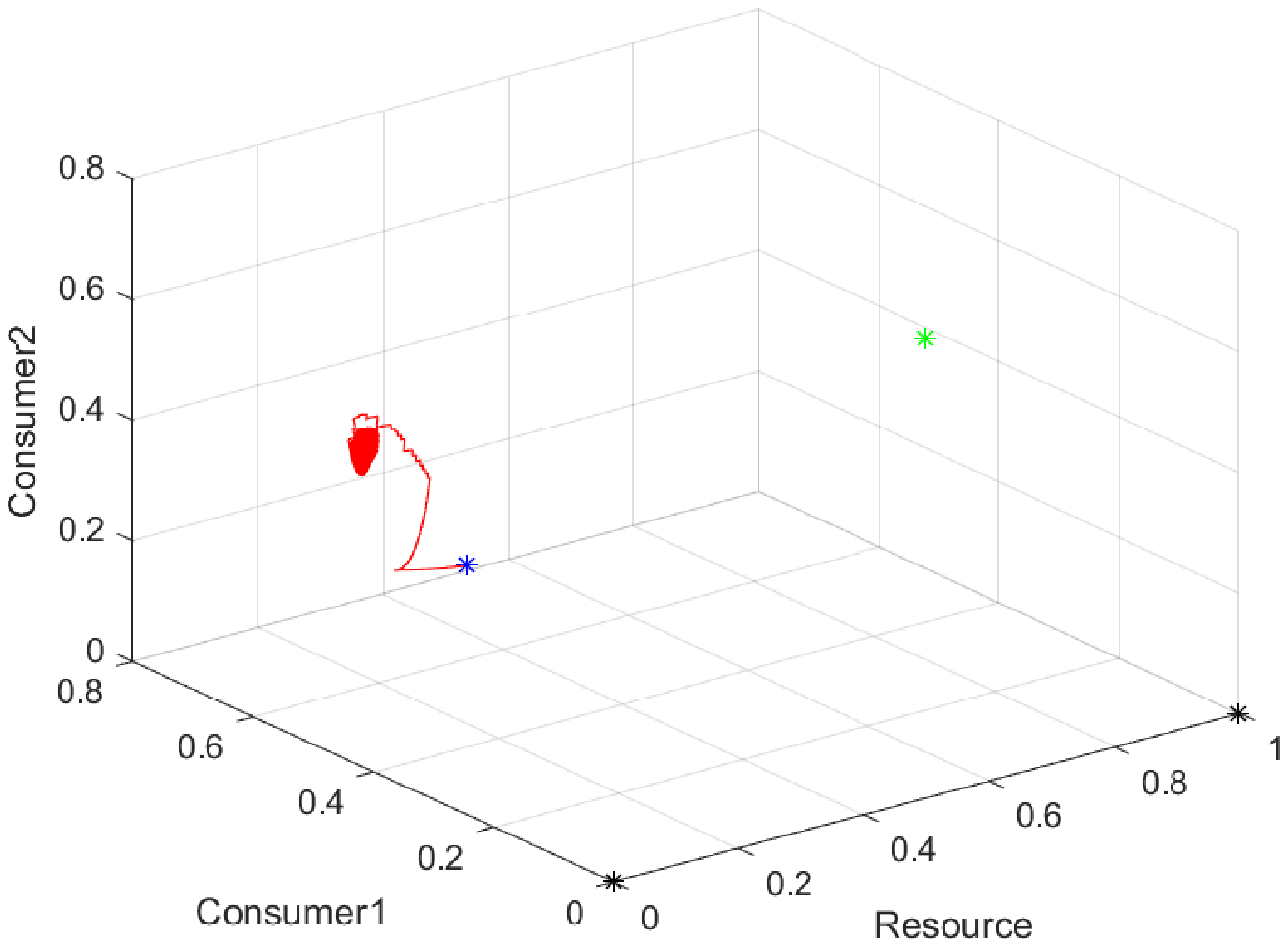}
		\caption{Dynamics in the region I in the form IVb}
	\end{subfigure}
	\begin{subfigure}[b]{0.45\linewidth}
		\includegraphics[width=\linewidth]{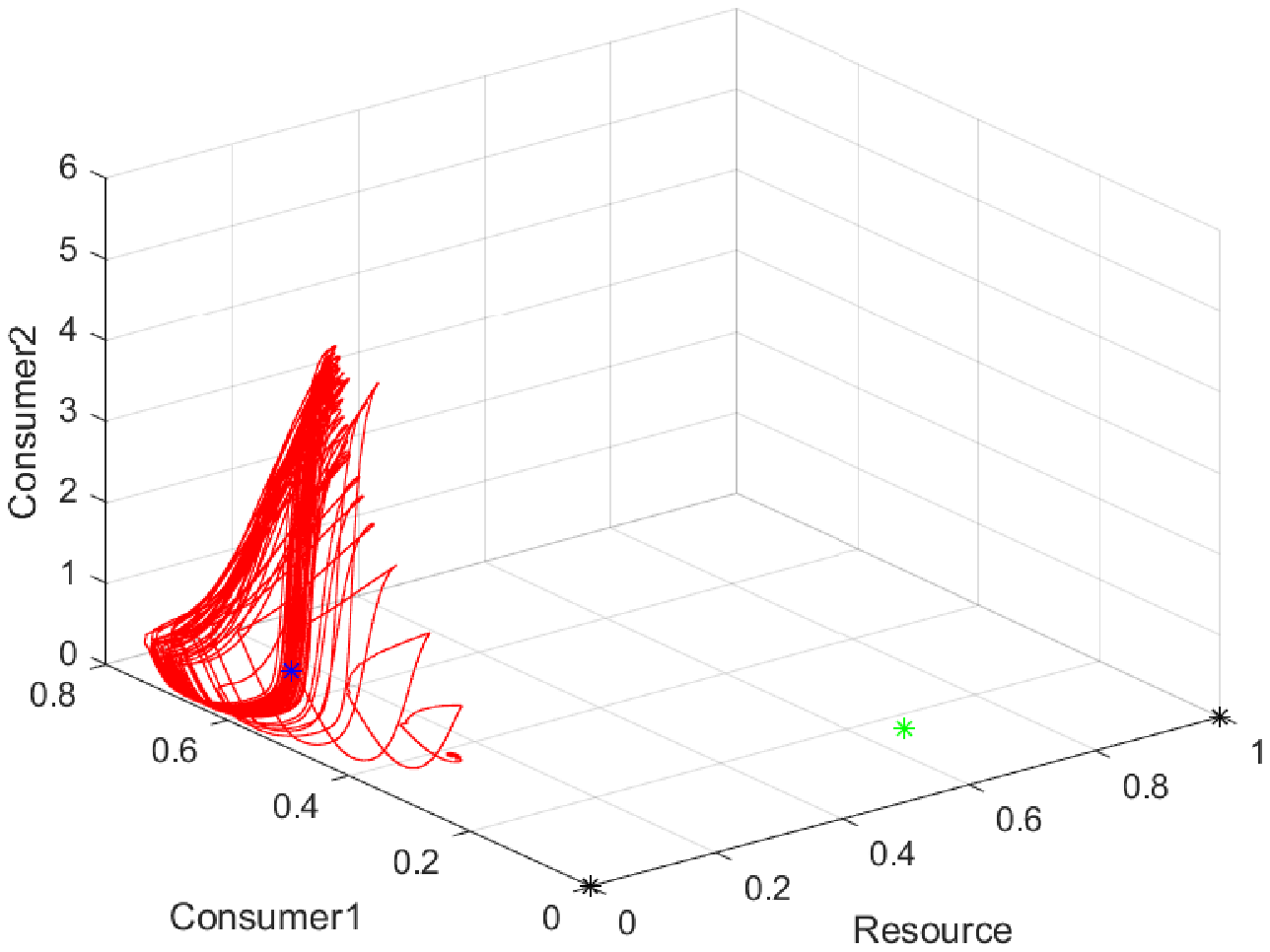}
		\caption{Dynamics in the region II in the form IVb}
	\end{subfigure}
	\caption{Three dimensional phase diagram of the system~\eqref{eq6} for type IVb of zero$-$Hopf bifurcation with $a<-1,b=-1$. The plot shows the three dimensional phase diagram of the system in region$-$I (a) and region$-$II (b), respectively. Dynamics in the region$-$III is symmetric to the I one. Blue, green and black stars are initial values, positive, and boundary equilibrium, respectively. }
	\label{fig21}
\end{figure}

For this type of zero$-$Hopf normal form bifurcation, we used parametric values $\alpha=0.1783636$ and $\beta=4.170286$ to simulate the system \eqref{eq6}. Figure \ref{fig14} (e) shows the typical dynamics pattern of type IVb in which three regions are determined that show some dynamics according to respective initial values. Figure \ref{fig21} (a) shows three dimensional phase diagram of system \eqref{eq6} components that visualizes the dynamics for the above competitive coefficients values in region I of the typical pattern leading to normal form $a= -4.36788$, $b=-1$ and the real root $\lambda_3=0$. The $\mathrm{FRR}$ threshold value is $0.32$ with $+16$ improvements of attack rate after relaxation. In this case, consumer 1 is more effective in competition than consumer 2 with $\eta_{zy}=23.38$. Figure \ref{fig21} (b) shows the three dimensional phase diagram of the system considering the above competitive coefficients in the region$-$II in the typical dynamics pattern. Initial values for the simulation in the regions I and II are $(0.15,0.4,0.32)$ and $(0.15,0.4,0.3)$, respectively. Dynamics in the region$-$III is symmetric to the region$-$I. The threshold value of $\mathrm{FRR}$ is $0.313$ with $+66$ improvement of attack rate after relaxation for the dynamics in the region$-$I. Parameters symmetric values are $a=b=1.2$, $c=d=0.5$, and $\mu = \nu = 0.03$ in this simulation. Dynamics in the region$-$I is a conical attractor constructed by rolling trajectories. A limit cycle forms in the region$-$II. Temporal dynamics in region$-$I is stable and do not shape a regular relaxation oscillation. But, in region$-$II, the dynamics form regular relaxation$-$oscillations.

\begin{figure}[h!]
	\centering
	\begin{subfigure}[b]{0.3\linewidth}
		\includegraphics[width=\linewidth]{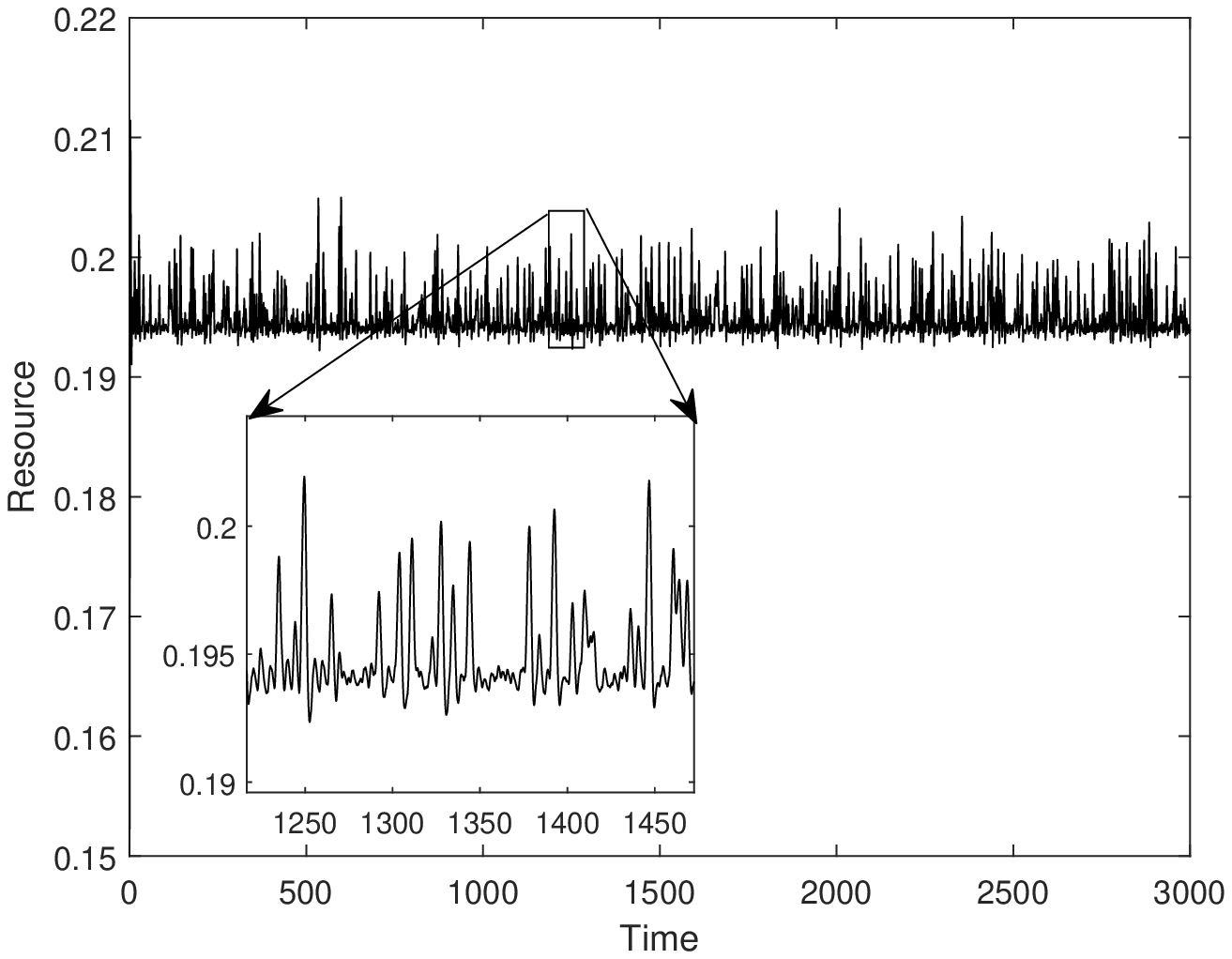}
		\caption{Resource success in region I}
	\end{subfigure}
	\begin{subfigure}[b]{0.3\linewidth}
		\includegraphics[width=\linewidth]{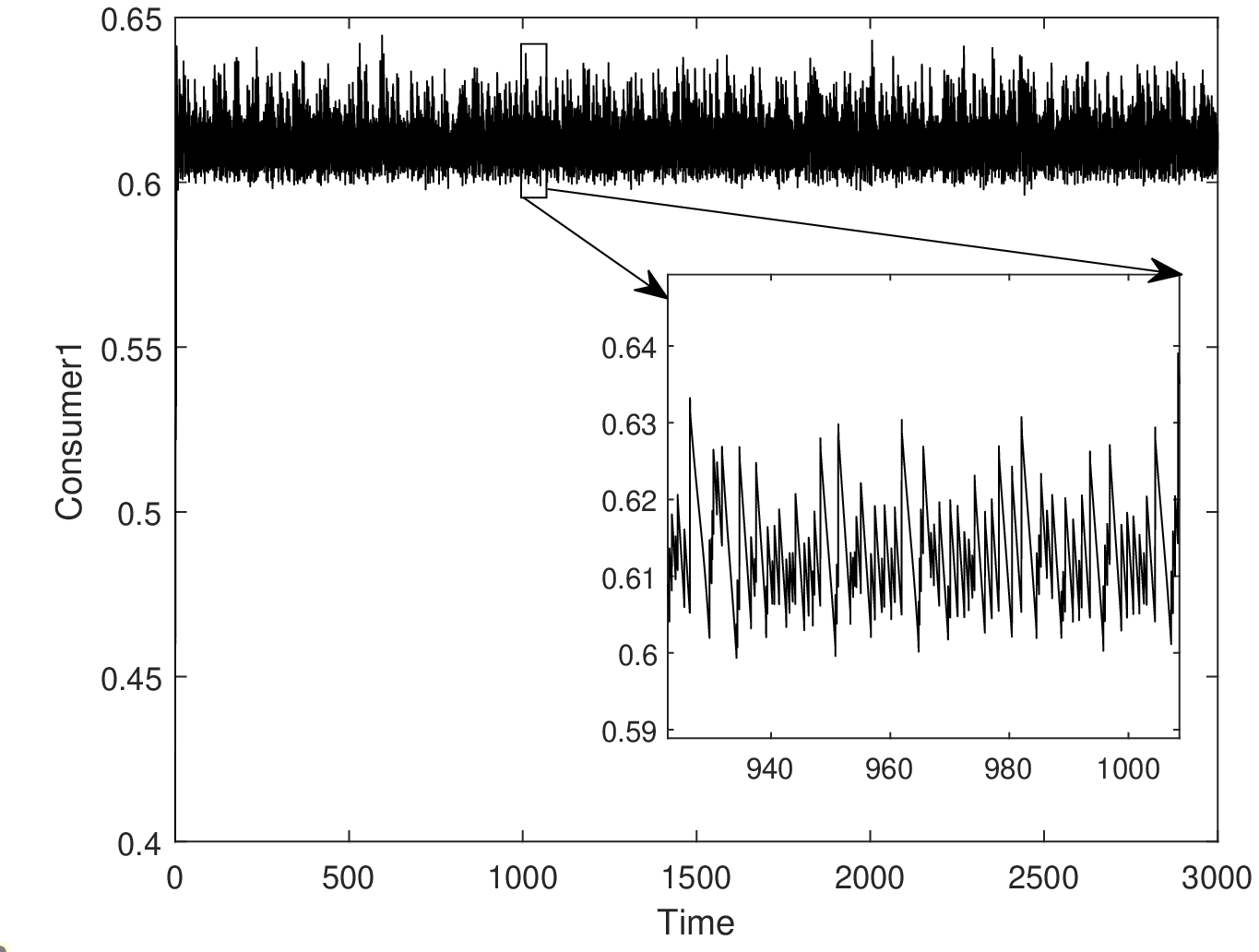}
		\caption{Consumer 1 success in region I}
	\end{subfigure}
	\begin{subfigure}[b]{0.3\linewidth}
		\includegraphics[width=\linewidth]{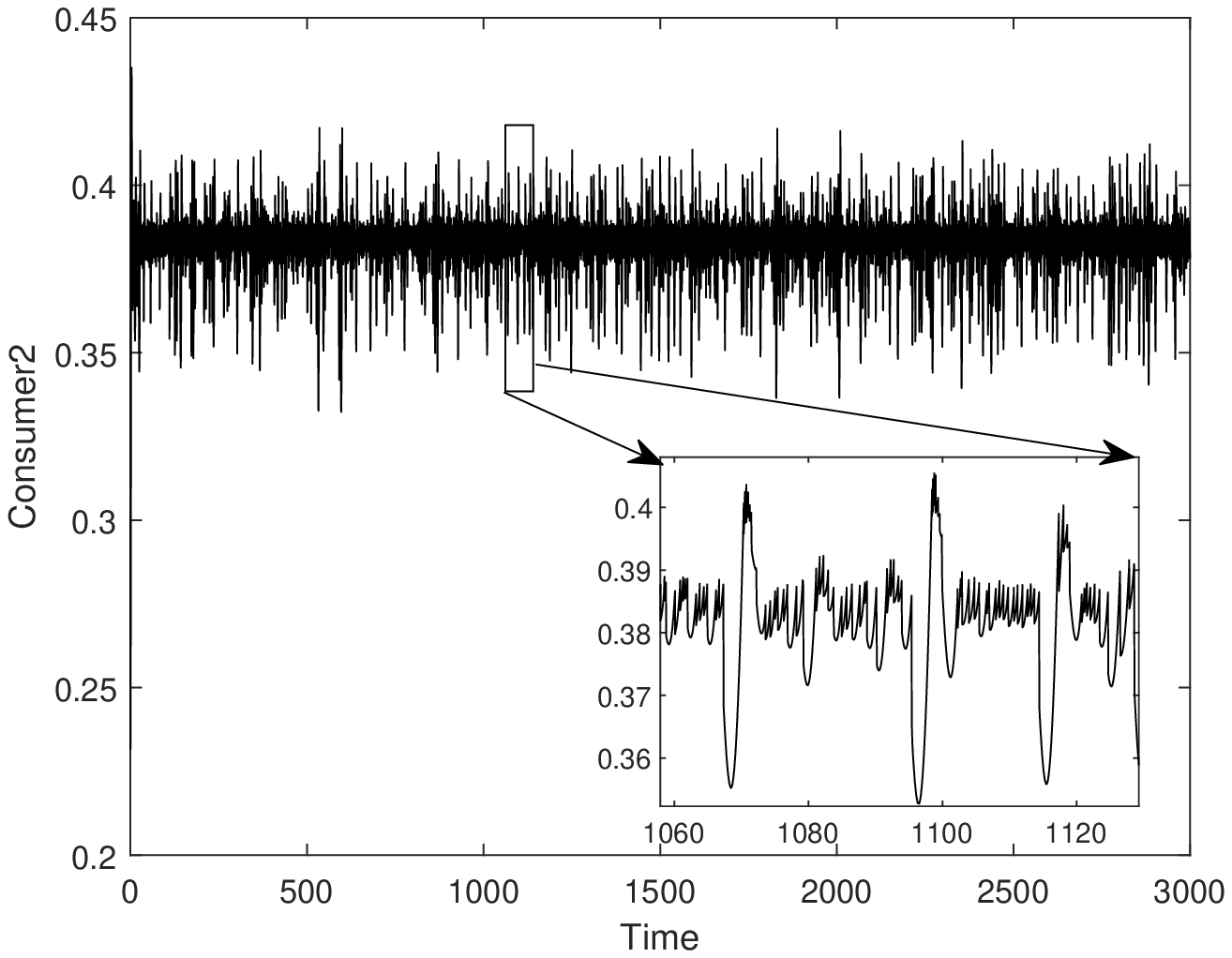}
		\caption{Consumer 2 success in region I}
	\end{subfigure}
	\begin{subfigure}[b]{0.3\linewidth}
		\includegraphics[width=\linewidth]{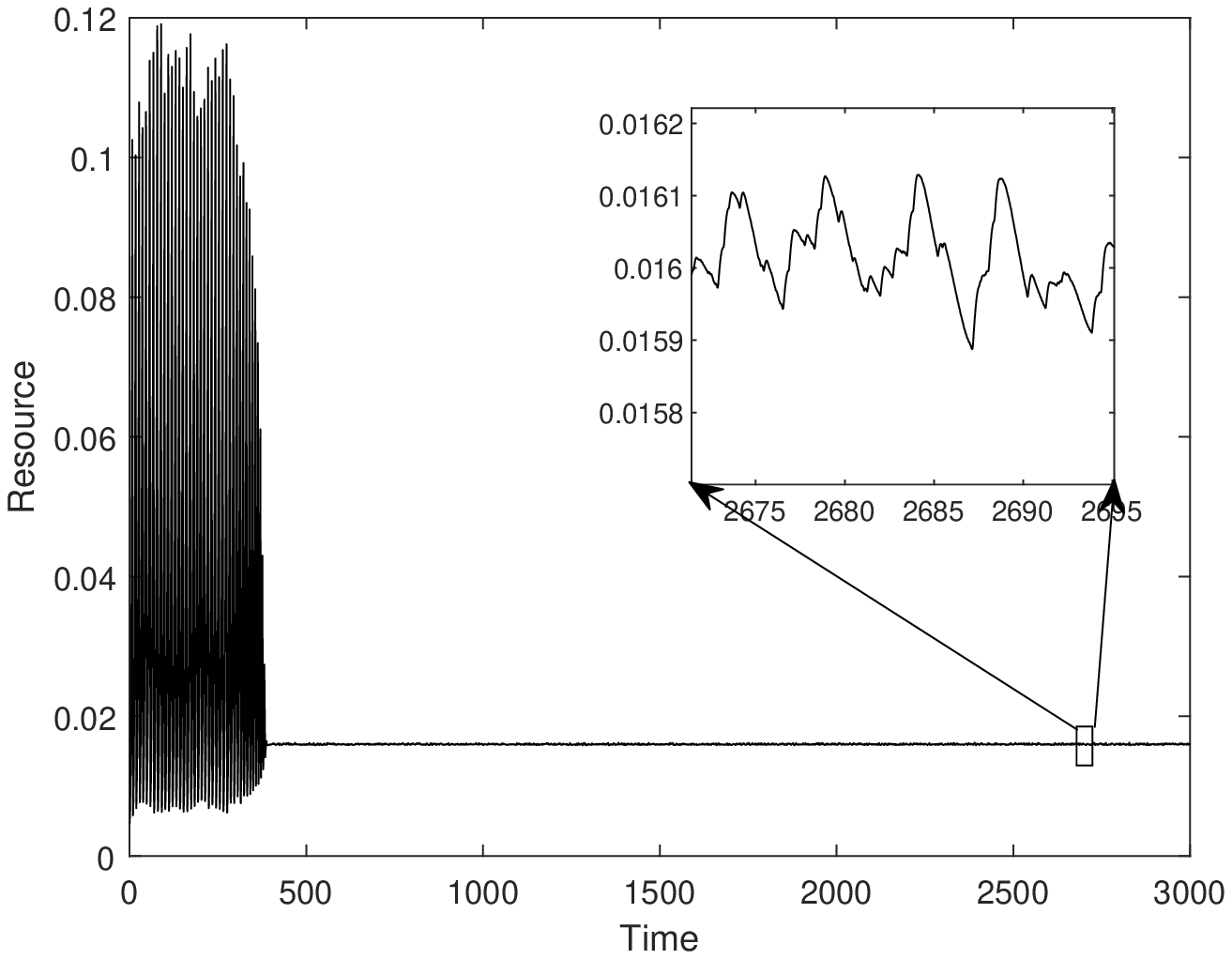}
		\caption{Resource success in region II}
	\end{subfigure}
	\begin{subfigure}[b]{0.3\linewidth}
		\includegraphics[width=\linewidth]{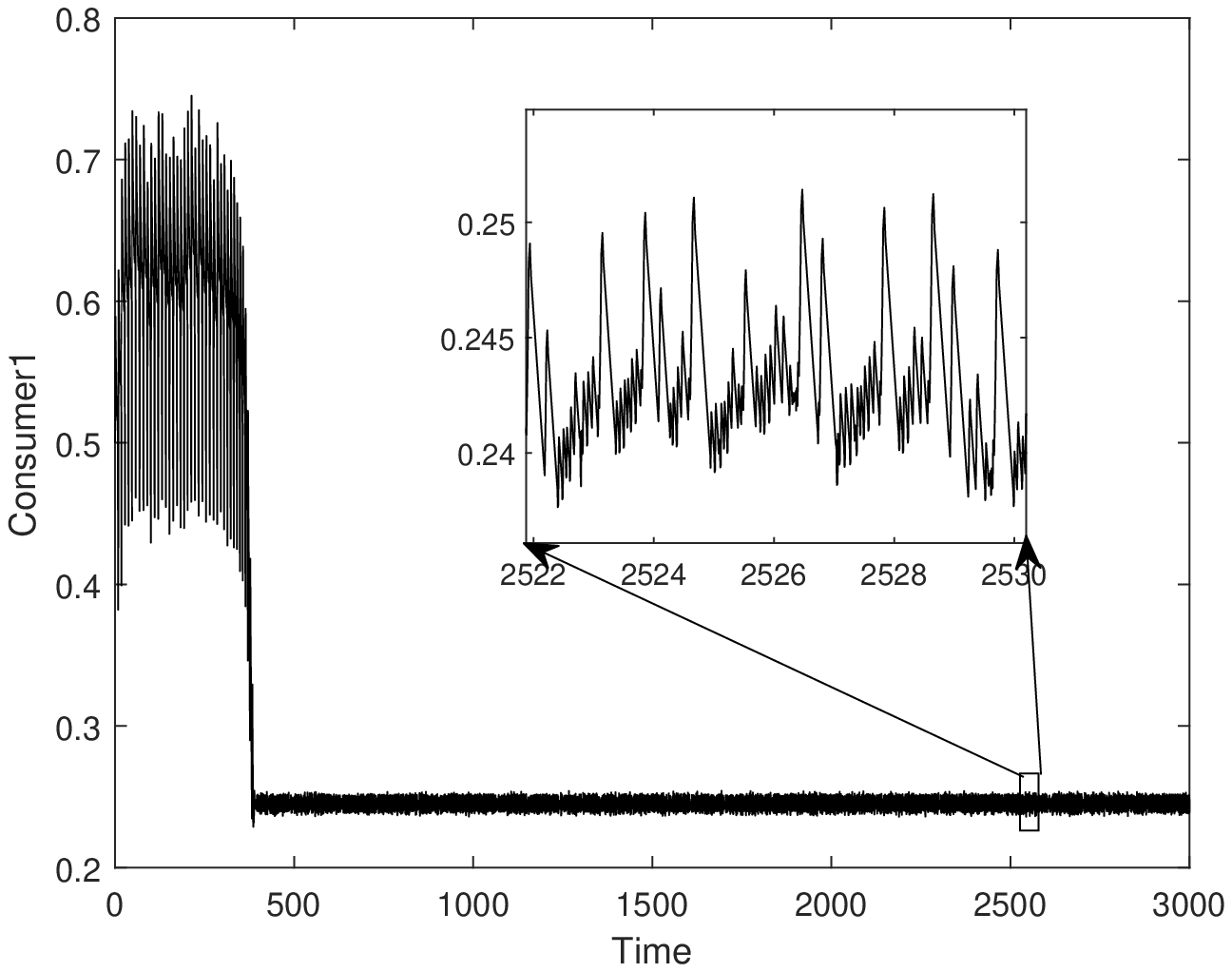}
		\caption{Consumer 1 success in region II}
	\end{subfigure}
	\begin{subfigure}[b]{0.3\linewidth}
		\includegraphics[width=\linewidth]{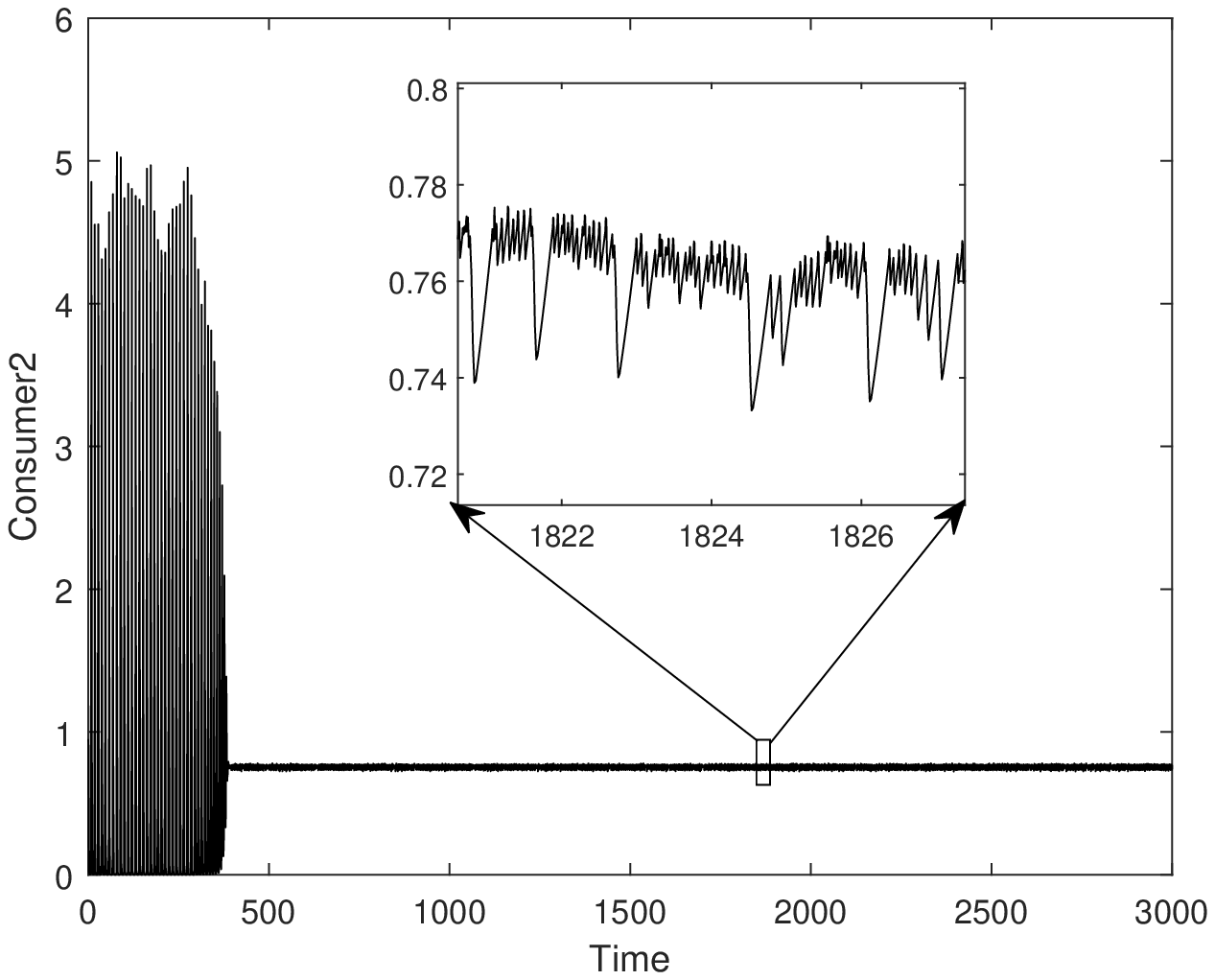}
		\caption{Consumer 2 success in region II}
	\end{subfigure}
	\caption{Time series of system \eqref{eq6} components success for the type IVb of zero$-$Hopf normal form catalog, $a<-1,b=-1$. a,b, and c show the time series plots in region$-$I of the typical dynamics pattern of type I, and d, e, and f show the time series in region$-$II of the typical pattern of type IVb. Temporal dynamics in region$-$III is symmetric to region$-$I.}
	\label{fig20}
\end{figure}

The time series of the system~\eqref{eq6} undergoing type IVb of zero-Hopf bifurcation is shown in figure \ref{fig20}.

\subsubsection{Type IIa of zero-Hopf normal form catalog, $-1<a<0,b=+1$}

Variation of parameters in their valid ranges is not lead to this pattern of dynamics. We conclude that this qualitative behavior is not feasible in system~\eqref{eq6}.

\section{Conclusions and Discussions}

{
Based on~\cite{denno}, among all the herbivorous insects examined from various guilds, competition is most frequently detected; the exception is chewing folivores, who are least likely to be interspecific competitors.  
They suggested that most interspecific competitive interactions were noticeably asymmetric. This asymmetry could be due to the lack of tendency for intraspecific competition that could diminish whole interspecific interactions.
To the best of our knowledge, in literature, expressing success as a phase variable and its dynamics in competition is a new view in our study. 
But, numerous studies focused on either interspecific or intraspecific competition for both exploitative and interference represent correlates biologically that are instances of success either of predation or escaping from being preyed on.
For instance, the copepod Diaptomus tyrrelli has a potential competitor and predator, namely Epischura nevadensis. In the presence of Epischura nevadensis, the filtering rate of the copepod Diaptomus tyrrelli reduces up to 60\%; see~\cite{foltGoldman}. 
They claimed that the reduction in filtering rate is related to interference competition between two species and may have roots in the evolution of a mechanism for avoiding predation. 
The reduction mentioned above in the filtering rate is the same as the reduction in success in consumption in this paper. 
They concluded that although the presence of Epischura can cause a substantial decrease in the filtering rate of Diaptomus, the ultimate effect on the fitness of Diaptomus may not be harmful. Thus, filtering rate reduction could be a mechanism for escaping predation that may have evolved because of predation pressure on early instar. 
The filtering rates were measured for each species monthly with the unit of ml/day per animal. Patterns obtained in two-species experiments show a kind of fluctuations in~\cite[Fig 2-B, C]{foltGoldman}) similar to ones depicted nearly in every time series of success dynamics in our study. 
  }\\

In two-dimensional Lotka-Volterra competitive systems, attractors are usually stable equilibria or stable limit cycles; see \cite{wang10}. 
The nonlinear functional response could be a factor that increases the complexity of such a competitive system. The general Rosenzweig-MacArthur prey-predator model is sensitive to the mathematical form of the predator functional response; see \cite{gunoc18}. They showed that among Holing type II, Ivlev, and Trigonometric functions, only the last form can give rise to subcritical Hopf bifurcation and has a saddle-node bifurcation of periodic orbits. 
The results in~\cite{soufbaf} confirm the competitive exclusion of species with low consumption rates, which are the initial population size $\times$ per capita consumption rate, in a well$-$mixed habitat patch. 
Those results examined the number of patches affecting the coexistence and spatial distribution of the two parasitoids using a stochastic cellular automaton (CA) discretized of the LV-RM model \eqref{eq1} developing various patches' numbers on a two$-$dimensional habitat with a periodic boundary condition. The in-patch dynamics ensure the coexistence of two parasitoids regarding the number of habitat patches. Also, it was evident that there existed a minimum habitat size to ensure species persistence and coexistence. However, in~\cite{soufbaf}, the dynamical analysis is not performed. \\

The model of this paper is improved biologically by incorporating the type II functional response of the Holling function for both competing consumers. 
Here, we considered the RM model, which is integrated into a competitive LV model and has type II of Holling functional response to study the coexistence. 
Here, we studied two competitors acting on one resource in a three-dimensional LV-RM system and depicted some attractors with interesting dynamics after Hopf and zero-Hopf bifurcation.
{
  Here, we showed that another class of competition models, LV-RM, generates oscillations and chaos if two consumers compete for one limited biotic resource. Competition models link competing species' dynamics to the resources' dynamics. An attractive feature of resource competition models is that they use species' biological traits to predict the time course of competition. Hence, it becomes possible to pinpoint the mechanisms that underlie nonequilibrium dynamics in competition models \cite{huisman2001}. In this paper, we analyzed the biological mechanism of two-species competition exhibiting oscillations and chaos.}
Current results show that Smale's complicated dynamics~\cite{smale76} can occur in a three-dimensional LV-RM system, also in which the competitive Lotka-Volterra subsystem is of dimension 2. 
We focused on the new concept of success as a virtual source of competition among predators and prey. Since we have two predators and one prey, we introduce a three-dimensional model for the dynamics of success. The quantitative value of success can be coupled with the population dynamics of each species to observe the related population dynamics, but that is beyond the scope of this paper. The success of each competitor (two predators and one prey) is not a fixed scalar number. Instead, they deviate according to their intrinsic time scale. Frequent experience of lower value of success opens the gateway to extinction regardless of the current number of that species' population. The conventional competition model for two predators and one prey predicts that only one predator is extinct unless in the presence of specific symmetry. We observed that the asymmetric action of predators is the key for both sides to survive. In particular, we observed that the inferior competitor survives and coexists with the superior competitor. Based on the model, the inferior consumer waits until the superior one performs its functional response (during its handling time), then the inferior consumer starts to attack with improved rates.\\

In integrated form with the RM model, we founded new attractors with complex dynamics. Dimensionless competition coefficients $\alpha$ and $\beta$ determine the sign of the first Lyapunov coefficient, $a$, and $b$ in Hopf bifurcation and zero-Hopf bifurcation, respectively. The attack rate amplifies when the other competing consumer falls into the relaxation phase and other parameters are valued symmetrically. From a practical point of view, the stability of competitors' interactions is significant. For example, in some biological control programs, a predatory species enters a new environment beside a native predator whose coexistence is desirable~\cite{soufbaf}. 
{
  Lawton and Hassell,~\cite{lawtonHassell}, gathered data on competition in various insect communities with high variation of feeding strategies and suggested that irrespective of taxonomic group or type of resource, observed competition is highly asymmetrical known as amensalism (strongly asymmetrical competition through which one competitor has nearly zero effect on the other). Using zero-growth isoclines of a two-dimensional competition LV model, they showed that $\alpha_{i,j}$ needs to be positive while  $\alpha_{j,i}$ should be nearly equal to zero. By reviewing many examples of competition in insect communities, Lawton and Hassell (1981) showed that amensalism is nearly two times more general than symmetrical competition. In the current study, either before or after Hopf bifurcation, we observed a sharp asymmetrical competition; for example, see Hopf case 4 before bifurcation. However, except for Hopf case 1, this asymmetricity decreased 221 times after bifurcations. The most asymmetricity in Zero-Hop bifurcation appeared in case IVa. Here, we observed an hourglass-shaped attractor in the simulation. However, in type IIb of zero-Hopf bifurcation, like case 2 of Hopf, a kind of Rossler attractor was depicted; for both cases, the asymmetricity was of small values in comparison, especially before bifurcation. So, it seems that the dynamics of success in less asymmetrical competitive interaction shows screw chaos (namely the Rossler attractor). While, in high asymmetricity, say amensalism or dominance relations, we found more regular attractors like repeated limit cycles in case 4 of Hop.}
The underlying mathematical formalisms of this study are codimension-one Hopf bifurcations and codimension-two zero-Hopf bifurcation. 
In the first case of Hopf bifurcation, decreasing competitive effects lead to changing the stability patterns from node to R\"ossler$-$type attractor. In this case, reduction in $\alpha$ and $\beta$, increase in the initial values of the shared resource, and $\mathrm{FRR}$ leads to this kind of bifurcation. But, in the second case of Hopf bifurcation, increasing competitive effects change the stability pattern from stable node to the R\"ossler type attractor. In this case, an increase in $\alpha$, $\beta$, the initial values of the shared resource, and $\mathrm{FRR}$ leads to this kind of bifurcation. In the third type of Hopf bifurcation, an increase in $\alpha$, a decrease in $\beta$ and the initial values of the shared resource, and the same $\mathrm{FRR}$ leads to dynamics change from unstable limit cycles to spirals. But, in the fourth case of Hopf bifurcation, increased $\alpha$ and decreased $\beta$ lead to changing the stability pattern. Asymmetric attack rates and relaxation-oscillation mechanisms assure coexistence. The asymmetry is detectable in various regions of the zero-Hopf bifurcation diagram. The dimensionless competition coefficients $\alpha$ and $\beta$ show that the per-capita inhibitory effect of consumer species on each other is inversely related to the attack rate and the initial population of the competitors. This issue is the main idea of introducing $\mathrm{FRR}$. The notable result of the present study is the complicated dynamics such as R\"ossler$-$type, hourglass, and tornado attractors that are depicted for the first time from a realistic deterministic model of LV-RM class after Hopf and zero-Hopf bifurcations. Also, the relationship between competitive coefficients in the form of the Lyapunov first coefficient formula is a precursor in the construction of competitive functional response in future studies.\\

{
Interaction strength $\eta_{i,j}$ in the current work varied from $2.06$, in type IIb of zero-Hopf, to $1.01 \times 10^14$, in case4 of Hopf bifurcation. This variation raises a question about the boundaries of competition strengths by which one can categorize various competitive behaviors. 
Moving from simple competition towards ammensalism and higher asymmetry via dominance relations should be quantified by boundaries in the spectrum of intensity of competitive interactions. 
Another issue worth studying in future work is the range of interplay between regular competition and lethal interference in which one natural enemy kills another \cite{collier2001}. One may expect that the strength difference between these two cases should be profoundly high. These high variations in interaction strengths might be more logical in interactions between competitors of highly distant taxa. For instance, 
Brown and Davidson (1977), using desert rodents and ants, showed that strong-competitive interactions among distantly related organisms are probably widespread and common in natural ecosystems. 
They showed that rodents and ants are the primary granivores in desert ecosystems and that their competitive interaction results in a compensatory reduction of seed resources when one taxon is absent. Similarly, Hochberg and Lawton (1990), under interkingdom competition title and using a system of insect host-parasitoid-pathogen, suggested that competitive coexistence between parasitoids and pathogens on a shared limited insect host may occur as either constant or oscillatory dynamics. 
They suggested that when parasitoids search randomly, the dynamics of long-period cycles appear, while when the survival of pathogens from one generation to another decreases, short-period cycles occur. They showed that if one predator oscillates, the other acts either periodically or aperiodically, after which three-species equilibrium is impossible, i.e., nonequilibrium dynamics for host-parasitoids-pathogens in nature.
In contrast, some researchers claimed that interspecific competition is not a major driving force of evolution in many insect communities \cite{shoroks84}. 
They suggested that the most significant parameters in competitive models are not per capita competition coefficient $\alpha_{i,j}$, while the carrying capacity $k$, as the spatial parameter, plays the core role in the dynamics. 
They claimed that the superior competitor's removal from simulation does not impact the dynamics, and therefore competition could not be a driving force of organization in community structure. 
Bolnick (2001) showed that frequency-dependent competition is crucial in maintaining stable population polymorphisms. There, competition also has a role in favoring the initial evolution of the greater phenotypic variance, and in the rate of niche expansion. 
As an ecologist, Bolnick (2001) suggested that when interspecific competition reduces, the intraspecific one becomes the driving force of rapid diversification \cite{bolnik2001}.
How interference competition affects community structure depends on the configuration of the dominance interactions among the species, and many organisms of a wide range of taxa, including many insects (especially in guilds of omnivorous ants) are in linear dominance relations \cite{lebrun2005}. 
}\\

We employed nondimensionalization and presented the complex dynamics of the model based on values of dimensionless bifurcation parameters and dimensionless phase variables. Reverting the dimensionless parameters and phase variables to their actual dimensionalized values requires extensive numerical computations, which could be a part of further research. The nontrivial point of this reverting is the existence of a biologically viable set of survival.   
As another possible further study, we propose the strong coupling of success dynamics with population dynamics. At the level of success dynamics, the integrators of population dynamics filter the frequency spectrum of success relaxation oscillation. This coupling and filtration may lead to new time scales at the level of population dynamics. This complementary analysis is capable of revealing the complex dynamics of species competitions. The possible degeneracies have codimension two and higher such as Hopf-Hopf, Pitchfork-Hopf, and double${}_{}$zero-Hopf. One may also consider the concept of delay as a source of fluctuations. Delay may cause more complexity in combinations with higher codimension bifurcations. It is also interesting to focus on differential tools for manifolds of population dynamics and consider distribution theory. As another direction of further study, one may focus on the effect of the measure of initial perturbations on survival and extinction.\\

{\bf Competing Interests:}
Authors disclose that there are no financial or non-financial interests that are directly or indirectly related to the work submitted for this publication.

\end{document}